\documentclass[acmtog]{acmart}

\AtBeginDocument{%
  \providecommand\BibTeX{{%
    \normalfont B\kern-0.5em{\scshape i\kern-0.25em b}\kern-0.8em\TeX}}}

\setcopyright{acmcopyright}
\acmJournal{TOG}
\acmYear{2022} \acmVolume{1} \acmNumber{1} \acmArticle{1} \acmMonth{1} \acmPrice{15.00}\acmDOI{10.1145/3508372}

\usepackage[l2tabu,orthodox]{nag}

\usepackage{graphicx}
\usepackage[T1]{fontenc}
\usepackage[utf8]{inputenc}
\usepackage{mathtools} % this package already loads amsmath
\usepackage{textcomp}
\usepackage[binary-units]{siunitx} % more powerful than gensymb
\usepackage{physics} % partial derivatives & other stuff
\usepackage{stmaryrd}
\usepackage{mathrsfs}
\usepackage{ifthen}
\usepackage{float} % explicitly loaded before hyperref (http://www.shawnlankton.com/2009/04/hyperref-graphicx-algorithm/)
\usepackage{hyperref}
\usepackage{cleveref}
\usepackage{csquotes}
\usepackage{calc}
\usepackage{enumitem}
\usepackage{wrapfig}

\usepackage{array}
\usepackage{booktabs}
\usepackage{multirow}

\usepackage[percent]{overpic}

\definecolor{Green}{rgb}{0.4,1,0.4}
\definecolor{Violet}{rgb}{0.6, 0.35, 0.7}

\newcommand{\revisionn}[1]{{#1}}
\newcommand{\revision}[1]{{#1}}

\providecommand{\finalversion}{0}
\ifthenelse{\equal{\finalversion}{1}}
{
	\renewcommand{\DZ}[1]{}
	\renewcommand{\DP}[1]{}
	\renewcommand{\HS}[1]{}
	\renewcommand{\TS}[1]{}
	\renewcommand{\MB}[1]{}
	\renewcommand{\jeremie}[1]{}

	\renewcommand{\todo}[1]{}
	\renewcommand{\draft}[1]{}
	\renewcommand{\toref}[1]{}
	\renewcommand{\tocite}[1]{}
}
{}

\newcommand\restr[2]{{% we make the whole thing an ordinary symbol
  \left.\kern-\nulldelimiterspace % automatically resize the bar with \right
  #1 % the function
  \vphantom{\big|} % pretend it's a little taller at normal size
  \right|_{#2} % this is the delimiter
  }}

\hypersetup{
	unicode=true,
	colorlinks=true,
}

\newcommand{\RR}{\mathbb{R}}

\definecolor{diagcoldef}{rgb}{0.8,0.8,0.8}
\newcommand{\diagcol}[1]{{\color{diagcoldef}#1}}

\definecolor{goodcoldef}{rgb}{0.1529411765, 0.6823529412, 0.3764705882}
\newcommand{\goodcol}[1]{{\color{goodcoldef}#1}}

\newcommand{\interactiveplot}{\href{https://polyfem.github.io/tet-vs-hex/plot.html}{interactive plot}}
\newcommand{\code}{\url{https://github.com/polyfem/polyfem/}}
\newcommand{\data}{\url{https://archive.nyu.edu/handle/2451/44221}}
\newcommand{\scripts}{\url{https://github.com/polyfem/tet-vs-hex}}

\begin{document}

%%
%% The "title" command has an optional parameter,
%% allowing the author to define a "short title" to be used in page headers.
\title{A Large-Scale Comparison of Tetrahedral and Hexahedral Elements for Solving Elliptic PDEs with the Finite Element Method}

%
% The "author" command and its associated commands are used to define
% the authors and their affiliations.
% Of note is the shared affiliation of the first two authors, and the
% "authornote" and "authornotemark" commands
% used to denote shared contribution to the research.
\author{Teseo Schneider}
\email{teseo@uvic.ca}
\orcid{1234-5678-9012}
\affiliation{%
  \institution{University of Victoria}
  \country{Canada}
}

\author{Yixin Hu}
\affiliation{%
  \institution{New York University}
  \country{USA}
}

\author{Xifeng Gao}
\affiliation{%
  \institution{Lightspeed \& Quantum Studio, Tencent America}
  \country{USA}
}

\author{J\'er\'emie Dumas}
\affiliation{%
  \institution{Adobe Research}
  \country{USA}
}

\author{Denis Zorin}
\affiliation{%
  \institution{New York University}
  \country{USA}
}

\author{Daniele Panozzo}
\affiliation{%
  \institution{New York University}
  \country{USA}
}

%%
%% By default, the full list of authors will be used in the page
%% headers. Often, this list is too long, and will overlap
%% other information printed in the page headers. This command allows
%% the author to define a more concise list
%% of authors' names for this purpose.
% \renewcommand{\shortauthors}{Trovato and Tobin, et al.}

%%
%% The abstract is a short summary of the work to be presented in the
%% article.
\begin{abstract}
  The Finite Element Method (FEM) is widely used to solve discrete Partial Differential Equations (PDEs) in engineering and graphics applications. The popularity of FEM led to the development of a large family of variants, most of which require a tetrahedral or hexahedral mesh to construct the basis. While the theoretical properties of FEM basis (such as convergence rate, stability, etc.) are well understood under specific assumptions on the mesh quality,  their practical performance, influenced \emph{both} by the choice of the basis construction and quality of mesh generation, have not been systematically documented for large collections of automatically meshed 3D geometries. 

%In particular, most of the research in FEM assumed that high-quality meshes are available, a requirement that is often unreachable in industrial settings.
We introduce a set of benchmark problems involving most commonly solved \revisionn{elliptic} PDEs, starting from simple cases with an analytical solution, moving to commonly used test problem setups, and using manufactured solutions for thousands of real-world, automatically meshed geometries. For all these cases, we use state-of-the-art meshing tools to create both tetrahedral and hexahedral meshes, and compare the performance of different element types for common elliptic PDEs.

The goal of his benchmark is to enable comparison of complete FEM pipelines, from mesh generation to algebraic solver, and exploration of relative impact of different factors on the overall system performance. 

As a specific application of our geometry and benchmark dataset,
we explore the question of relative advantages of unstructured (triangular/tetrahedral) and structured (quadrilateral/hexahedral) discretizations.  We observe that for Lagrange-type elements, while linear tetrahedral elements perform poorly, quadratic tetrahedral elements perform equally well or outperform hexahedral elements for our set of problems and currently available mesh generation algorithms. This observation suggests that for common problems in  structural analysis, thermal analysis, and low Reynolds number flows, high-quality results can be obtained with unstructured tetrahedral meshes, which can be created robustly and automatically.

We release the description of the benchmark problems, meshes, and reference implementation of our testing infrastructure to enable statistically significant comparisons between different FE methods, which we hope will be helpful in the development of new meshing and FEA techniques.

\end{abstract}

%%
%% The code below is generated by the tool at http://dl.acm.org/ccs.cfm.
%% Please copy and paste the code instead of the example below.
%%
\begin{CCSXML}
<ccs2012>
   <concept>
       <concept_id>10010147.10010371.10010396.10010397</concept_id>
       <concept_desc>Computing methodologies~Mesh models</concept_desc>
       <concept_significance>500</concept_significance>
       </concept>
   <concept>
       <concept_id>10010147.10010371.10010396.10010398</concept_id>
       <concept_desc>Computing methodologies~Mesh geometry models</concept_desc>
       <concept_significance>500</concept_significance>
       </concept>
   <concept>
       <concept_id>10010147.10010371.10010396.10010401</concept_id>
       <concept_desc>Computing methodologies~Volumetric models</concept_desc>
       <concept_significance>500</concept_significance>
       </concept>
   <concept>
       <concept_id>10010147.10010341.10010370</concept_id>
       <concept_desc>Computing methodologies~Simulation evaluation</concept_desc>
       <concept_significance>500</concept_significance>
       </concept>
 </ccs2012>
\end{CCSXML}

\ccsdesc[500]{Computing methodologies~Mesh models}
\ccsdesc[500]{Computing methodologies~Mesh geometry models}
\ccsdesc[500]{Computing methodologies~Volumetric models}
\ccsdesc[500]{Computing methodologies~Simulation evaluation}

%%
%% Keywords. The author(s) should pick words that accurately describe
%% the work being presented. Separate the keywords with commas.
\keywords{Tetrahedral Mesh, Hexahedral Mesh, Finite Element, Comparison}

%%
%% This command processes the author and affiliation and title
%% information and builds the first part of the formatted document.
\maketitle

\section{Introduction}

% \DZ{I think the black-box part is a distraction; and will confuse the reader, as we are not really
%   addressing meshing in this paper; I would suggest reducing this greatly, not doing this now as
%   it requires further discussion; the main consequence of the black-box view of things is that we
%   measure accuracy vs time, and even this is not quite true -- meshing time is not included}

The finite element method (FEM) is commonly used to discretize partial differential equations (PDEs), due to its generality, rich selection of elements adapted to specific problem types, and wide availability of commercial implementations. At a high level, a FE analysis code takes as input the domain boundary, the boundary conditions, and the governing equations of the phenomena of interest, and computes the solution everywhere in the domain.

As an initial step in this procedure, the domain typically has to be discretized in a finite collection of elements. Many choices are possible, ranging from unstructured grids of tetrahedra to perfectly regular grids of cubes. Despite the large amount of research on mesh generation, we were unable to find a systematic study answering a basic question: ``What are the practical pros and cons of using unstructured (triangular/tetrahedral) or  structured (quadrilateral/hexahedral/grids) discretizations for commonly used elliptic PDEs?''.

%\DZ{the text below sounded as if preference for hexes is a uniformly held belief, while in reality lots of people simply use tets and quite happy with these}
This question is critical to inform the development of meshing algorithms: while tetrahedral meshes are easier to generate automatically, hexahedral meshes (i.e., meshes that are composed of only deformed cubes) are much more difficult to adapt to objects with complex geometries, while maintaining high mesh quality. One of the arguments motivating development of these more complex algorithms is a common belief that hexahedral elements  yield better accuracy for a given computational cost (see the introduction of, e.g., \cite{Lyon:2016,Guo:2020:PolyCube,Bernard:2016}).

The overall aim of our work is to provide an extensive benchmark for comparing the performance of FE pipelines, including \revisionn{automatic} meshing, FE basis construction, and algebraic system solvers, on a set of most common \revisionn{elliptic} PDEs and a set of realistic geometries. As an immediate application, we explore the performance of widely used families of elements, coupled with standard solvers, on a large set of meshes generated using currently available meshing algorithms.

More specifically, we compare the \emph{efficiency} of different elements, that is, how much time is typically required to obtain a solution with a given accuracy for different element types on automatically generated unstructured meshes, on manually and automatically generated semi-structured meshes, and on regular lattices.
%We can thus measure the pros and cons of using manually designed structured, semi-structured, or automatically generated unstructured meshes for the solution of common elliptic PDEs using the finite element method.

We consider standard Lagrangian bases~\cite{szabo1991finite,ciarlet2002finite} of varying degrees, as well as  serendipity~\cite{Zienkiewicz:2005:TFE} elements (for hexahedra only), which are by far the most popular brick element. Finally, we perform several comparisons using spline-based elements~\cite{Hughes:2005:IAC}, which have recently gained popularity in the IsoGeometric Analysis (IGA) community. While this clearly does not reflect the broad range of existing element types and PDEs in the literature, it includes the most popular general-purpose elements currently used in commercial and open-source FE systems.  \revisionn{For solving the resulting linear systems, we consider both state of the art direct~\cite{code:pardiso-a} and iterative solvers~\cite{code:hypre}}.

We collected a set of test problems of varying complexity for elliptic PDEs (including Poisson, linear elasticity, Neo-Hookean elasticity, incompressible elasticity, and incompressible Stokes equations). Our set includes common simple test problems \revisionn{(where most of the hex-meshes are grids)}: beam bending, beam twisting, driven cavity flow, planar domain with a hole, elasticity problems with singular solutions, as well as a large-scale benchmark of manufactured solutions~\cite{SALARI:2000:CVB} on $3\,200$ \revisionn{automatically meshed}, real-world, complex 3D models. Our model collection includes both CAD models and scanned geometries, providing a realistic sampling of analysis scenarios. \revisionn{We use TetWild~\cite{Hu:2018:TMI,Hu2019fast} and MeshGems~\cite{code:meshgems} to generate the tetrahedral and hexahedral meshes respectively (we also included the state-of-the-art meshes from Hexalab~\cite{Bracci:2019:HNA}).}

This combination of test models, 3D meshes for these models, elements and solvers is representative of many common FE application scenarios.

We quantify (to the best of our knowledge, for the first time) the overall performance differences between these two families of elements. Our main conclusion is that, while linear elements on triangular/tetrahedral meshes exhibit well-known problems, quadratic tetrahedral elements  perform similarly or better (i.e., require similar or less time to compute a solution with a given accuracy) than Lagrangian elements on semi-structured hexahedral meshes, and are somewhat inferior (but still competitive, especially considering tetrahedral meshing is much faster and more robust) to the performance of spline elements on regular lattices when a direct solver is used. Combined with available state-of-the-art robust meshing techniques, quadratic tetrahedral elements are a good choice to realize a fully automatic pipeline,.e.g., for SciML applications, or shape optimization, without sacrificing performance compared to hexahedral elements, which require far more complex and less robust mesh generation.
More detailed conclusions are presented in Section~\ref{sec:conclusions}.

We emphasize that our study is limited to a specific set of  PDEs, commonly used geometry-agnostic linear solvers, \revisionn{and state-of-the-art meshing algorithms}; we leave adding dynamic scenarios, multi-physics, \revisionn{different linear solvers}, and other extensions as future work -- the provided framework can be readily extended to these cases.  We also note that adaptive refinement is simpler for hexahedral meshes and, as a consequence, adaptive geometric multigrid solvers are more readily available \cite{code:dealii}, although it is possible to develop similar solvers for tetrahedral meshes \cite{kohl2019hyteg}.  While the outcome of our study should not be interpreted as a reason to favor tetrahedral discretizations in all situations (and there are applications of hexahedral meshes outside the scope of FEM discretizations, such as lattice structure design), it does point to the need for direct experimental evaluation of meshing strategies, in the context of specific target applications.

We provide the complete source code\footnote{\code{}} for the integrated analysis pipelines we tested, the dataset we used\footnote{\data{}},
the benchmark solutions, and the scripts to reproduce all results\footnote{\scripts{}}, to enable researchers and practitioners to easily expand this study to additional mesh types (such as polyhedral meshes) and bases.

This study is divided into five sections: we first introduce the closest related works on meshing and analysis (Section \ref{sec:related}). We then overview the background on mesh types, basis, and the model PDE that we consider in this study (Section \ref{sec:background}). We divide the experimental evaluation into a set of individual experiments, targeting a set of common test problems (including problems with singularities) in Section \ref{sec:classical-problems}, and then perform a large scale analysis on thousands of automatically generated meshes in Section \ref{sec:large-dataset}. We finally draw  conclusions and identify open challenges in Section \ref{sec:conclusions}.

%!TEX root = ./00-main.tex

\section{Related Work}
\label{sec:related}
We first review existing comparisons of different types of finite elements (Section~\ref{sec:related_others}), then briefly discuss commonly used finite element software (Section~\ref{sec:related_fea}) and the state-of-the-art  meshing algorithms (Section~\ref{sec:related_meshing}).

\subsection{FEA on Unstructured and Structured Meshes}
\label{sec:related_others}

To the best of our knowledge, our study is the first large-scale comparison between different commonly used types of elements in FEM.  However, there are multiple existing comparisons focused on specific models and physics.

In \cite{Cifuentes:1992:APS}, the authors conclude that quadratic tetrahedral meshes lead to roughly the same accuracy and time as linear hexahedral meshes, by comparing solutions for several simple structural problems. By evaluating the eigenvalues of the stiffness matrices of various nonlinear and elastoplastic problems,~\cite{Benzley:1995:ACO} reports that, in their study, linear hexahedral meshes are superior to linear tetrahedral meshes. The authors also show that linear hexahedral meshes are slightly superior to quadratic tetrahedral meshes in the nonlinear elastoplastic analysis experiment.

A more recent work,~\cite{Tadepalli:2010:ACO,Tadepalli:2011:COH}, focuses on modeling footwear with a nonlinear incompressible material model under shear force loading conditions. The conclusion of these works is that trilinear hexahedral meshes are superior to linear tetrahedral meshes, and that quadratic tetrahedral elements are computationally more expensive compared to trilinear hexahedral elements, but have higher accuracy. \cite{Wang:2004:BTE} compares tetrahedral and hexahedral meshes on linear static problems, modal and nonlinear analysis. The study concludes that quadratic tetrahedral and hexahedral elements have similar performance, but quadratic hexahedra are computationally more expensive. The same study also confirms that linear tetrahedra are too stiff for large deformations, and linear hexahedra with large corner angles should be avoided in regions of stress concentration. \revision{The study is restricted to a small set of geometries and focuses on manual hexahedral mesh generation. Our study instead focuses on automatic meshing algorithms for both tetrahedral and hexahedral meshes, and we provide experimental results on thousand of complex geometric models and a wide array of elliptic PDEs.}

In medical applications, results for femur models~\cite{Ramos:2006:TVH} show that linear tetrahedral meshes of the simplified femur model lead to a closer agreement with the theoretical ones, while quadratic hexahedral meshes are more stable and the result is less affected by mesh refinement. On a kidney model,~\cite{Bourdin:2007:COT} observes that both linear and quadratic tetrahedral meshes are slightly stiffer than hexahedral meshes, but are more stable when high impact energies are present in the simulation. For heart mechanics and electrophysiology,~\cite{Oliveira:2016:COT} notes that quadratic hexahedra are slightly better than quadratic tetrahedra in the mechanics regime, while linear tetrahedral meshes are the best choice for the electrophysiology problem.

\subsection{Finite Element Analysis Software}
\label{sec:related_fea}
There exists a large number of libraries and software for finite-element analysis, both open-source and commercial. An exhaustive comparison of all existing packages\footnote{\footnotesize A non-exhaustive list of open-source FEA packages known to the authors include, in alphabetical order, code\_aster~\cite{code:codeaster}, Deal.II~\cite{code:dealii}, DOLFIN (FEniCS)~\cite{code:fenics}, ElmerFEM~\cite{code:elmerfem},  FEATool Multiphysics (MATLAB)~\cite{code:featool}, Feel++~\cite{code:feelpp}, FEI (Trilinos)~\cite{code:trilinos}, Firedrake~\cite{code:firedrake}, FreeFEM++~\cite{code:freefempp}, GetDP~\cite{code:getdp}, GetFEM++~\cite{code:getfempp}, libMESH~\cite{code:libmesh}, MFEM~\cite{code:mfem}, Nektar++~\cite{code:nektarpp}, NGSolve~\cite{code:ngsolve}, OOFEM~\cite{code:oofem}, PolyFEM~\cite{code:polyfem}, Range~\cite{code:range}, SOFA~\cite{code:sofa}, and VegaFEM~\cite{code:vegafem}.} is beyond the scope of this paper, therefore we discuss only several representative packages. We point out an interesting project~\cite{web:feacompare} attempting to maintain a complete list of FEA packages with a list of characteristics.

Our goal is to investigate and compare the performance of FEM on meshes with tetrahedral and hexahedral elements, using the standard Lagrangian basis functions and serendipity elements commonly used in engineering applications, as well as spline elements used in IGA.

Open-source packages such as FEniCS~\cite{code:fenics}, GetFEM++~\cite{code:getfempp}, libMesh~\cite{code:libmesh}, and MFEM~\cite{code:mfem} support both tetrahedral and hexahedral meshes, although very few (e.g., libMesh) implement both the 20-(serendipity) and 27-nodes variant for quadratic hexahedral elements. Deal.II~\cite{code:dealii} is another popular open-source FEA library, however it only supports quadrilateral and hexahedral elements.
Commercials packages such as ANSYS~\cite{code:ansys}, Abaqus~\cite{code:abaqus}, COMSOL Multiphysics~\cite{code:comsol} support Lagrangian tetrahedral elements, but surprisingly often implement only serendipity elements for hexahedra~\cite[Chapter~6]{Zienkiewicz:2005:TFE}. Given their popularity, we included serendipity elements in our study in addition to traditional Lagrangian elements.

Another increasingly popular choice of bases for hexahedral meshes are B-splines and NURBS, most commonly used in the context of isogeometric analysis (IGA)~\cite{Hughes:2005:IAC}. The popularity of spline bases stems from the fact that they have only one dof per element independently of the degree (however, the support of each basis function grows accordingly, and, as a consequence the stiffness matrices become less sparse). Defining this type of element on fully general hexahedral domains is an open problem~\cite{Aigner:2009:SVP, Martin:2010:VPO,Li:2013:SMT}. Due to their rising popularity, we deem important to include experiments with these elements in our study, but restrict them to cases where a regular lattice mesh is used.

Since none of these libraries implements both Lagrangian (tetrahedral and hexahedral), serendipity, and spline basis functions (hexahedral only) in the same framework, we added all the elements and basis used in this study to our own open-source FEA library~\cite{code:polyfem} to ensure a fair comparison. PolyFEM~\cite{code:polyfem} supports all these element types and interfaces with Hypre~\cite{code:hypre} and PARDISO~\cite{code:pardiso-a,code:pardiso-b,code:pardiso-c} for the solver and Eigen~\cite{code:eigen} for linear algebra.

\subsection{Meshing}
\label{sec:related_meshing}
Three-dimensional mesh generation has been thoroughly studied in multiple communities~\cite{Shewchuk:2012:DMG,Carey:1997:CGG,Owen:1998:ASO,Tautges:2001:TGO}. For the sake of brevity, we restrict our review to the techniques generating pure tetrahedral or pure hexahedral meshes, which are the focus of our study, with an emphasis on methods implemented in readily available open-source or commercial libraries.

\paragraph{Tetrahedral Meshing}
The most efficient, popular, and well-studied family of algorithms tackles the generation of meshes satisfying the Delaunay condition~\cite{Chew:1993:GQM,Shewchuk:1996:triangle,Shewchuk:1998:TMG,Ruppert:1995:ADR,Shewchuk:2012:DMG,Sheehy:2012:NBO,Remacle:2017:ATL,Du:2003:TMG,Alliez:2005:VTM,Tournois:2009:IDR,MURPHY:2001:APP,CohenSteiner:2002:CDT,Chew:1987:CDT,Si:2005:MPL,Shewchuk:2002:CDT,Si:2014:ICA,Hang:2015:tetgen,Cheng:2008:APD,Boissonnat:2005:PGS,Jamin:2015:CAG,Dey:2008:DAD,Chen:2004:ODT}. These methods are robust if the input is a point cloud, but might fail if the boundary of a shape has to be preserved exactly~\cite{Hu:2018:TMI,Hu2019fast}.

To overcome these robustness limitations, alternative approaches are based on a background grid~\cite{Molino:2003:TMG,Bronson:2013:LCC,Labelle:2007:ISF,Doran:2013:ISI,code:quartet}. The idea is to fill the bounding box of the 3D input surface with either a uniform grid or an adaptive octree, whose convex cells are trivial to tetrahedralize. These methods achieve high quality in the interior of the mesh (where the grid is regular), but introduce badly shaped elements near the boundary, which is often the region of interest in many practical simulations. On the other hand, front-advancing methods \cite{Cuilliere:2013:ADM,Alauzet:2014:ACA,Haimes:2014:MMO} start by marching from the boundary to the interior, adding one element at a time, pushing the problematic elements into the interior where the advancing fronts meet.

All these methods are unable to handle commonly occurring input surfaces which contain degenerated faces, gaps, and self-intersections. These types of defects are, unfortunately, common in CAD models, due to the NURBS representation (with a fixed degree) not being closed under boolean operations. To the best of our knowledge, the only method that was demonstrated to be capable of handling these cases \emph{robustly} is TetWild~\cite{Hu:2018:TMI}. It is based on a hybrid numerical representation to ensure correctness, and it allows a small, controlled deviation from the input surface to achieve a good element quality. We used this technique to generate all unstructured tetrahedral meshes in this study.

\paragraph{Hexahedral Meshing} aims at filling the volume enclosed by an input surface with hexahedra. Hexahedra also need to have a good shape to ensure good solution approximation.  The natural tensor-product structure of a hexahedron enables to define tensor-product bases, and, e.g., use spline-based elements,  but dramatically increases the complexity of meshing algorithms. Semi-manual or interactive approaches are usually employed, such as sweeping and advancing front methods \cite{Shepherd:2008:HMG,Gao:2016:SVD,Livesu:2016:SDA}, which are used in commercial software such as \cite{code:ansys,code:cubit}.

By allowing lower element quality, one can design automatic approaches based on regular lattices~\cite{Schneiders:1995:AGO,Schneiders:1996:AGB,Su:2004:AHM,Zhang:2006:AAQ,Zhang:2007:AGO} or on octrees~\cite{Schneiders:1996:OBG,Zhang:2006:AAQ,Ito:2009:OBR,Marechal:2009:AIO,Qian:2010:SFP,Ebeida:2011:ICR,Zhang:2013:ARR,Elsheikh:2014:ACO,Owen:2017:ATB,code:meshgems}.

Polycube methods~\cite{Gregson:2011:AHM,Livesu:2013:PMG,Huang:2014:BCO,Fang:2016:AHM,Fu:2016:EVP,Li:2013:SMT,Zhao:2018:REP} and field-aligned parameterization-based methods~\cite{Nieser:2011:CPO,Huang:2011:BAS,Li:2012:AHM,Jiang:2014:FFS,Solomon:2017:BEO,Liu:2018:SCO} aim at producing hexahedral meshes with as few irregular edges and vertices as possible, but designing robust algorithms of this type is still an open problem. Sample results from some of the previous methods have been recently collected into a single repository~\cite{Bracci:2019:HNA}, which we use in our study. We also generate a new dataset composed of 3200 hexahedral meshes using the commercial MeshGems-Hexa software \cite{code:meshgems}.

\section{Background}
\label{sec:background}

\subsection{FEM bases}
There is a multitude of different definitions of bases for both tetrahedral (or triangular) and hexahedral (or quadrilateral) element shapes, with different elements tailored to specific types of problems (e.g.,  axisymmetric elements, shell elements, plasticity elements, etc.). In our comparison, we target the most common choices: we use the standard  linear and quadratic Lagrange bases for tetrahedra,
\revision{which we denote $P_1$ and $P_2$ respectively,}
and hexahedra, \revision{with $Q_1$ denoting linear tensor-product basis and $Q_2$ quadratic tensor-product basis} \cite{szabo1991finite, Ciarlet:2002:fem}. We also use the serendipity basis \cite{Zienkiewicz:2005:TFE}, commonly used in commercial software, and spline basis \cite{Hughes:2005:IAC} for hexahedral elements. We use the standard Galerkin formulation~\cite{szabo1991finite, Ciarlet:2002:fem} with Gaussian quadrature for all our experiments, avoiding non-standard quadrature.

\subsection{Mesh and solution characterization}

We  use the number of vertices as a measure of the resolution of tetrahedral and hexahedral meshes, as the number of vertices is often used by the meshing algorithms as the ``budget'' that the meshing algorithms can use to create the best possible mesh, and the number of vertices is equal to the number of degrees of freedom in the case of linear (or tri-linear) elements.

In addition to this particular choice, we also investigate other metrics for a specific example (Table~\ref{tab:p-vs-q-fixed}), and provide an \interactiveplot{} that allows one to compare our results using 24 different measures:
solution error measured using $H^1$, $H^1$ semi-norm, $L_2$, $L^\infty$, $L^\infty$ of gradient, and $L^8$ norms;
mesh average edge length,
minimum edge length and
number of vertices;
the system matrix size and the
number of non zero entries,
the numbers of basis functions, dofs, elements,
and pressure basis functions;
timings for loading mesh data, building basis functions, computing the right-hand side, assembling the system matrix, solving the system, computing the errors, total time and
time  without right-hand side assembly.

\subsection{Model PDEs}

%Another largely represented topic in the literature is the definition of the physical model which is expressed in the form of a partial differential equation (PDE). For instance, just to model the reaction of an elastic body one can choose between more than seven different material models (and PDEs) depending on the regime (large/small deformations, elastic/plastic, etc.) or accuracy wanted.
We selected the following set of representative \textbf{elliptic} problems: (1) Poisson; (2) incompressible stationary Stokes fluid flow equations; (3) elasticity with linear Hooke's law as the constitutive equation; (4) Neo-Hookean elasticity (5) incompressible linear elasticity. We list the corresponding PDEs for completeness.

Let $\Omega\subset \RR^d$, $d\in\{2,3\}$ be the domain with boundary $\partial\Omega$.
We aim to solve
\begin{gather*}
    \mathcal{F}(x, u, \nabla u, D^2 u) = b,\text{ subject to}\\
    u  = d \text{ on } \partial\Omega_D\quad\text{and}\quad
    \nabla u\cdot n  = f \text{ on } \partial\Omega_N
\end{gather*}
for the function
\[
    u\colon \Omega \to \RR^n,
\]
where $D^2$ is the matrix of second derivatives, $b$ is the right-hand side, $\partial\Omega_D\subset\partial\Omega$ is the part of the boundary with Dirichlet boundary conditions, and $\partial\Omega_N\subset\partial\Omega$ is the part of the boundary with Neumann boundary conditions. Since we consider second-order PDEs only, $\partial\Omega_D\cap\partial\Omega_N = \emptyset$. The form of $\mathcal{F}$ and the role of $u$ depends on the specific PDE.

We consider polygonal and polyhedral domains $\partial \Omega$ (possibly non-convex). The right-hand side $b$ in our test examples is analytic, the boundary $d$ is continuous and piecewise-smooth, and $f$ is piecewise smooth (but possibly with finite-jump discontinuities); under these assumptions, the weak solutions of the equations we consider are (at least) continuous, but the solution derivatives may be singular. We primarily focus on the error in the solution itself, rather than the derivative error, although consider the stress for some elasticity examples.  We state the model problems in the strong form, but only the weak solutions exist for many of the test cases.

\paragraph{Poisson Equation} This problem is given by
\begin{equation}\label{eq:poisson}
    \begin{cases}
        -\Delta u = b        & \text{on } \Omega            \\
        u  = d               & \text{on } \partial\Omega_D  \\
        \nabla u\cdot n  = f & \text{on } \partial\Omega_N.
    \end{cases}
\end{equation}

\paragraph{Incompressible \revisionn{Steady} Stokes Equations} The Stokes equations provide the relationship between the velocity $u$ and the pressure $p$ for an incompressible fluid with viscosity $\mu$.
\begin{equation} \label{eq:stokes}
    \begin{cases}
        -\mu\Delta u + \nabla p = b                 & \text{on } \Omega           \\
        -\div{u} = 0                                & \text{on } \Omega           \\
        u  = d                                      & \text{on } \partial\Omega_D \\
        \mu(\nabla u + \nabla^T u)\cdot n - pn  = f & \text{on } \partial\Omega_N
    \end{cases}
\end{equation}

\paragraph{Elasticity} Elasticity PDEs are formulated in terms of the stress tensor $\sigma[u]$ (which depends on the displacement $u$) as
\begin{equation}\label{eq:elasticity}
    \begin{cases}
        -\div{\sigma[u]} = b & \text{on } \Omega            \\
        u  = d               & \text{on } \partial\Omega_D  \\
        \sigma[u] \, n  = f  & \text{on } \partial\Omega_N.
    \end{cases}
\end{equation}
In this case the right-hand side $b$ plays the role of a body force, the Dirichlet boundary conditions are fixed displacement, and the Neumann ones are surface tractions.

Material models define how the stress $\sigma$ is related to the displacement field $u$. For the linear Hookean model,
\begin{equation}\label{eq:lin-elast}
    \sigma^L[u] = 2 \mu  \epsilon[u]+ \lambda \trace{\epsilon[u]} I
    \qquad
    \epsilon[u] = \frac 1 2 \left(\nabla u^T + \nabla u\right),
\end{equation}
where $\epsilon[u]$ is the strain tensor, $\lambda$ is the first Lam\'e parameter, and  $\mu$ is the shear modulus.
\revision{There are two common assumptions reducing the elasticity problem to a 2D problem, plane \emph{stress} and plane \emph{strain}; in our experiments we are using plane stress. In this case, the elasticity equation has the same form but with different constants~\cite{hughes:2012:finite}:}
\[
    \mu = \frac{E} {2 (1 + \nu)},\quad
    \lambda_{\mathrm{3D}} = \frac{E \nu}{(1 + \nu) (1 - 2 \nu)},\quad\mathrm{and}\quad
    \lambda_{\mathrm{2D}} = \frac{\nu  E} {1 - \nu^2}.
\]

Incompressible materials form a separate class: in 3D, an isotropic material has Poisson ratio equal to 0.5, and the previous equation is not well-defined, as $\lambda$ becomes infinite.
While isotropic materials in plane stress state cannot have this problem, as the isotropic Poisson ratio cannot exceed 0.5, anisotropic materials can have Poisson ratio 1 for in-plane deformations, and thus can be 2D-incompressible, which geometrically corresponds to the area of the cross-section of a material element preserved under deformations ~\cite{Lee:1997:aniso}). As a consequence, equations for 2D-incompressible materials in plane stress state are also of interest.
Additionally, when $\lambda$ grows, the linear system arising from the discretization of the PDE becomes unstable. A common way to avoid such problem is to introduce a Lagrange-multiplier-like function in the form of the pressure $p$. This leads to a mixed formulation of elasticity similar to Stokes equations which is stable for large $\lambda$s, and reduces to incompressible elasticity for $\lambda^{-1} \rightarrow 0$.
\begin{equation}\label{eq:incompressible-elast}
    \begin{cases}
        -\div(2\mu\epsilon[u] + p I) = b & \text{on } \Omega          \\
        \div u - \lambda^{-1}p = 0       & \text{on } \Omega          \\
        u  = d                           & \text{on }\partial\Omega_D \\
        \sigma^N[u] \cdot n  = f         & \text{on }\partial\Omega_N
    \end{cases}
\end{equation}

Finally, in the Neo-Hookean material model the stress is a nonlinear function of strain.
\begin{equation}\label{eq:nl-elast}
    \sigma[u] = \mu (F[u] - F[u]^{-T}) + \lambda \ln(\det F[u]) F[u]^{-T}
    \qquad
    F[u] = \nabla u + I,
\end{equation}
where $F[u]$ is the deformation gradient.

For elasticity problems, we often use the von Mises stresses
\begin{equation}\label{eq:von-mises}
    \begin{split}
        S_{\text{2D}}^2 =& \sigma_{0, 0}^2 -\sigma_{0, 0} \sigma_{1, 1} + \sigma_{1, 1}^2 + 3\sigma_{0, 1}\sigma_{1, 0}\\
        S_{\text{3D}}^2 =& %\frac{3}{2}\|\sigma\|_F,
        \frac{
            (\sigma_{0, 0} - \sigma_{1, 1})^2 +
            (\sigma_{2, 2} - \sigma_{1, 1})^2 +
            (\sigma_{2, 2} - \sigma_{0, 0})^2}{2} +\\
        &3(\sigma_{0, 1} \sigma_{1, 0}+\sigma_{2, 1} \sigma_{1, 2}+ \sigma_{2, 0} \sigma_{0, 2}).
    \end{split}
\end{equation}
\revision{Note that the stresses are discontinuous since they depend on the gradient of the displacement which is only $C^0$ for our discretizations. To mitigate visual artefacts we average the stresses around vertices in our plots. }

\subsection{Linear Solvers}
All FEM problems we consider require to solve a linear system, which, as the mesh size grows, dominates the running time. A vast amount of research has been invested in developing efficient and robust linear solvers. In our study we use two state-of-the art solvers: Pardiso~\cite{code:pardiso-a} a direct solver using the Cholesky factorization, which we use for smaller problems, and Hypre~\cite{code:hypre} an algebraic multigrid solver, which we use for larger problems. Direct solvers work particularly well in 2D, but scale poorly for 3D problems. We leave as future work a more detailed study on the effect of the linear solver on the solution time. The conclusions of this study hold for both types of solvers for our experimental setup and test problems.

\section{Common Test Problems}\label{sec:classical-problems}
% \todo{pick an 88 model do 3-4 refinements and plot ask Yixin for tet mesh}

We collected a number of standard test cases to cover different physical phenomena and different scenarios: fluid simulation (Section~\ref{sec:stokes}), linear elastic time dependent (Section~\ref{sec:time_dependent}), linear elastic bars (Section~\ref{sec:bened_bar}), linear orthotropic material models (Section~\ref{sec:ortho}), meshes with high aspect-ratio for linear elastic bars (Section~\ref{sec:high-aspect}), classical plane with hole with symmetric boundary conditions for compressible and nearly incompressible material (Section~\ref{sec:plate-hole}), nearly incompressible linear material (Section~\ref{sec:incompressible}), nonlinear Neo-Hookean material (Section~\ref{sec:torsion}), and nonlinear Neo-Hookean material with high stresses (Section~\ref{sec:high-stresses}).

% for fluid simulation (Section~\ref{sec:stokes}), linear Hookean material deformation (sections~\ref{sec:time_dependent}, \ref{sec:bened_bar}, \ref{sec:ortho}, \ref{sec:high-aspect}, and \ref{sec:plate-hole}), nearly incompressible linear material (Section~\ref{sec:incompressible}), and nonlinear Neo-Hookean material (Sections~\ref{sec:torsion} and~\ref{sec:high-stresses}).

Most of the solution domains are chosen to simplify manual creation of hexahedral meshes: the simulations will be performed on an unstructured tetrahedral mesh and a nearly regular lattice with the same number of vertices.
Experiments in Sections~\ref{sec:time_dependent} to~\ref{sec:incompressible} are run on a MacBook Pro 3.1GHz Intel Core i7, 16GB of RAM, and 8 threads. Experiments in Sections~\ref{sec:torsion} and~\ref{sec:high-stresses} are run on a cluster node with 2 Xeon E5-2690v4 2.6GHz CPUs and 250GB memory, each with max 128GB of reserved memory and 8 threads. For all experiments, we use the PolyFEM library \cite{code:polyfem}, which uses the Pardiso~\cite{code:pardiso-a,code:pardiso-b,code:pardiso-c} direct solver, and Newton iterations for the nonlinear problems.

Note that, for completeness, we also validated PolyFEM on the example in Figure~\ref{fig:square_beam} for linear and quadratic tetrahedra and serendipity hexahedra on Hooke material against Abaqus. The results are identical up to numerical precision.

\subsection{Incompressible Stokes}\label{sec:stokes}
\begin{figure}
    \centering\footnotesize
    \parbox{0.42\linewidth}{\centering
        \parbox{0.08\linewidth}{\rotatebox{90}{\centering $y$ veclocity}}
        \parbox{0.9\linewidth}{\includegraphics[width=\linewidth]{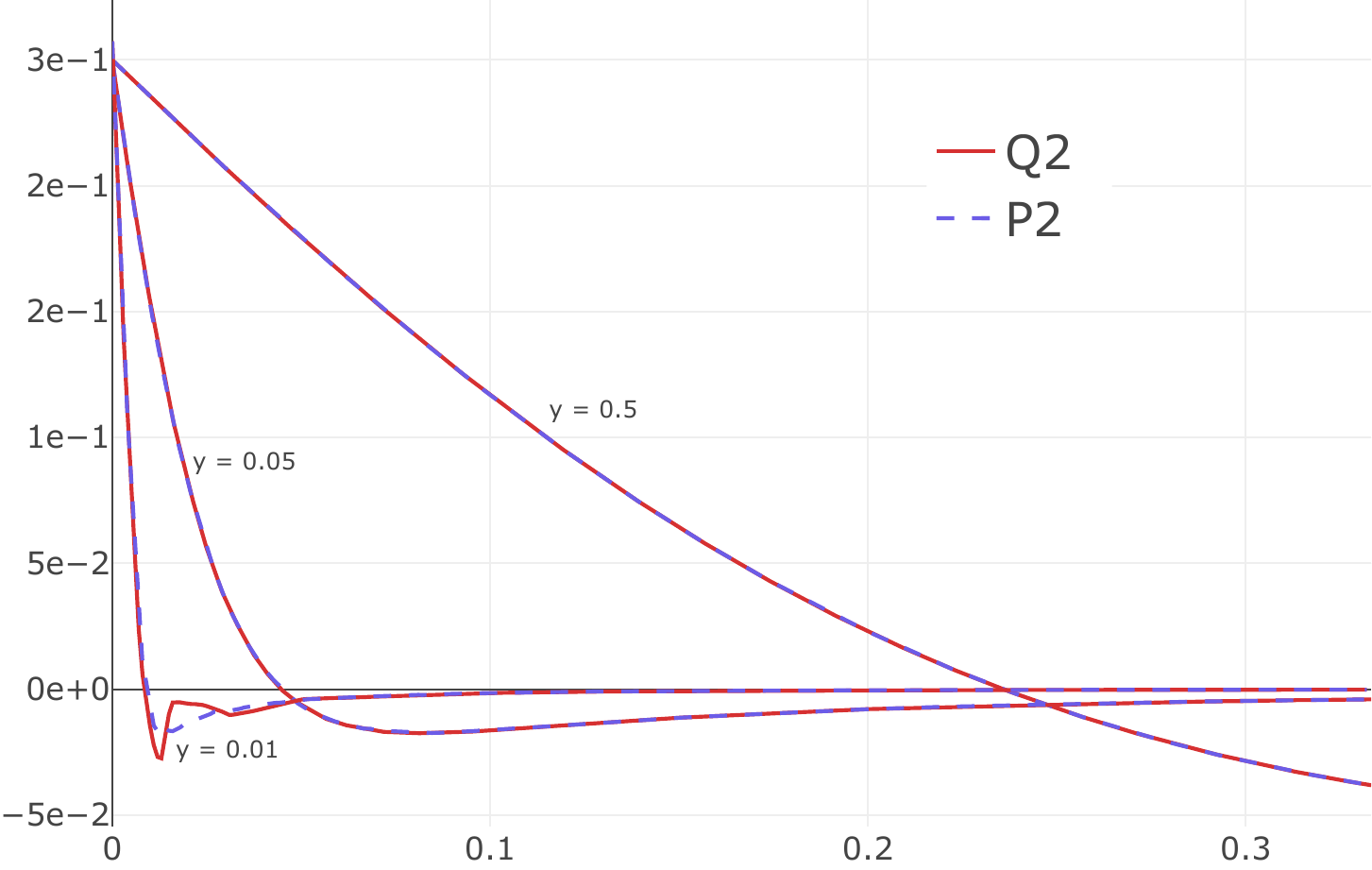}}\\
        \revisionn{$x$}
    }\hfill
    \parbox{0.54\linewidth}{\centering\footnotesize
        \includegraphics[width=0.4\linewidth]{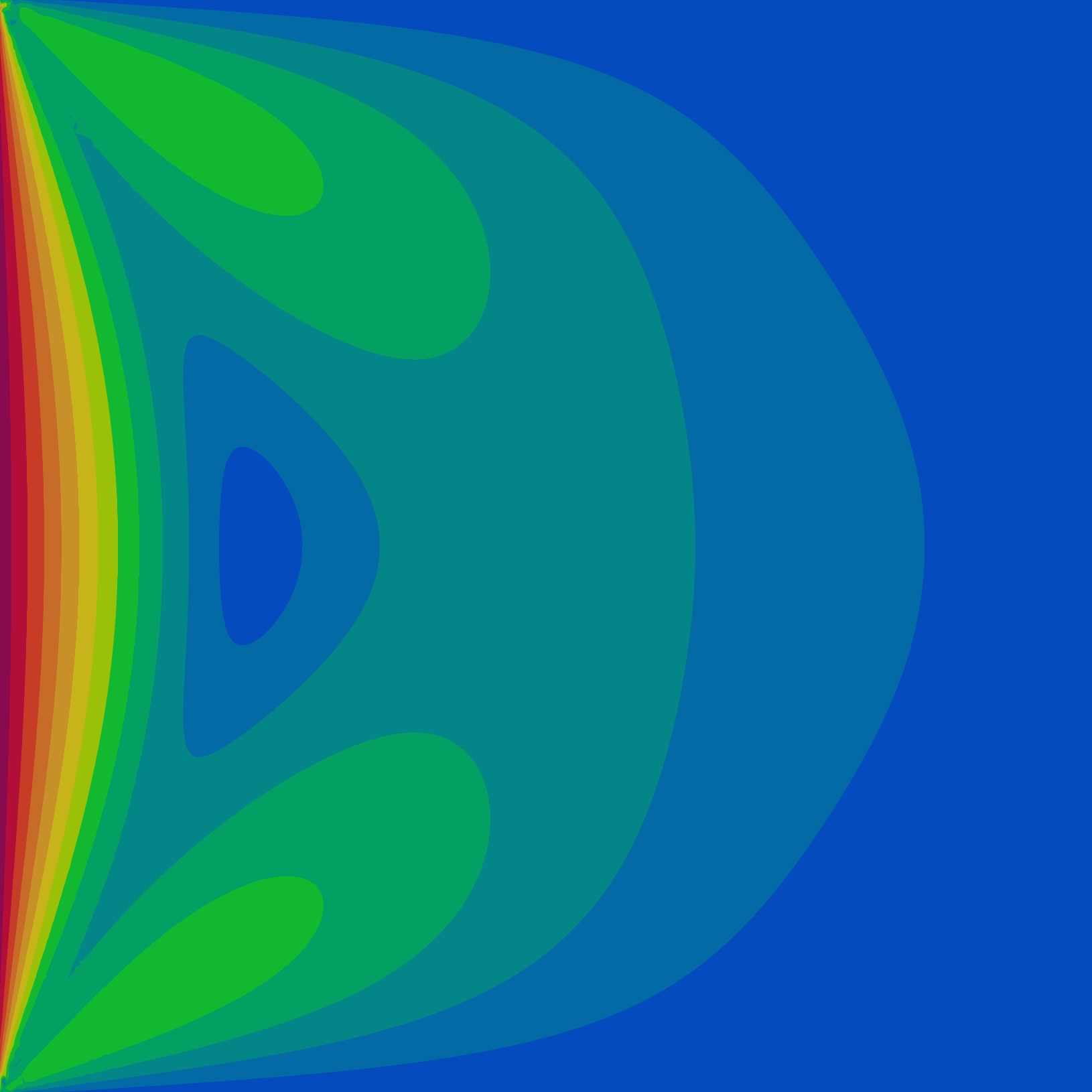}\hfill
        \includegraphics[width=0.4\linewidth]{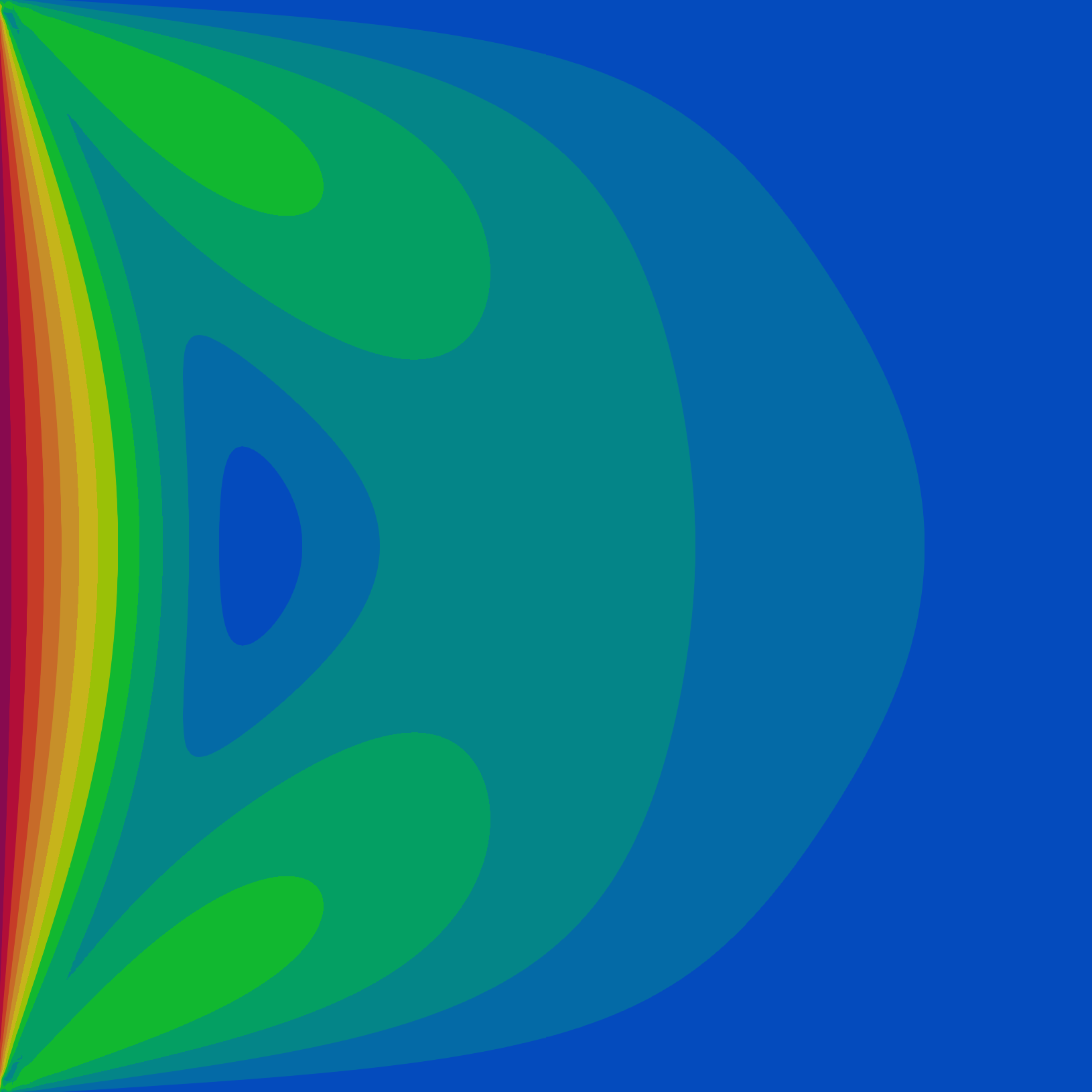}\hfill
        \includegraphics[width=0.1\linewidth]{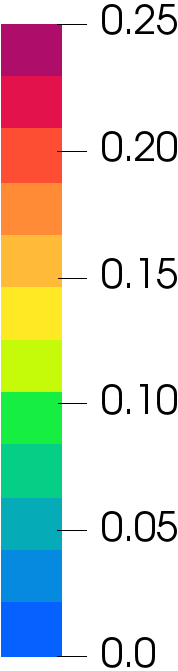}\par
        \parbox{0.4\linewidth}{\centering Tri. mesh}\hfill
        \parbox{0.4\linewidth}{\centering Quad. mesh}\hfill
        \parbox{0.1\linewidth}{~}\par
    }
    \caption{The velocity magnitude for a Stokes problem discretized with mixed elements. \revision{The plot shows the velocity in $y$-direction along horizontal lines $y=0.01$, $0.05$, and $0.5$ parametrized by \revisionn{$x$}.}}
    \label{fig:stokes}
\end{figure}

We use a planar square domain mesh with 4\,229 vertices for the triangle mesh and 4\,225 vertices for the regular grid. We simulate the Stokesian fluid~\eqref{eq:stokes} with viscosity $\mu=1$ in the standard ``driven cavity'' example: the fluid has zero boundary conditions on 3 of the 4 sides and a tangential velocity of 0.25 on the left side. Figure~\ref{fig:stokes} shows the results for mixed linear (for the pressure) and quadratic (for the velocity) elements: the results are indistinguishable between hexahedral and tetrahedral elements.

\subsection{Time-Dependent Linear Elasticity}\label{sec:time_dependent}

\begin{figure}\centering\footnotesize
    \parbox{.32\linewidth}{\centering
        \includegraphics[width=.45\linewidth]{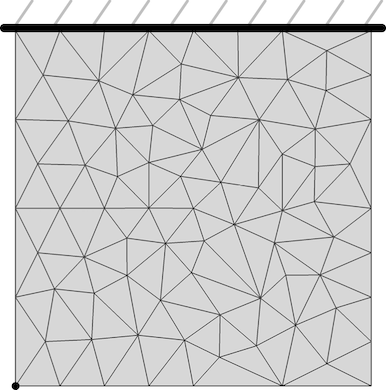}\hfill
        \includegraphics[width=.45\linewidth]{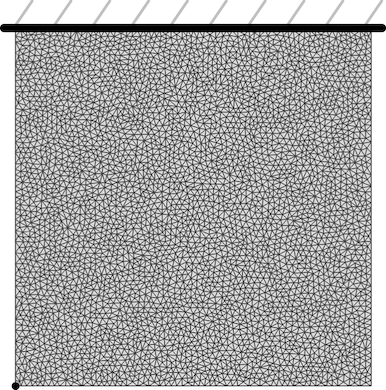}\par%\vspace{-0.7em}
        \parbox{.45\linewidth}{\centering \#v = 82}\hfill
        \parbox{.45\linewidth}{\centering \#v = 4\,229}\\[1em]
        \includegraphics[width=.45\linewidth]{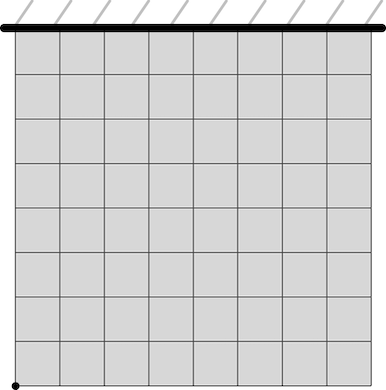}\hfill
        \includegraphics[width=.45\linewidth]{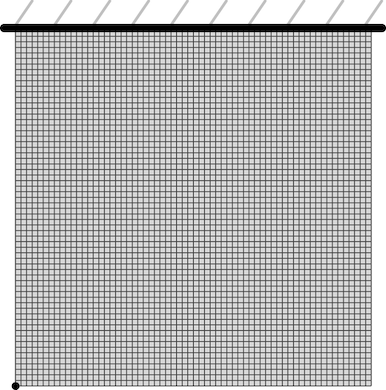}\par%\vspace{-0.7em}
        \parbox{.45\linewidth}{\centering \#v = 81}\hfill
        \parbox{.45\linewidth}{\centering \#v = 4\,225}\par
    }\hfill
    \parbox{.65\linewidth}{\centering
        \parbox{0.08\linewidth}{\centering\rotatebox{90}{\centering $x$ displacement}}
        \parbox{0.9\linewidth}{\includegraphics[width=\linewidth]{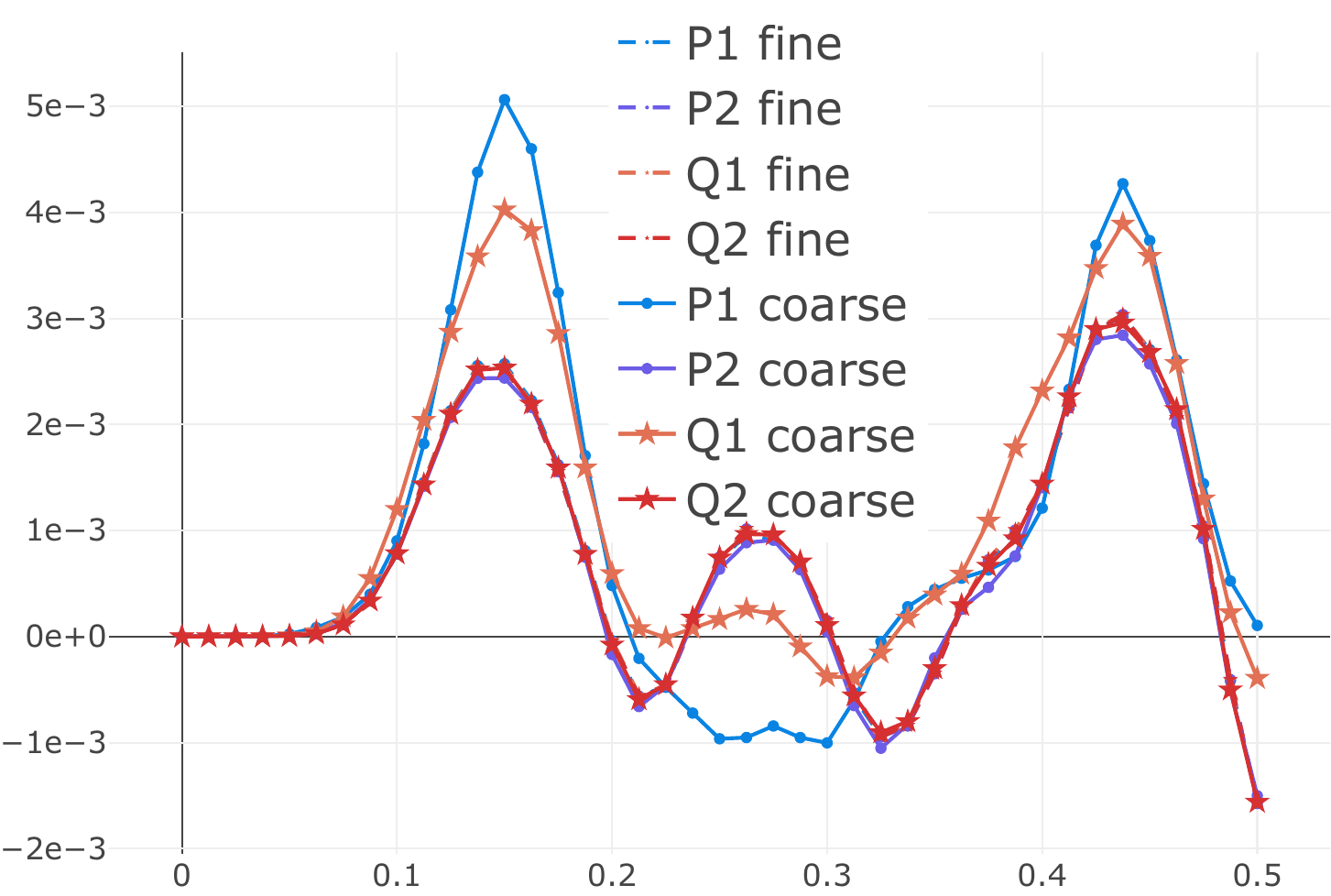}}\\
        Time
    }
    \caption{\revision{$x$ displacement of the bottom-left corner (black dot) of a unit square.}}
    \label{fig:gravity}
\end{figure}

\revision{We consider the dynamics of a suspended object under gravity: we fix the top part of a unit square with material parameters $E=200$ and $\nu=0.35$ and apply a constant body force of $20$ in the $y$ direction. We integrate the dynamic simulation for $t$ from $0$ to $0.5$ with 40 time steps integrated with Newmark~\cite{Newmark:1959:int}. We mesh the domain at a coarse and fine resolution, both for triangles and for quads. Figure~\ref{fig:gravity} shows the displacement in the $x$ direction of the bottom left corner for the 4 discretizations, using linear and quadratic elements.}

\subsection{Transversally Loaded Beam}\label{sec:bened_bar}

In this experiment, we consider beams with different cross-sections (square, circular, and I-like) in the $xy$-plane of length $L$. The beam is fixed (i.e., zero Dirichlet conditions are applied) at the end ($z=L$), and different tangential forces $f=[0, -f_y, 0]^T$, $f_y\in[-0.1, -2]$, are applied at $z=0$, opposite to the fixed side. The rest of the boundary is left free and we do not apply any body force. For these experiments we use linear isotropic material model~\eqref{eq:lin-elast} with Young's modulus $E=210\,000$ and Poisson's ratio $\nu = 0.3$. We study the displacement at the bottom corner of the moving end ($z=0$) in the $y$ direction and compare it with a dense solution to compute the error $e$ \revision{(note that the solution is singular only at $z=L$, far from the evaluation points)}. We report as $e_f$ the slope of the linear fit of the error as a function of the force magnitude. We also report the basis construction time $t_b$, assembly time $t_a$, solve time $t_s$, and total time $t$. Note that all the timings reported are averaged over $10$ different runs per force sample.

\paragraph{Square Cross-section}
\begin{figure}\centering\footnotesize
    \parbox{0.7\linewidth}{\centering\footnotesize
        \begin{tabular}{c|lll|l||c}
                  & $t_b$             & $t_a$             & $t_s$             & $t$               & $e_f$             \\
            \hline
            $P_1$ & 8.07e-3           & \goodcol{1.88e-2} & \goodcol{5.60e-2} & \goodcol{8.29e-2} & 6.14e-3           \\
            $P_2$ & 2.30e-2           & 1.80e-1           & 3.43e-1           & 5.47e-1           & 9.19e-5           \\
            $Q_1$ & \goodcol{5.96e-3} & 3.36e-2           & 6.39e-2           & 1.03e-1           & 1.27e-3           \\
            $Q_2$ & 1.46e-2           & 4.61e-1           & 4.34e-1           & 9.10e-1           & \goodcol{4.66e-5} \\
        \end{tabular}}\hfill
    \parbox{.28\linewidth}{\centering\footnotesize
        \includegraphics[width=.48\linewidth]{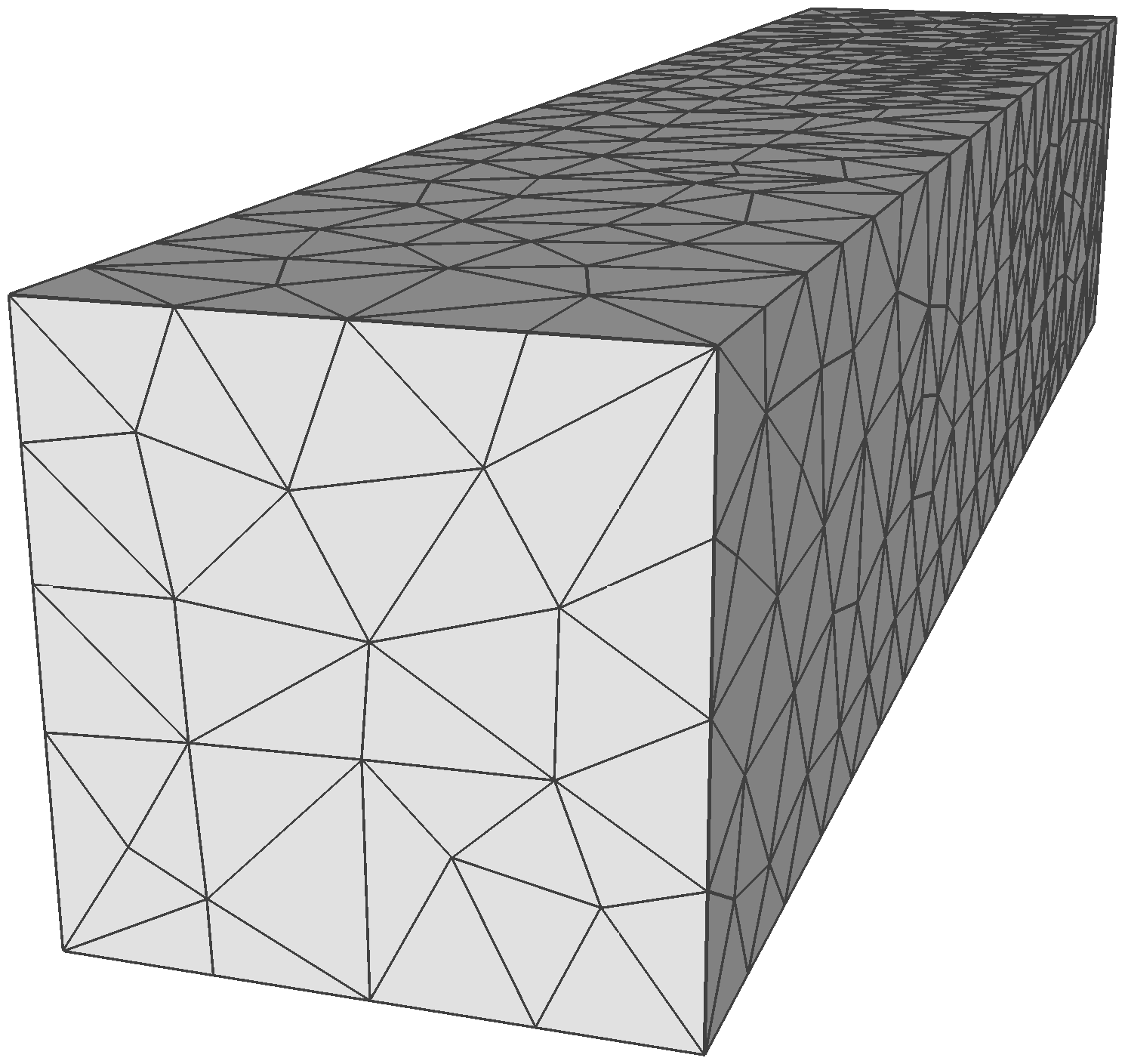}\hfill
        \includegraphics[width=.48\linewidth]{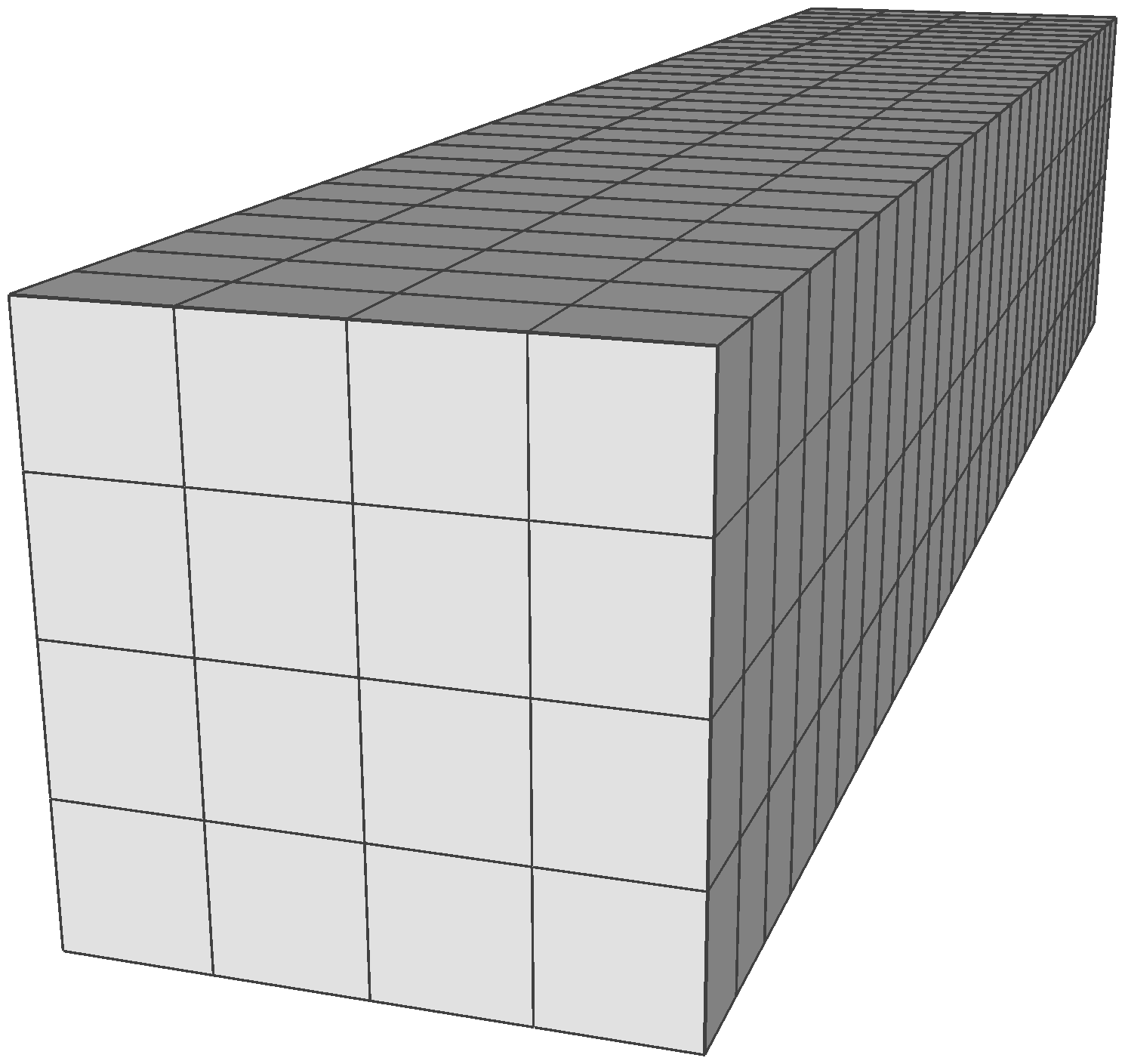}\par
        \parbox{0.48\linewidth}{\centering Tet. mesh}\hfill
        \parbox{0.48\linewidth}{\centering Hex. mesh}\par
    }
    \caption{Displacement error in the $y$ displacement of the moving endpoint compared with a dense solution for a unit force applied at the endpoint of a beam with a square cross-section. The times are averaged over $10$ runs per force sample.}
    \label{fig:square_beam}
\end{figure}
For running the simulation, we use a square cross-section of side $s=20$, length $L=100$ and mesh it with a  tetrahedral mesh with $739$ vertices and a hexahedral mesh (regular grid) with $750$ vertices. Figure~\ref{fig:square_beam} shows the errors compared with the dense solution, where trilinear hexahedral elements outperform linear tetrahedral elements but the quadratic counterparts are indistinguishable. Timing-wise, the quadratic tetrahedra are \revisionn{slightly} better.

\begin{figure}\centering\footnotesize
    \rotatebox{90}{\centering Time}
    \parbox{.48\linewidth}{\includegraphics[width=\linewidth]{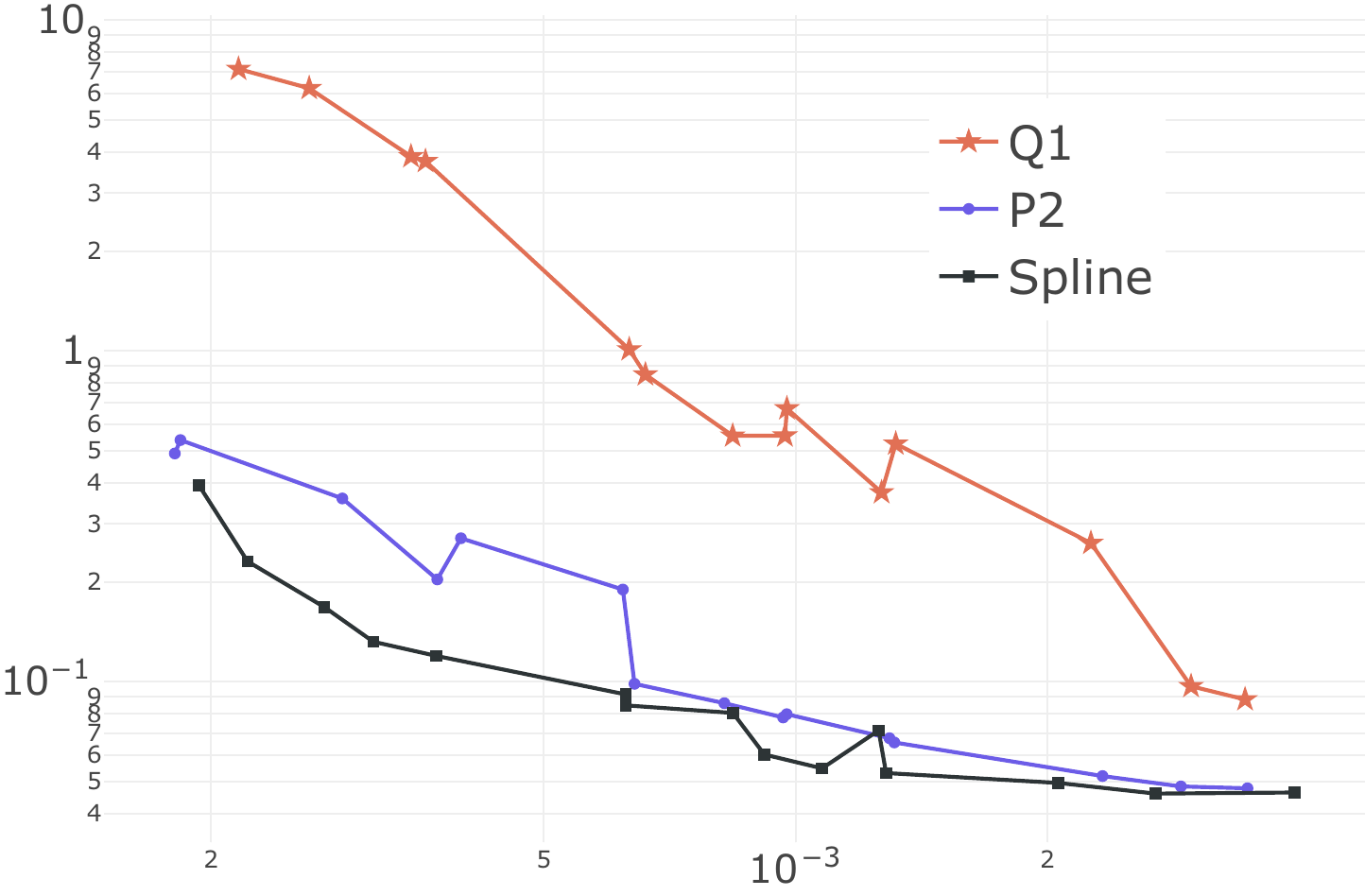}}\hfill
    \parbox{.48\linewidth}{\includegraphics[width=\linewidth]{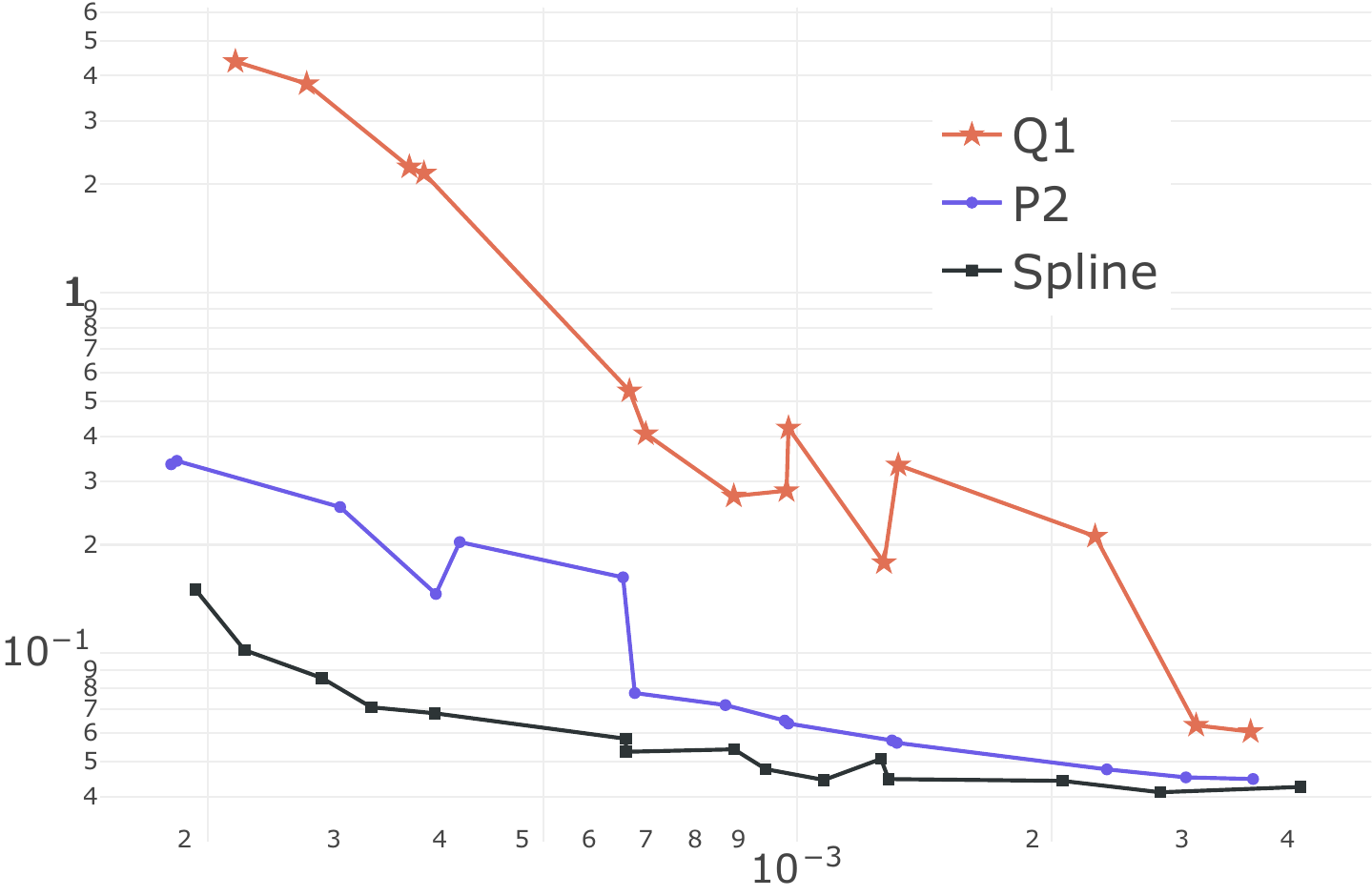}}\\
    Error
    \caption{Time vs. error (with respect to a dense solution) \revision{(total time on the left, and solve time only on the right)} for  $P_2$, $Q_1$, and \revision{quadratic spline elements}.}
    \label{fig:q1-vs-p2}
\end{figure}

We created a sequence of hexahedral and tetrahedral meshes \revisionn{with similar errors} for a force $f=[0, -2, 0]^T$. Figure~\ref{fig:q1-vs-p2} shows that for a given error, $P_2$ discretization is around four times faster than $Q_1$, \revision{and $\mathrm{SPLINE}_2$ have a slight advantage over $P_2$. Note that both $Q_1$ and $\mathrm{SPLINE}_2$ are constructed over a perfectly regular grid, while the $P_2$ elements are defined over an unstructured tetrahedral mesh.}

Finally, we created a sequence of hexahedral meshes that matches total time, total memory, total number of degrees of freedom, and error of the tetrahedral mesh in Figure~\ref{fig:square_beam} for both linear and quadratic elements.  Table~\ref{tab:p-vs-q-fixed} summarizes our findings: $P_1$ is significantly worse than $Q_1$ but the two quadratic discretizations produce similar results. $P_2$ overall performs  better than $Q_1$.

\begin{table}\centering\scriptsize
    \begin{tabular}{lr|cccc}
                                   &        & time (s) & memory (MB) & DOF & error \\
        \hline
        \multirow{4}{*}{\rotatebox{90}{$P_1$/$Q_1$}}
        %%%%%%%%%%%%%%%%%%%%%%%%%%%%%%%%%%%%%
                                   & time   &
        \diagcol{1.01~/~0.98}      &
        \goodcol{125}~/~132        &
        8,258~/~\goodcol{6,413}    &
        1.98e-03~/~\goodcol{6.93e-04}                                              \\
        %%%%%%%%%%%%%%%%%%%%%%%%%%%%%%%%%%%%%
                                   & memory &
        1.01~/~\goodcol{0.91}      &
        \diagcol{125~/~125}        &
        8,258~/~\goodcol{6,050}    &
        1.98e-03~/~\goodcol{7.49e-04}                                              \\
        %%%%%%%%%%%%%%%%%%%%%%%%%%%%%%%%%%%%%
                                   & DOF    &
        \goodcol{1.01}~/~1.07      &
        \goodcol{125}~/~149        &
        \diagcol{8,258~/~7,139}    &
        1.98e-03~/~\goodcol{6.06e-04}                                              \\
        %%%%%%%%%%%%%%%%%%%%%%%%%%%%%%%%%%%%%
                                   & error  &
        1.01~/~\goodcol{0.12}      &
        125~/~\goodcol{18}         &
        8,258~/~\goodcol{1,224}    &
        \diagcol{1.98e-03~/~1.87e-03}                                              \\
        %%%%%%%%%%%%%%%%%%%%%%%%%%%%%%%%%%%%%
        %%%%%%%%%%%%%%%%%%%%%%%%%%%%%%%%%%%%%
        \hline
        \multirow{3}{*}{\rotatebox{90}{$P_2$/$Q_2$}}
        %%%%%%%%%%%%%%%%%%%%%%%%%%%%%%%%%%%%%
                                   & time   &
        \diagcol{16.86~/~17.56}    &
        2,236~/~\goodcol{2,033}    &
        59,885~/~\goodcol{44,541}  &
        1.85e-05~/~\goodcol{1.24e-05}                                              \\
        %%%%%%%%%%%%%%%%%%%%%%%%%%%%%%%%%%%%%
                                   & memory &
        \goodcol{16.86}~/~18.77    &
        \diagcol{2,236~/~2,241}    &
        59,885~/~\goodcol{48,951}  &
        1.85e-05~/~\goodcol{8.40e-06}                                              \\
        %%%%%%%%%%%%%%%%%%%%%%%%%%%%%%%%%%%%%
                                   & DOF    &
        \goodcol{16.86}~/~24.44    &
        \goodcol{2,236}~/~2,988    &
        \diagcol{59,885~/~59,777}  &
        1.85e-05~/~\goodcol{5.43e-06}                                              \\
        %%%%%%%%%%%%%%%%%%%%%%%%%%%%%%%%%%%%%
                                   & error  &
        16.86~/~\goodcol{11.22}    &
        2,236~/~\goodcol{1,451}    &
        59,885~/~\goodcol{35,017}  &
        \diagcol{1.85e-05~/~1.70e-05}                                              \\
        %%%%%%%%%%%%%%%%%%%%%%%%%%%%%%%%%%%%%
        %%%%%%%%%%%%%%%%%%%%%%%%%%%%%%%%%%%%%
        \hline
        \multirow{3}{*}{\rotatebox{90}{$P_2$/$Q_1$}}
        %%%%%%%%%%%%%%%%%%%%%%%%%%%%%%%%%%%%%
                                   & time   &
        \diagcol{16.86~/~16.74}    &
        \goodcol{2,236}~/~2,630    &
        59,885~/~\goodcol{58,719}  &
        \goodcol{1.85e-05}~/~1.58e-04                                              \\
        %%%%%%%%%%%%%%%%%%%%%%%%%%%%%%%%%%%%%
                                   & memory &
        16.86~/~\goodcol{14.26}    &
        \diagcol{2,236~/~2,226}    &
        59,885~/~\goodcol{52,272}  &
        \goodcol{1.85e-05}~/~1.70e-04                                              \\
        %%%%%%%%%%%%%%%%%%%%%%%%%%%%%%%%%%%%%
                                   & DOF    &
        \goodcol{16.86}~/~17.52    &
        \goodcol{2,236}~/~2,669    &
        \diagcol{59,885~/~59,777}  &
        \goodcol{1.85e-05}~/~1.54e-04                                              \\
        %%%%%%%%%%%%%%%%%%%%%%%%%%%%%%%%%%%%%
                                   & error  &
        \goodcol{16.86}~/~170.29   &
        \goodcol{2,236}~/~11,805   &
        \goodcol{59,885}~/~180,774 &
        \diagcol{1.85e-05~/~6.36e-05}                                              \\
        %%%%%%%%%%%%%%%%%%%%%%%%%%%%%%%%%%%%%
        %%%%%%%%%%%%%%%%%%%%%%%%%%%%%%%%%%%%%
    \end{tabular}
    \caption{Comparison of performance of tetrahedral and hexahedral elements on several measures: time, memory, DOF and error, with one of the measures matched \revisionn{(marked in gray)}: for the first row of each comparison, we match time, second memory, etc. The best-performing (according to each measure) element is shown in green.  For instance, by comparing $P_2$ with $Q_1$ for the same error (last section of the table), $P_2$ is faster (fist column), it uses less memory (second column), and it has less DOFs (third column).}
    \label{tab:p-vs-q-fixed}
\end{table}

\begin{figure}\centering\footnotesize
    \parbox{0.7\linewidth}{\centering\footnotesize
        \begin{tabular}{c|lll|l||c}
                  & $t_b$             & $t_a$             & $t_s$             & $t$               & $e_f$             \\
            \hline
            $P_1$ & 3.52e-2           & \goodcol{7.52e-2} & \goodcol{1.21e-1} & \goodcol{2.31e-1} & 3.50e-3           \\
            $P_2$ & 9.88e-2           & 8.58e-1           & 1.79              & 2.75              & \goodcol{5.21e-5} \\
            $Q_1$ & \goodcol{2.22e-2} & 1.03e-1           & 1.76e-1           & 3.02e-1           & 9.82e-4           \\
            $Q_2$ & 5.78e-2           & 1.71              & 2.77              & 4.54              & 8.38e-5           \\
        \end{tabular}}\hfill
    \parbox{.28\linewidth}{
        \includegraphics[width=0.48\linewidth]{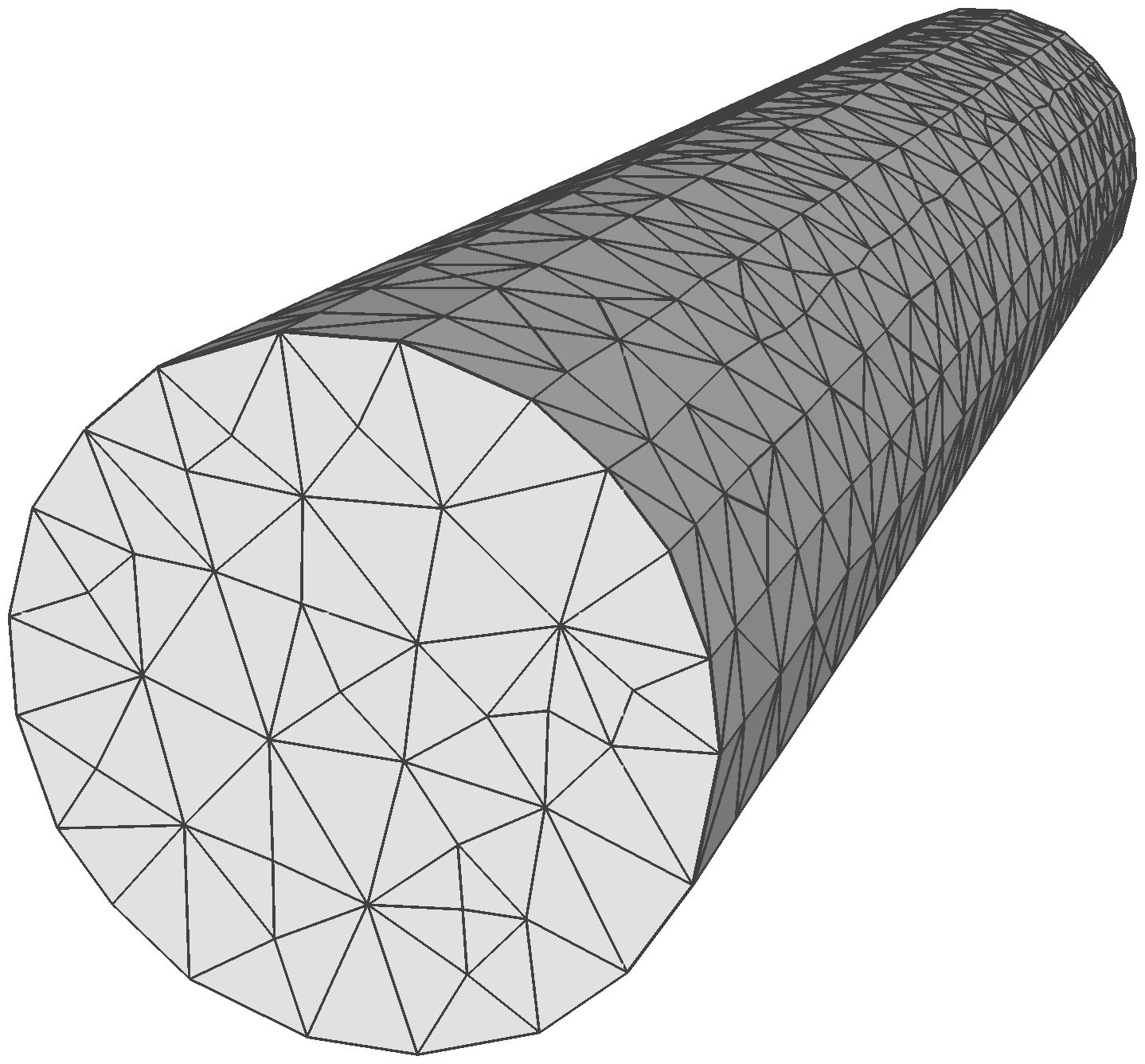}\hfill
        \includegraphics[width=0.48\linewidth]{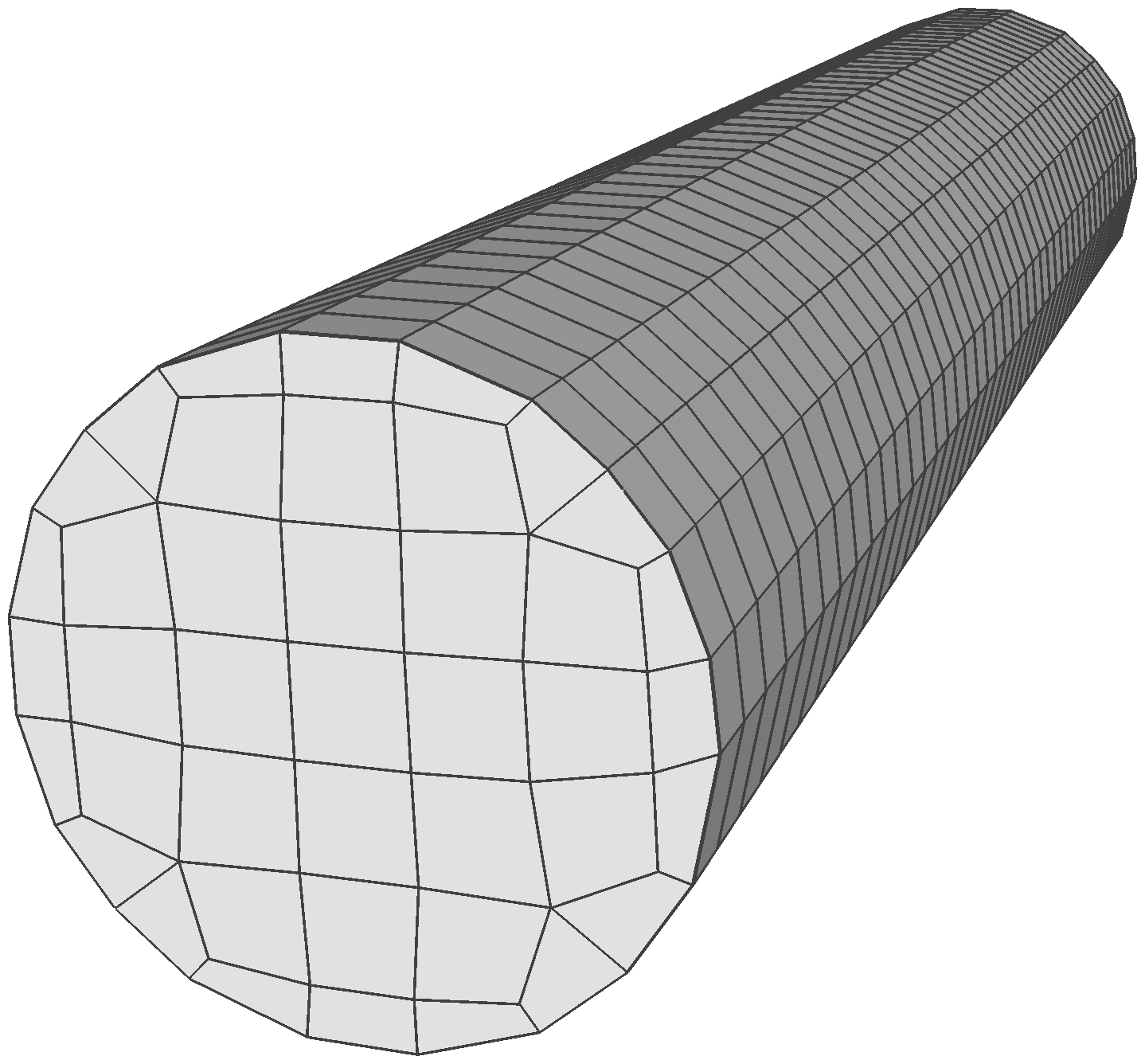}\par
        \parbox{0.48\linewidth}{\centering Tet. mesh}\hfill
        \parbox{0.48\linewidth}{\centering Hex. mesh}\par
    }
    \caption{Displacement error with respect to a dense solution per force unit at the endpoint of a beam with a circular cross-section. The times are averaged over $10$ runs per force sample.}
    \label{fig:circle_beam}
\end{figure}

\paragraph{Circular Cross-section} We consider a beam of length $L=100$ with a circular cross-section of diameter $d=20$. We created a tetrahedral mesh with $2\,252$ vertices and a hexahedral mesh with $2\,288$ vertices (note that in this case the mesh is not a regular grid anymore), \revisionn{by extruding a quad mesh generated with~\cite{Jakob2015Instant}.}
Figure~\ref{fig:circle_beam} shows similar $y$-displacement errors as for the square cross-section, $P_1$ produces low-quality results, while $P_2$ and $Q_2$ are similar.

\begin{figure}\centering\footnotesize
    \parbox{0.7\linewidth}{\centering\footnotesize
        \begin{tabular}{c|lll|l||c}
                  & $t_b$             & $t_a$             & $t_s$             & $t$               & $e_f$             \\
            \hline
            $P_1$ & 9.29e-2           & \goodcol{1.99e-1} & \goodcol{2.67e-1} & \goodcol{5.58e-1} & 1.85e-3           \\
            $P_2$ & 2.68e-1           & 1.94              & 4.59              & 6.80              & \goodcol{7.71e-5} \\
            $Q_1$ & \goodcol{6.16e-2} & 3.12e-1           & 5.28e-1           & 9.01e-1           & 6.77e-4           \\
            $Q_2$ & 1.54e-1           & 4.58              & 9.63              & 1.44e1            & 1.05e-4           \\
        \end{tabular}}\hfill
    \parbox{.28\linewidth}{
        \includegraphics[width=0.48\linewidth]{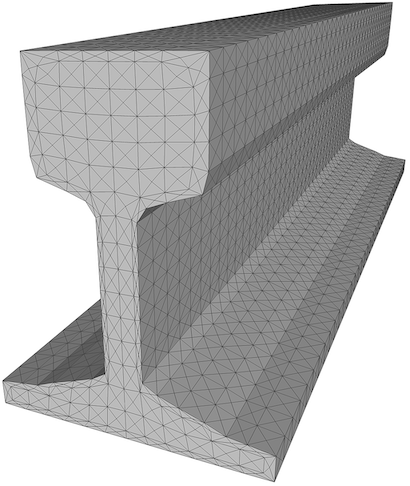}\hfill
        \includegraphics[width=0.48\linewidth]{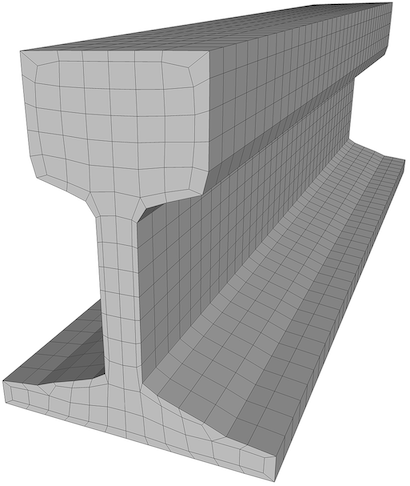}\par
        \parbox{0.48\linewidth}{\centering Tet. mesh}\hfill
        \parbox{0.48\linewidth}{\centering Hex. mesh}\par
    }
    \caption{Displacement errors for a unit force applied at the endpoint of an I-beam. The times are averaged over $10$ runs per force sample.}
    \label{fig:rail_beam}
\end{figure}

\paragraph{I-beam Cross-section}
We use an I-beam (the bounding box of the cross-section is $125\times 154$) of length $L=473.11$. The tetrahedral mesh has $6\,102$ while the hexahedral mesh, \revisionn{generated by extruding a quad mesh,} has $6\,080$ vertices, results are shown in  Figure~\ref{fig:rail_beam}.

\subsection{\revision{Orthotropic Material}}\label{sec:ortho}

\begin{figure}\centering\footnotesize
    \parbox{0.7\linewidth}{\centering\footnotesize
        \begin{tabular}{c|lll|l||c}
                  & $t_b$             & $t_a$             & $t_s$             & $t$               & $e_f$             \\
            \hline
            $P_1$ & 8.00e-3           & \goodcol{2.87e-1} & \goodcol{2.11e-2} & \goodcol{3.16e-1} & 2.47              \\
            $P_2$ & 2.48e-2           & 1.71              & 3.37e-1           & 2.08              & 4.81e-2           \\
            $Q_1$ & \goodcol{6.41e-3} & 9.65e-1           & 3.40e-2           & 1.01              & 1.35              \\
            $Q_2$ & 1.66e-2           & 9.80              & 4.73e-1           & 1.03e1            & \goodcol{2.17e-2} \\
        \end{tabular}}\hfill
    \parbox{.28\linewidth}{\centering\footnotesize
        \includegraphics[width=.48\linewidth]{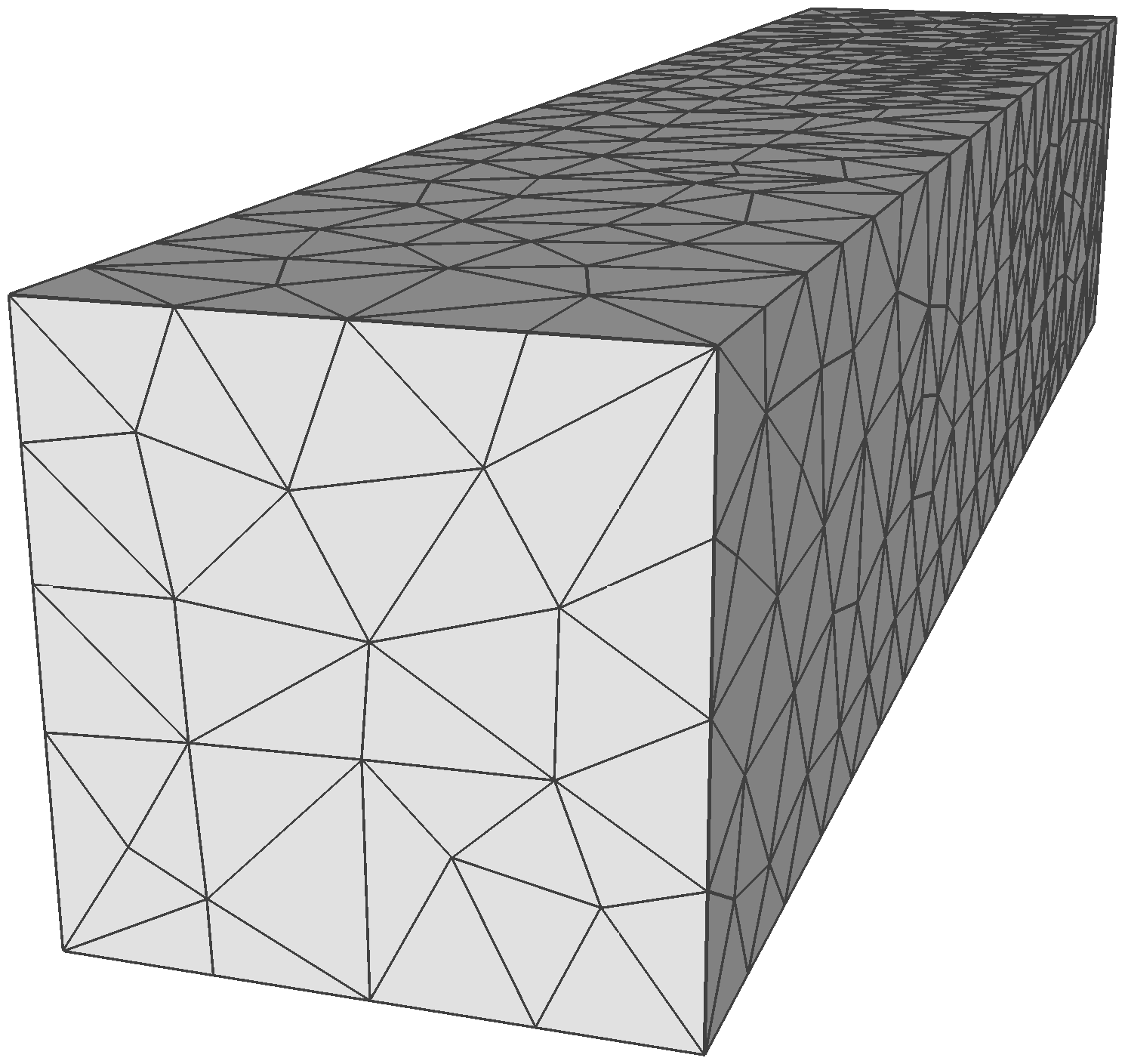}\hfill
        \includegraphics[width=.48\linewidth]{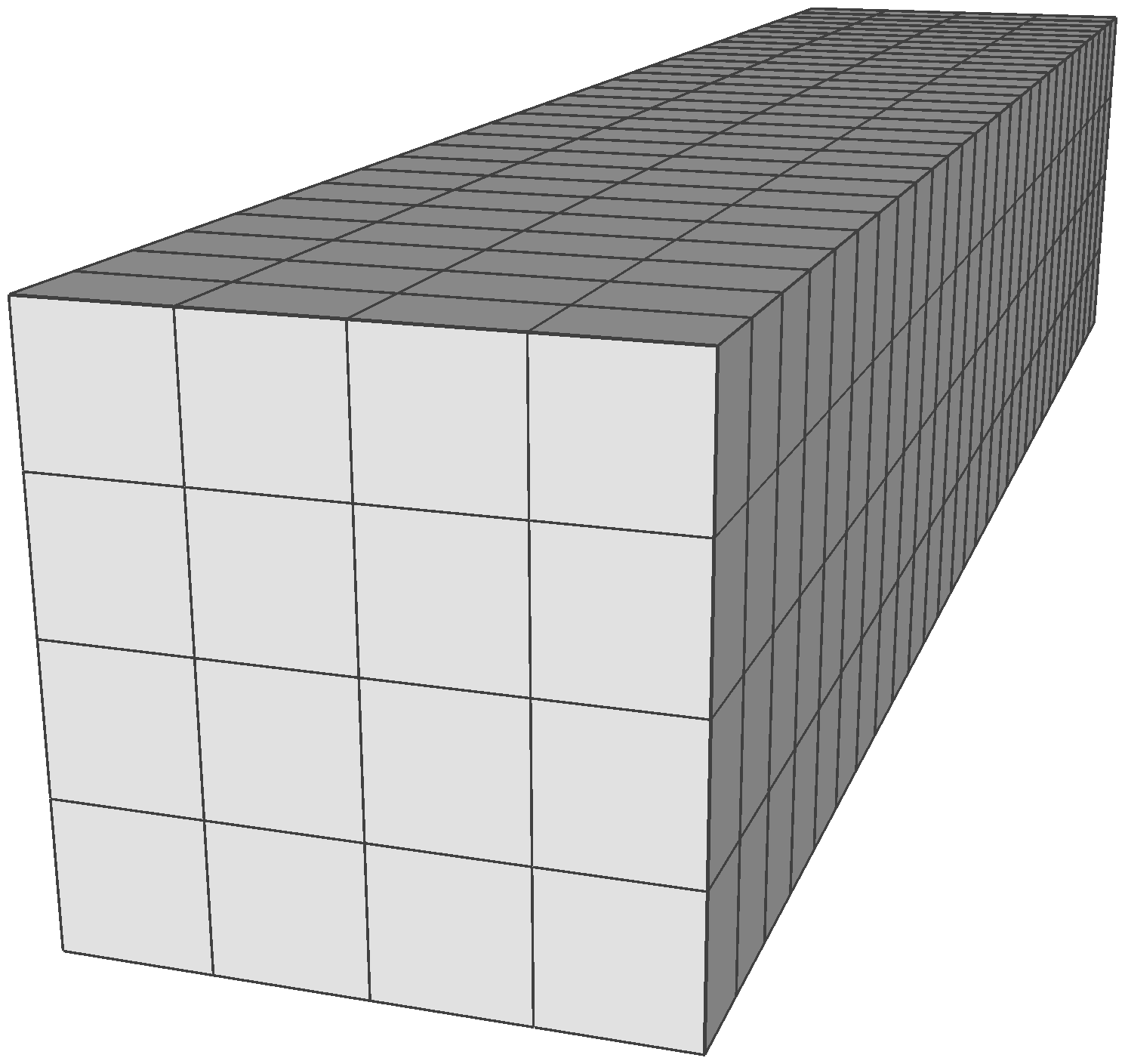}\par
        \parbox{0.48\linewidth}{\centering Tet. mesh}\hfill
        \parbox{0.48\linewidth}{\centering Hex. mesh}\par
    }
    \caption{\revision{
            Displacement error (compared with $P_4$)  for a force of magnitude $1$e$-5$ applied at the endpoint of a beam with orthotropic material and a square cross-section. The times are averaged over $10$ runs per force sample.}}
    \label{fig:ortho}
\end{figure}

\revision{We repeated the previous experiment using linear orthotropic material parameters (carbon fiber). The material parameters are obtained from~\cite{Pardini:2010:cf}. Three Young moduli are  167, 33, and 33, The Poisson ratios are  0.18, 0.25, and 0.18, shear moduli are  13, 21, and 21. Figure~\ref{fig:ortho} shows that the $y$-displacement error with respect to different discretizations has the same behavior as for isotropic materials (Section~\ref{sec:bened_bar}).}

\subsection{\revision{High Aspect-Ratio}}\label{sec:high-aspect}

\begin{figure}
    \centering\footnotesize
    \parbox{.02\linewidth}{\rotatebox{90}{\centering\scriptsize $P_1$, $P_2$}}\hfill
    \parbox{.97\linewidth}{
        \parbox{.22\linewidth}{\centering $4\times 20$}
        \parbox{.22\linewidth}{\centering $2\times 20$}\hfill
        \parbox{.22\linewidth}{\centering $1\times 20$}\hfill\hfill\hfill
        \parbox{.22\linewidth}{\centering $1\times 20^*$}\par
        %%%%%%%%%%%%%%%%%%%%%%%%%%%%%%%%%%%%%%%%%%%%%%%%%%%%%%%%%%%%%%%%%%%%%%%%%%
        \includegraphics[width=.22\linewidth]{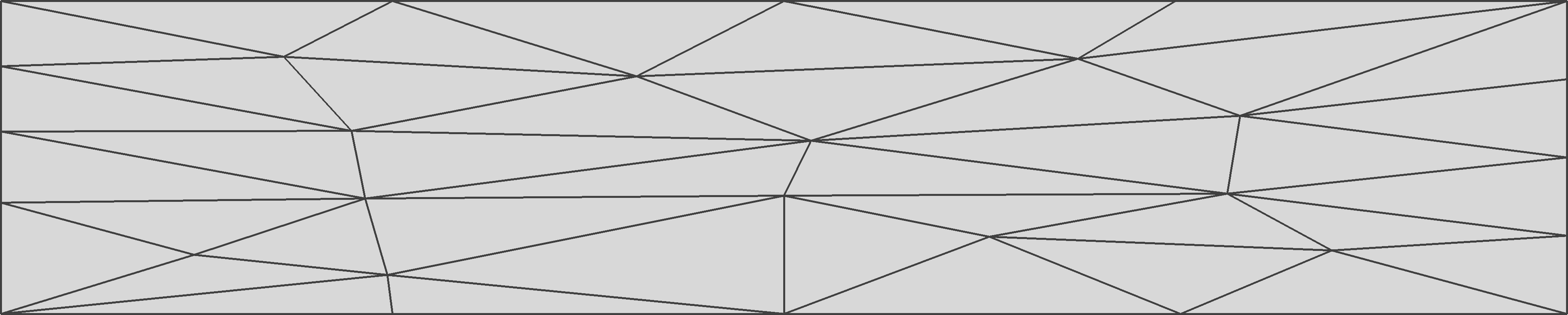}\hfill
        \includegraphics[width=.22\linewidth]{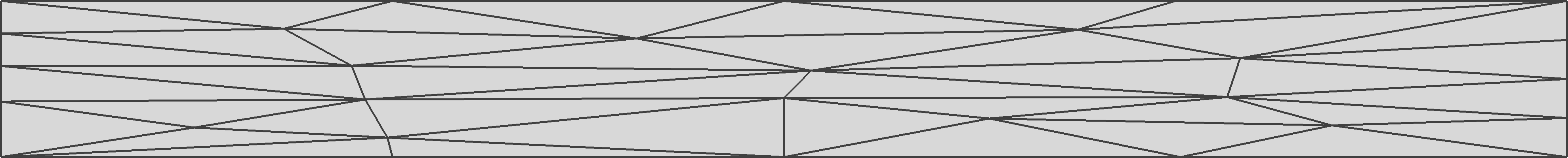}\hfill
        \includegraphics[width=.22\linewidth]{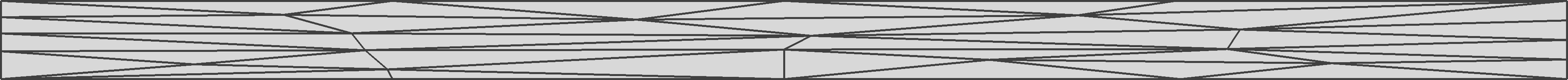}\hfill\hfill\hfill
        \includegraphics[width=.22\linewidth]{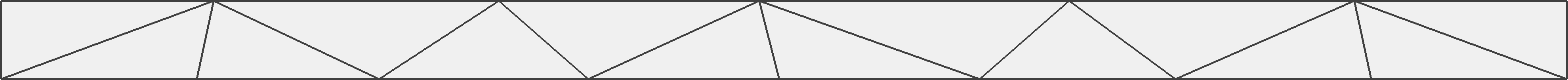}\par\vspace{-0.4em}
        \parbox{.22\linewidth}{\centering AR = 24}\hfill %578
        \parbox{.22\linewidth}{\centering AR = 48}\hfill %2\,313
        \parbox{.22\linewidth}{\centering AR = 96}\hfill\hfill\hfill %9\,251
        \parbox{.22\linewidth}{\centering AR = 9}\par %79
    }\\[1em]
    %%%%%%%%%%%%%%%%%%%%%%%%%%%%%%%%%%%%%%%%%%%%%%%%%%%%%%%%%%%%%%%%%%%%%%%%%
    \parbox{.02\linewidth}{\rotatebox{90}{\centering\scriptsize $P_1^\star$, $P_2^\star$}}\hfill
    \parbox{.97\linewidth}{
        \parbox{.22\linewidth}{~}\hfill
        \parbox{.22\linewidth}{~}\hfill
        \includegraphics[width=.22\linewidth]{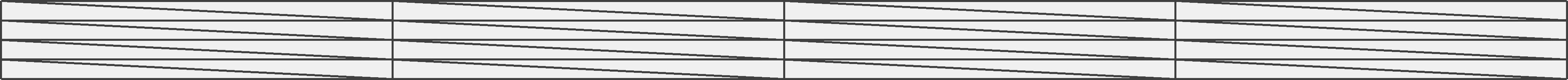}\hfill\hfill\hfill
        \includegraphics[width=.22\linewidth]{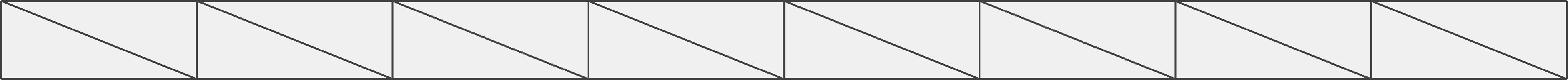}\par\vspace{-0.4em}
        \parbox{.22\linewidth}{~}\hfill
        \parbox{.22\linewidth}{~}\hfill
        \parbox{.22\linewidth}{\centering AR = 40}\hfill\hfill\hfill %1\,604
        \parbox{.22\linewidth}{\centering AR = 5}\par %29
    }\\[1em]
    %%%%%%%%%%%%%%%%%%%%%%%%%%%%%%%%%%%%%%%%%%%%%%%%%%%%%%%%%%%%%%%%%%%%%%%%%
    \parbox{.02\linewidth}{\rotatebox{90}{\centering\scriptsize $Q_1$, $Q_2$}}\hfill
    \parbox{.97\linewidth}{
        \includegraphics[width=.22\linewidth]{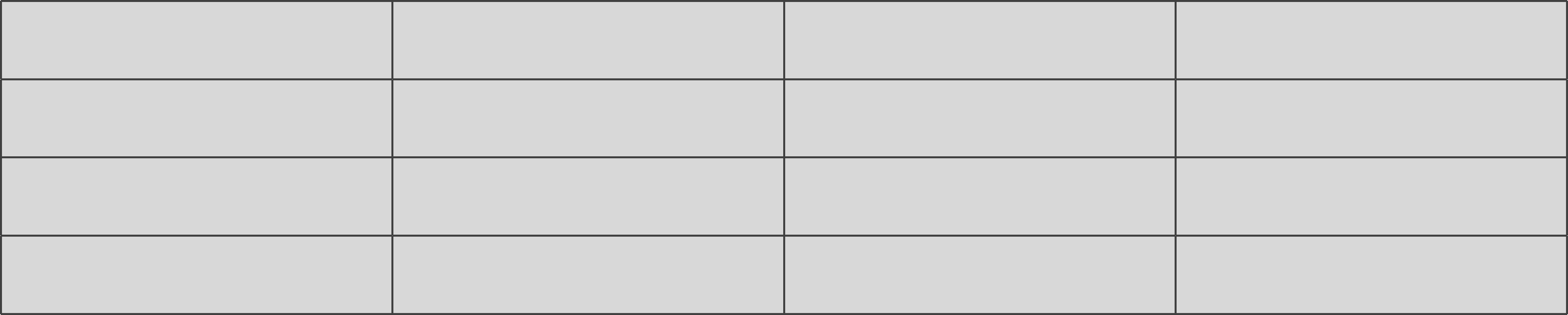}\hfill
        \includegraphics[width=.22\linewidth]{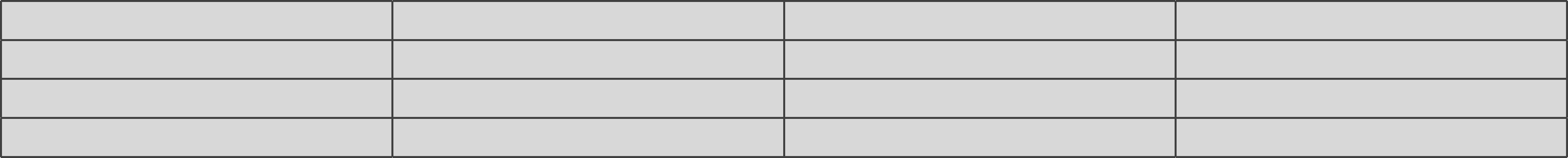}\hfill
        \includegraphics[width=.22\linewidth]{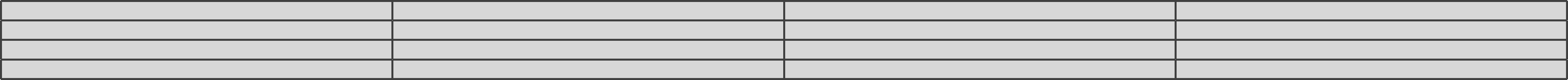}\hfill\hfill\hfill
        \includegraphics[width=.22\linewidth]{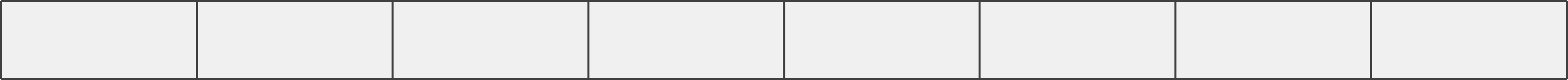}\par\vspace{-0.4em}
        \parbox{.22\linewidth}{\centering AR = 5}\hfill %25
        \parbox{.22\linewidth}{\centering AR = 10}\hfill %100
        \parbox{.22\linewidth}{\centering AR = 20}\hfill\hfill\hfill %400
        \parbox{.22\linewidth}{\centering AR = 2.5}}\\[1em] %6.25
    % \parbox{.33\linewidth}{\centering\footnotesize
    % \parbox[t]{.49\linewidth}{\centering
    %     \includegraphics[width=\linewidth]{figs/aspect/square_beam05}\\ \vspace{-0.8em}
    %     AR = 578\\[0.5em]
    %     \includegraphics[width=\linewidth]{figs/aspect/square_beam10}\\ \vspace{-0.8em}
    %     AR = 2\,313\\[0.5em]
    %     \includegraphics[width=\linewidth]{figs/aspect/square_beam20}\\ \vspace{-0.8em}
    %     AR = 9\,251\\[0.5em]
    %     \includegraphics[width=\linewidth]{figs/aspect/square_beam20_nice}\\ \vspace{-0.8em}
    %     AR = 79\\[0.5em]
    %     \includegraphics[width=\linewidth]{figs/aspect/square_beam20_nice}\\ \vspace{-0.8em}
    %     AR = \todo{xxx + mesh}
    % }\hfill
    % \parbox[t]{.49\linewidth}{\centering
    %     \includegraphics[width=\linewidth]{figs/aspect/square_beam_h05}\\ \vspace{-0.8em}
    %     AR = \todo{xxx}\\[0.5em]
    %     \includegraphics[width=\linewidth]{figs/aspect/square_beam_h10}\\ \vspace{-0.8em}
    %     AR = \todo{xxx}\\[0.5em]
    %     \includegraphics[width=\linewidth]{figs/aspect/square_beam_h20}\\ \vspace{-0.8em}
    %     AR = \todo{xxx}
    % }}\hfill
    \parbox{\linewidth}{\centering\scriptsize
        \begin{tabular}{c|c|lll|l||c}
            \multicolumn{2}{c|}{~}
             & $t_b$       & $t_a$   & $t_s$   & $t$     & $e_f$             \\
            \hline
            \multirow{4}{*}{\rotatebox[origin=c]{90}{$4\times20$}}
             & $P_1$       & 8.78e-3 & 1.96e-2 & 5.72e-2 & 8.56e-2 & 9.91e-1 \\
             & $P_2$       & 2.33e-2 & 1.76e-1 & 3.42e-1 & 5.41e-1 & 5.78e-3 \\
             & $Q_1$       & 5.70e-3 & 3.17e-2 & 6.34e-2 & 1.01e-1 & 2.87e-1 \\
             & $Q_2$       & 1.51e-2 & 4.88e-1 & 4.50e-1 & 9.53e-1 & 1.60e-3 \\
            \hline
            \multirow{4}{*}{\rotatebox[origin=c]{90}{$2\times20$}}
             & $P_1$       & 8.25e-3 & 1.91e-2 & 5.74e-2 & 8.47e-2 & 4.58    \\
             & $P_2$       & 2.30e-2 & 1.71e-1 & 3.27e-1 & 5.21e-1 & 8.56e-2 \\
             & $Q_1$       & 5.69e-3 & 3.10e-2 & 6.29e-2 & 9.95e-2 & 2.54    \\
             & $Q_2$       & 1.48e-2 & 4.66e-1 & 4.33e-1 & 9.14e-1 & 8.40e-3 \\
            \hline\multirow{6}{*}{\rotatebox[origin=c]{90}{$1\times20$}}
             & $P_1$       & 8.94e-3 & 1.93e-2 & 5.87e-2 & 8.69e-2 & 1.86e1  \\
             & $P_2$       & 2.25e-2 & 1.72e-1 & 3.31e-1 & 5.25e-1 & 1.65    \\
             & $Q_1$       & 5.85e-3 & 3.26e-2 & 6.53e-2 & 1.04e-1 & 1.53e1  \\
             & $Q_2$       & 1.47e-2 & 5.03e-1 & 4.67e-1 & 9.85e-1 & 5.58e-2 \\
            \cline{2-7}
             & $P_1^\star$ & 6.31e-3 & 1.50e-2 & 5.54e-2 & 7.67e-2 & 1.69e1  \\
             & $P_2^\star$ & 1.85e-2 & 1.31e-1 & 3.74e-1 & 5.24e-1 & 1.03    \\
            \hline
            \hline\multirow{6}{*}{\rotatebox[origin=c]{90}{$1\times20^*$}}
             & $P_1$       & 5.15e-3 & 1.40e-2 & 5.39e-2 & 7.30e-2 & 1.70e1  \\
             & $P_2$       & 1.55e-2 & 1.18e-1 & 2.28e-1 & 3.62e-1 & 6.54e-2 \\
             & $Q_1$       & 3.90e-3 & 2.47e-2 & 6.56e-2 & 9.42e-2 & 1.33e1  \\
             & $Q_2$       & 9.50e-3 & 3.43e-1 & 3.10e-1 & 6.62e-1 & 2.96e-2 \\
            \cline{2-7}
             & $P_1^\star$ & 4.68e-3 & 1.24e-2 & 6.29e-2 & 8.00e-2 & 1.67e1  \\
             & $P_2^\star$ & 1.42e-2 & 1.09e-1 & 3.14e-1 & 4.37e-1 & 1.17e-1 \\

            % \multirow{4}{*}{\rotatebox[origin=c]{90}{$4\times20$}}
            % &$P_1$&	8.93e-3&	2.08e-2&	6.23e-2&	9.21e-2&	9.91e-1\\
            % &$P_2$&	2.38e-2&	1.92e-1&	3.62e-1&	5.78e-1&	5.78e-3\\
            % &$Q_1$&	6.32e-3&	3.55e-2&	6.96e-2&	1.11e-1&	2.87e-1\\
            % &$Q_2$&	1.61e-2&	5.26e-1&	4.74e-1&	1.02&	1.60e-3\\
            % \hline
            % \multirow{4}{*}{\rotatebox[origin=c]{90}{$2\times20$}}
            % &$P_1$&	9.85e-3&	2.21e-2&	6.45e-2&	9.64e-2&	4.58\\
            % &$P_2$&	2.36e-2&	1.83e-1&	3.41e-1&	5.48e-1&	8.56e-2\\
            % &$Q_1$&	5.86e-3&	3.26e-2&	6.50e-2&	1.03e-1&	2.54\\
            % &$Q_2$&	1.48e-2&	4.90e-1&	4.48e-1&	9.53e-1&	8.40e-3\\
            % \hline
            % \multirow{4}{*}{\rotatebox[origin=c]{90}{$1\times20$}}
            % &$P_1$&	8.57e-3&	1.95e-2&	5.91e-2&	8.72e-2&	1.86e1\\
            % &$P_2$&	2.48e-2&	1.98e-1&	3.75e-1&	5.98e-1&	1.65\\
            % &$Q_1$&	5.74e-3&	3.31e-2&	6.59e-2&	1.05e-1&	1.53e1\\
            % &$Q_2$&	1.51e-2&	4.82e-1&	4.41e-1&	9.37e-1&	5.58e-2\\
            % \hline
            % \multirow{4}{*}{\rotatebox[origin=c]{90}{$1\times20$}}
            % &$P_1$&	5.15e-3&	1.40e-2&	5.39e-2&	7.30e-2&	1.70e1\\
            % &$P_2$&	1.55e-2&	1.18e-1&	2.28e-1&	3.62e-1&	6.54e-2\\
            % &\todo{$P_1$}&	5.15e-3&	1.40e-2&	5.39e-2&	7.30e-2&	1.70e1\\
            % &\todo{$P_2$}&	1.55e-2&	1.18e-1&	2.28e-1&	3.62e-1&	6.54e-2\\
        \end{tabular}
    }
    \caption{\revision{
            Displacement errors with respect to a $P_4$ solution for a unit force applied at the endpoint for different aspect ratios. The aspect ratio $1 \times 20^\star$ is the same domain as $1 \times 20$ remeshed with optimized element aspect ratio. The results $P_1^\star$ and $P_2^\star$ are obtained by splitting the hexahedra into 5 tetrahedra.}}
    \label{fig:high-aspect}
\end{figure}

\revision{To analyze the effect of using elements with high aspect ratio, we repeated the previous experiment for 3 different domains by shrinking the height of the square cross-section from 20, to 4, 2, and 1, while keeping the connectivity identical. This procedure  introduces artificial high aspect-ratio elements.  To obtain the same anisotropy measure for hexahedral and tetrahedral meshes, we define the aspect ratio between the largest and smallest eigenvalue of the covariance matrix of its vertices. Figure~\ref{fig:high-aspect} shows the error in the $y$-direction displacement per unit of force for the hexahedral and tetrahedral meshes. The tetrahedral mesh suffers from the low element quality much more than the hexahedral mesh. However, if we regenerate the meshes for the thin domain with the same number of degrees of freedom, but with element quality optimization (last row in Figure~\ref{fig:high-aspect}), the high error of $P_2$ disappears and the errors are similar as for $Q_2$, as shown in last four rows in Figure~\ref{fig:high-aspect}. For comparison, we also generate a tetrahedral mesh by simply splitting the hexahedra into six tetrahedra. Note that this kind of aspect ratios are extreme and do not appear in any \emph{automatically} meshed model in our data set (Figure~\ref{fig:dataset-aspect}).}

\subsection{2D Domain with a Hole}\label{sec:plate-hole}

\begin{figure}\centering\footnotesize
    \parbox{.88\linewidth}{
        \includegraphics[width=0.49\linewidth]{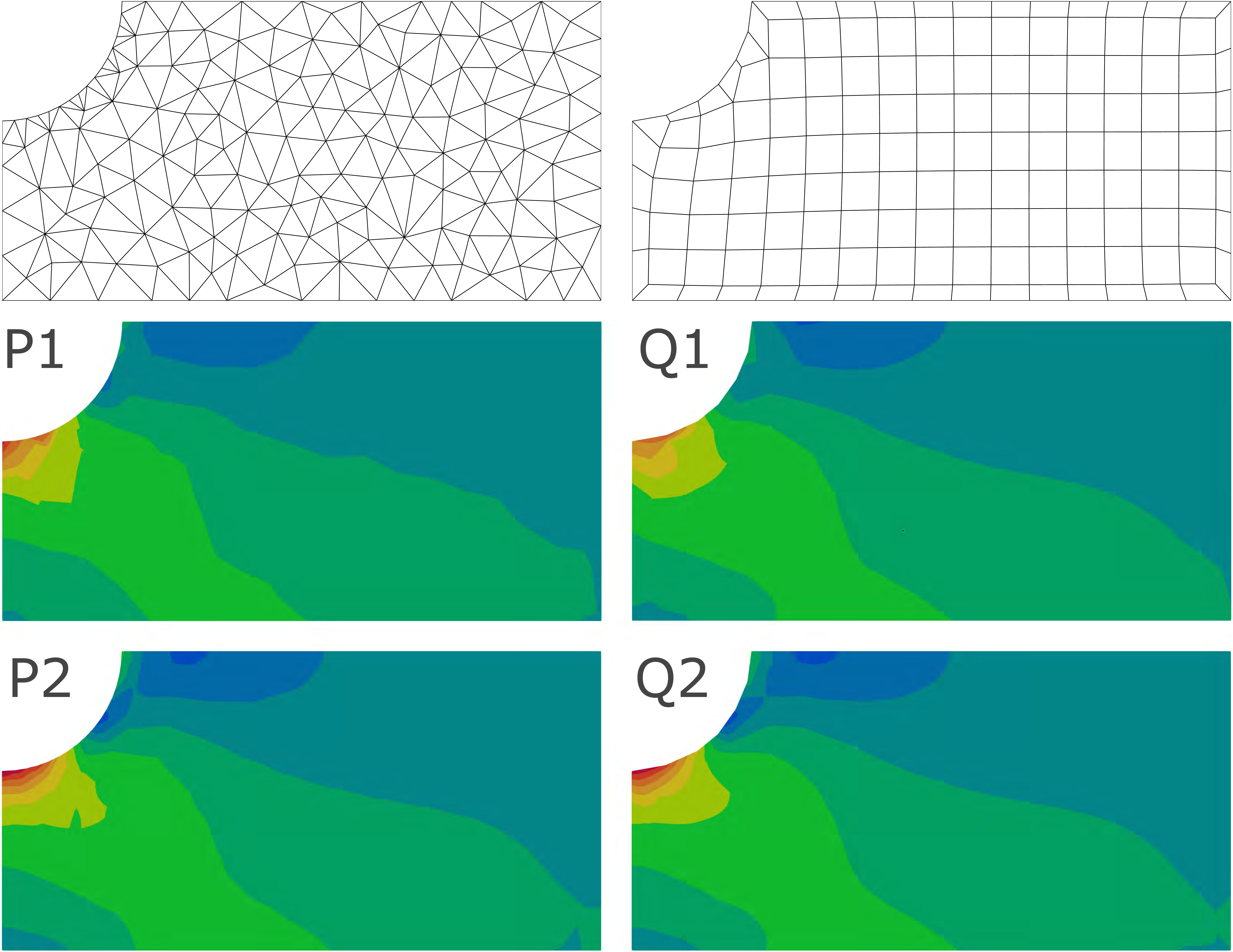}\hfill
        \includegraphics[width=0.49\linewidth]{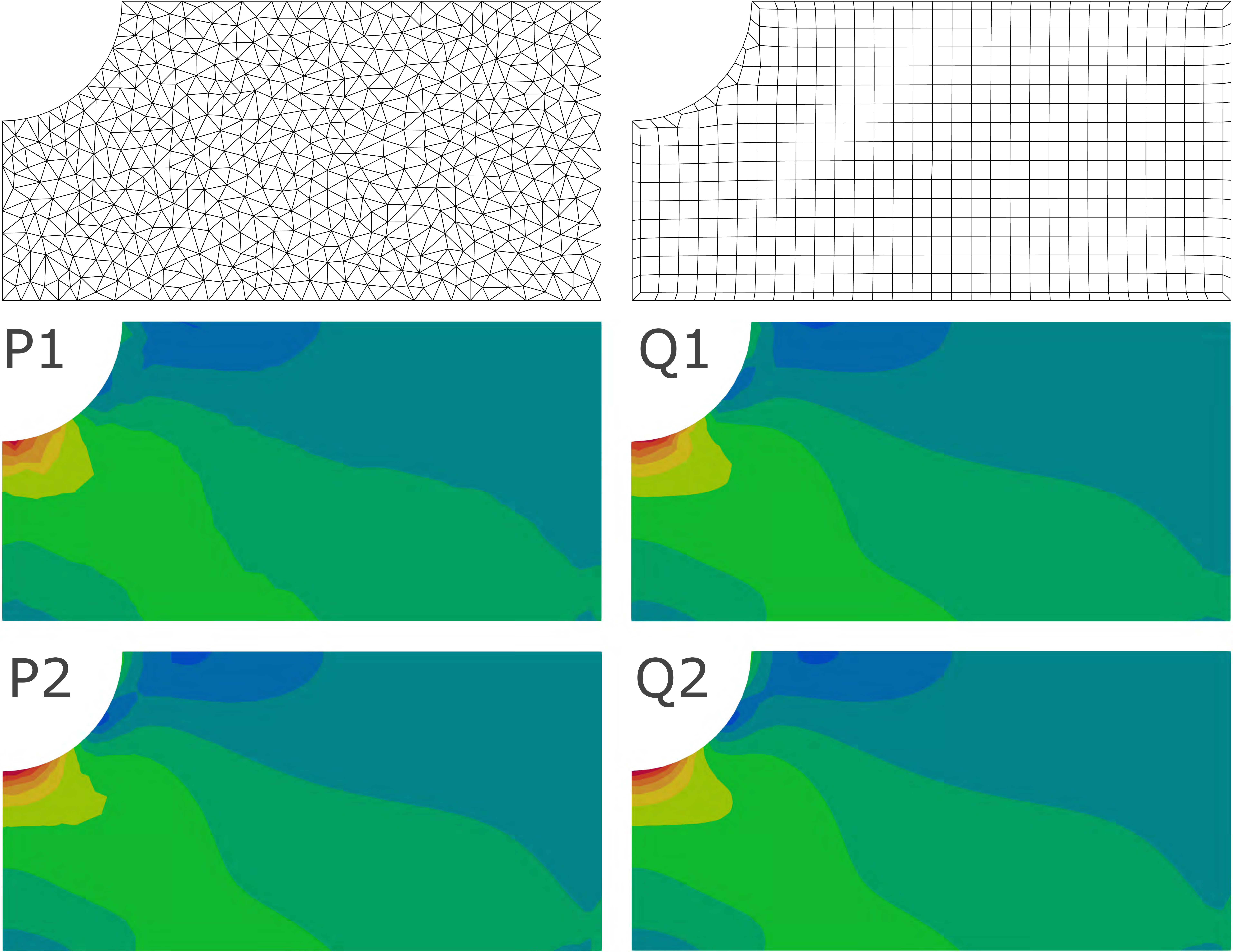}\par
        \parbox{.49\linewidth}{\centering $169$~/~$167$ vertices}\hfill
        \parbox{.49\linewidth}{\centering $622$~/~$558$ vertices}\\[2ex]
        \includegraphics[width=0.49\linewidth]{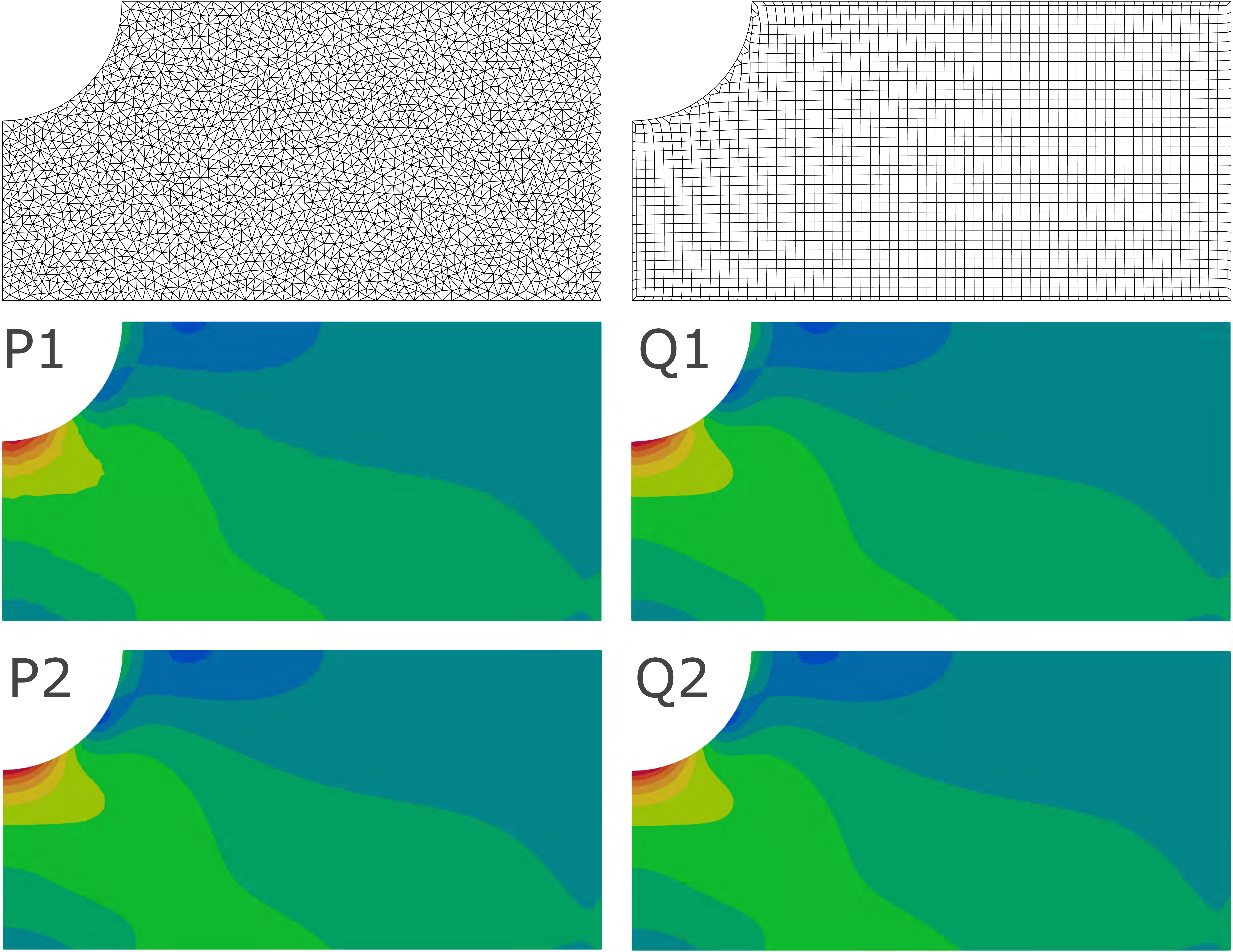}\hfill
        \includegraphics[width=0.49\linewidth]{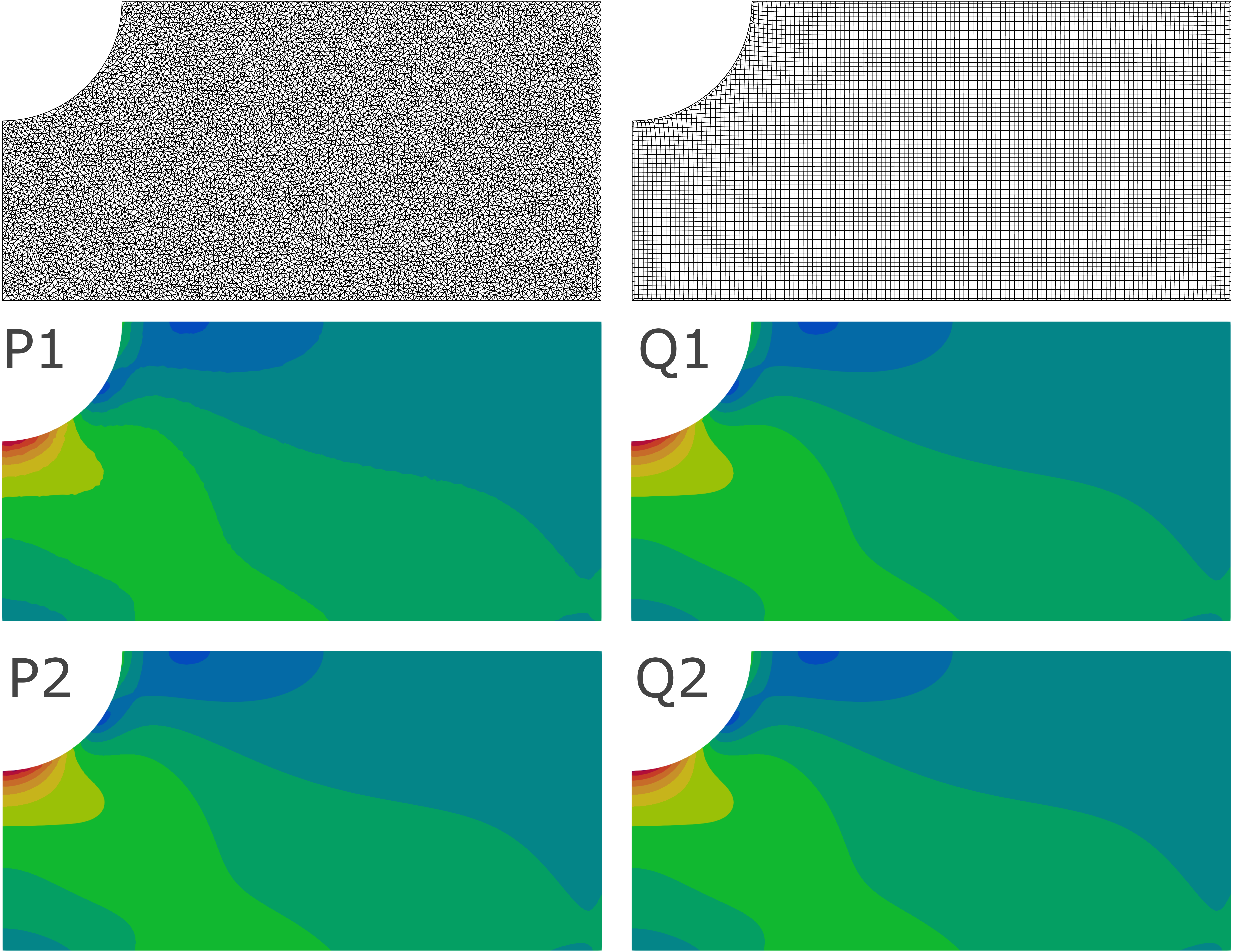}\par
        \parbox{.49\linewidth}{\centering $2\,176$~/~$2\,117$ vertices}\hfill
        \parbox{.49\linewidth}{\centering $8\,549$~/~$8\,504$ vertices}\par
    }\hfill
    \parbox{.1\linewidth}{
        \includegraphics[width=\linewidth]{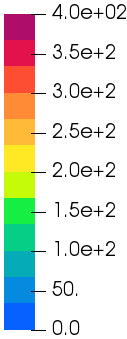}
    }
    \caption{\revisionn{Visualization of the von Mises stresses for four different mesh resolutions. Each figure shows $P_1$ (top left), $P_2$ (bottom left), $Q_1$ (top right), and $Q_2$ (bottom right). The numbers below the figures represent the number of vertices of the tri~/~quad-mesh.}}
    \label{fig:plane-hole}
\end{figure}

Another commonly used test problem is a 2D domain  with a hole in the middle. For our experiments we use a square domain of size $200\times100$ with a hole in the center of radius $20$, the same material model (linear elasticity~\eqref{eq:lin-elast}) and same material parameters $E=210\,000$ and $\nu=0.3$.  The experiment consists of applying an opposite in-plane force on the left and right boundary of $100$, that is, stretching the plane horizontally. This problem is obviously ill-posed because of the lack of Dirichlet boundary conditions. We use a standard approach to eliminate the null-space of solutions by exploiting symmetry, and simulating on a quarter of the domain. This leads to a domain with a ``corner'' cut with two symmetric boundary Dirichlet conditions (displacement is constrained only in the orthogonal direction),  a zero Neumann condition, and a Neumann condition corresponding to the original force. We solve this particular benchmark problem on  four meshes with different resolutions. Figure~\ref{fig:plane-hole} shows the von Mises stresses~\eqref{eq:von-mises} on the top for a triangle mesh and bottom for a quadrilateral mesh, left linear and right quadratic elements. As expected, for a sufficiently dense mesh, all methods converge to similar results. The interesting result is that $Q_2$ elements produce visually better results even at really low resolution (first image and second image). In contrast, for linear triangular elements, we need to increase the mesh resolution up to $8\,500$ vertices (last image) for the artifacts to disappear.

\begin{figure}
    \centering\footnotesize
    \includegraphics[width=.22\linewidth]{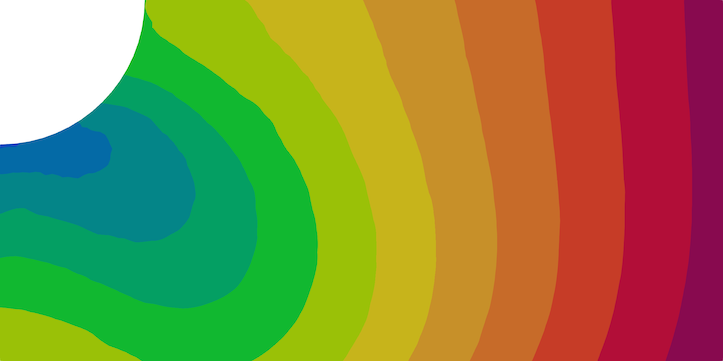}\hfill
    \includegraphics[width=.22\linewidth]{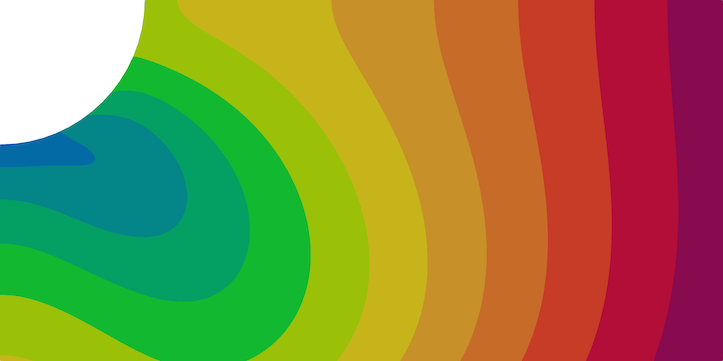}\hfill
    \includegraphics[width=.22\linewidth]{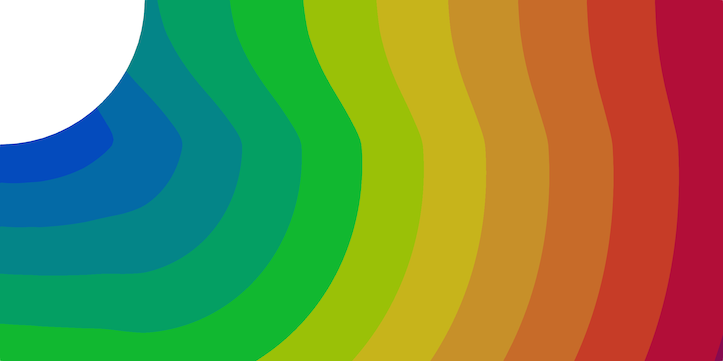}\hfill
    \includegraphics[width=.22\linewidth]{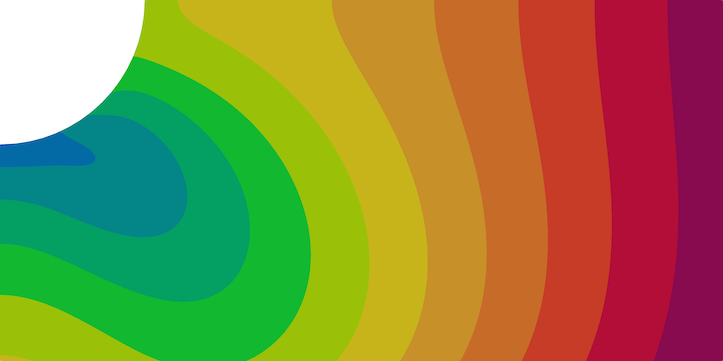}\hfill\hfill
    \includegraphics[width=.04\linewidth]{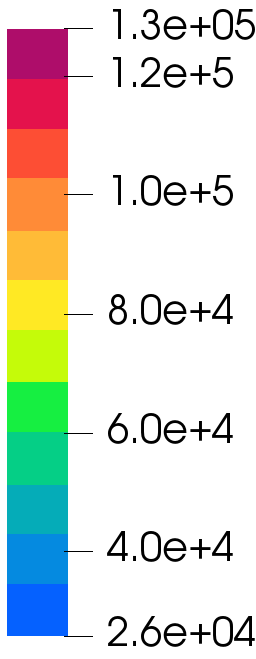}\par
    \parbox{.22\linewidth}{\centering$P_1$}\hfill
    \parbox{.22\linewidth}{\centering$P_2$}\hfill
    \parbox{.22\linewidth}{\centering$Q_1$}\hfill
    \parbox{.22\linewidth}{\centering$Q_2$}\hfill\hfill
    \parbox{.04\linewidth}{~}\par
    \caption{\revision{Displacement norm for a nearly incompressible 2D domain with a hole \revisionn{for a mesh with $8\,549$~/~$8\,504$ vertices}.}}
    \label{fig:inc-plate-hole}
\end{figure}

\revision{This particular problem is also a standard benchmark for incompressible material simulation. We performed the same experiment for a nearly incompressible material: $E=0.1$ and $\nu=0.9999$. Figure~\ref{fig:inc-plate-hole} shows the norm of the displacement: as for the compressible case, $P_2$ and $Q_2$ have a similar behavior. Interestingly, for this case, $Q_1$ produces a very different solution.}

\subsection{Nearly Incompressible Material}\label{sec:incompressible}
\begin{figure}\centering\footnotesize
    \parbox{.91\linewidth}{\centering\footnotesize
        % \parbox{.24\linewidth}{\centering $\nu = 0.9$}\hfill
        % \parbox{.24\linewidth}{\centering $\nu = 0.99$}\hfill
        % \parbox{.24\linewidth}{\centering $\nu = 0.999$}\hfill
        % \parbox{.24\linewidth}{\centering $\nu = 0.9999$}\par
        %
        \parbox{.32\linewidth}{\centering\footnotesize
            \includegraphics[width=.24\linewidth]{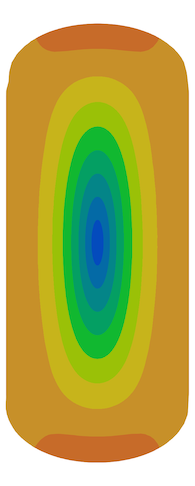}\hfill
            \includegraphics[width=.24\linewidth]{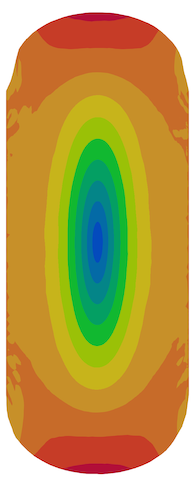}\hfill
            \includegraphics[width=.24\linewidth]{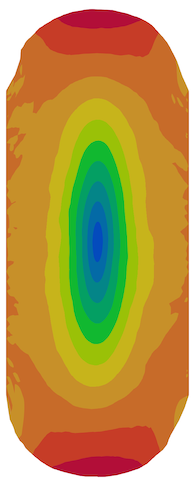}\hfill
            \includegraphics[width=.24\linewidth]{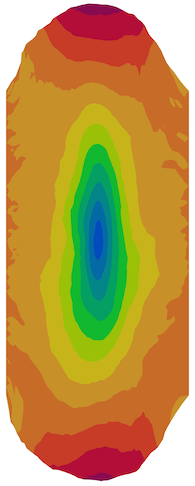}\par
            $P_1$, avg time: 0.10s
        }\hfill
        \parbox{.32\linewidth}{\centering\footnotesize
            \includegraphics[width=.24\linewidth]{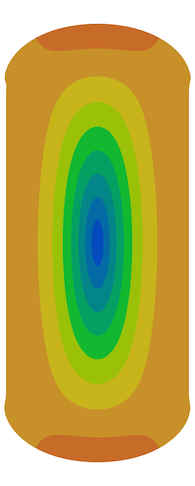}\hfill
            \includegraphics[width=.24\linewidth]{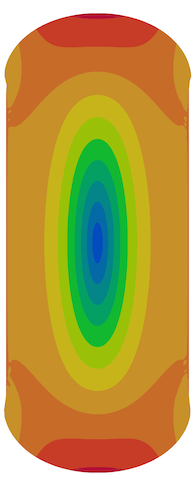}\hfill
            \includegraphics[width=.24\linewidth]{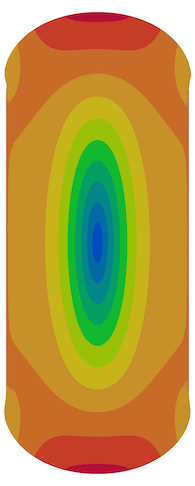}\hfill
            \includegraphics[width=.24\linewidth]{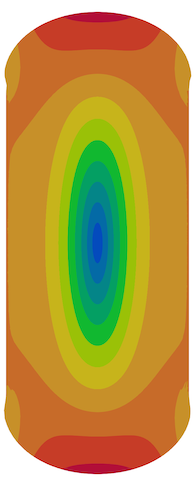}\par
            $P_2$, avg time: 0.56s
        }\hfill
        \parbox{.32\linewidth}{\centering\footnotesize
            \includegraphics[width=.24\linewidth]{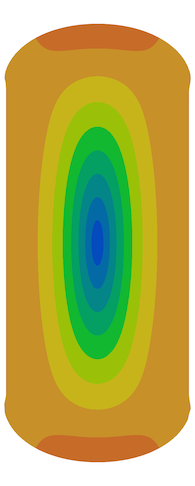}\hfill
            \includegraphics[width=.24\linewidth]{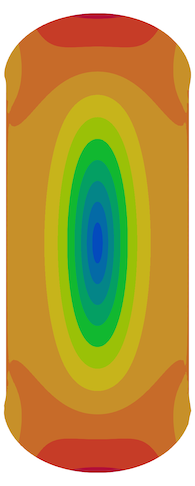}\hfill
            \includegraphics[width=.24\linewidth]{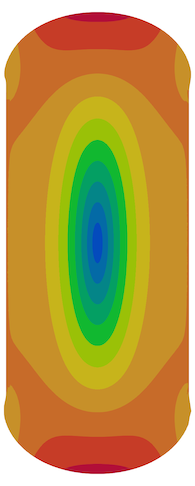}\hfill
            \includegraphics[width=.24\linewidth]{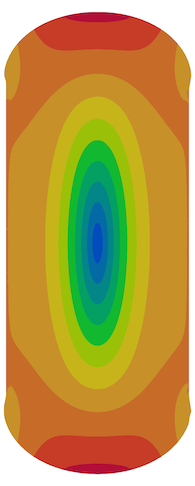}\par
            Mixed, avg time: 1.13s
        }\\[2ex]
        %%%%%%%%%%%%%%%%%%%%%%%%%%%%%%%%%%%%%%%%%%%%%%%%%%%%%%%%%%%%%%%%%%%%%%%%%%%%%%%%%
        \parbox{.32\linewidth}{\centering\footnotesize
            \includegraphics[width=.24\linewidth]{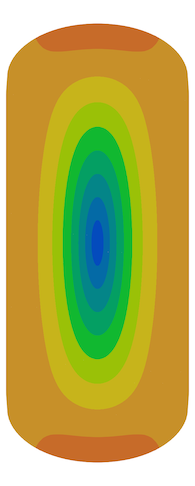}\hfill
            \includegraphics[width=.24\linewidth]{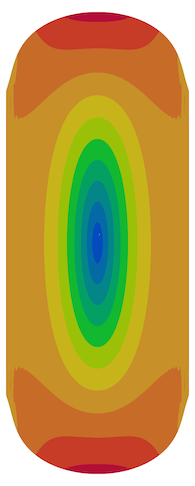}\hfill
            \includegraphics[width=.24\linewidth]{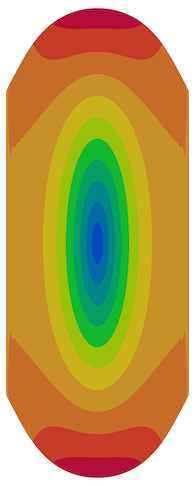}\hfill
            \includegraphics[width=.24\linewidth]{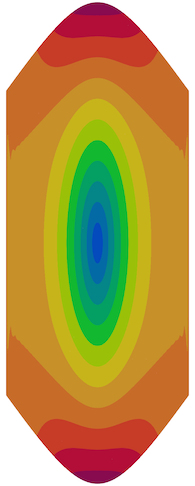}\par
            $Q_1$, avg time: 0.14s
        }\hfill
        \parbox{.32\linewidth}{\centering\footnotesize
            \includegraphics[width=.24\linewidth]{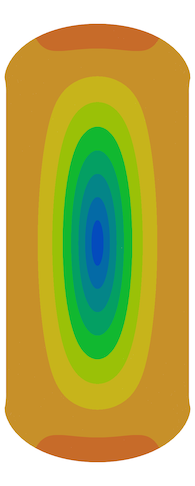}\hfill
            \includegraphics[width=.24\linewidth]{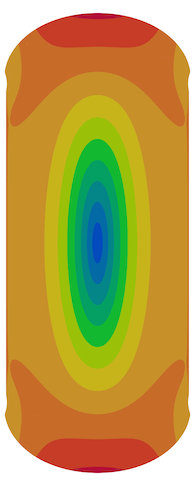}\hfill
            \includegraphics[width=.24\linewidth]{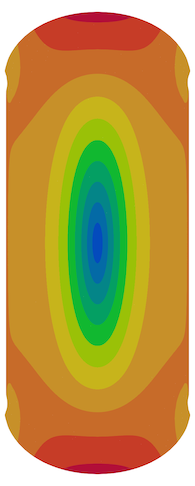}\hfill
            \includegraphics[width=.24\linewidth]{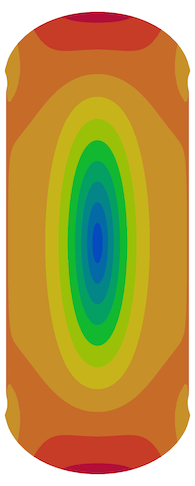}\par
            $Q_2$, avg time: 0.76s
        }\hfill
        \parbox{.32\linewidth}{\centering\footnotesize
            \includegraphics[width=.24\linewidth]{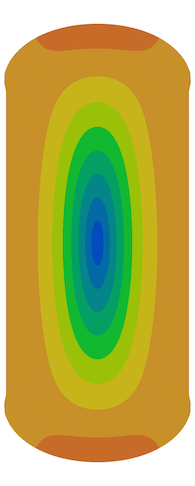}\hfill
            \includegraphics[width=.24\linewidth]{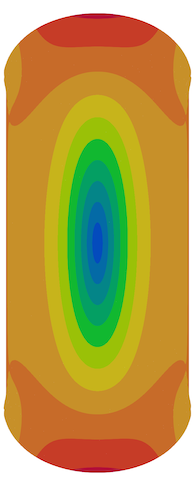}\hfill
            \includegraphics[width=.24\linewidth]{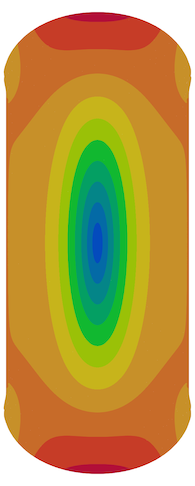}\hfill
            \includegraphics[width=.24\linewidth]{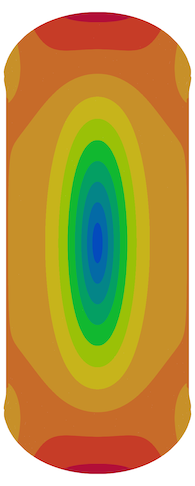}\par
            Mixed, avg time: 2.02s
        }
    }\hfill
    \parbox{.07\linewidth}{
        \includegraphics[width=\linewidth]{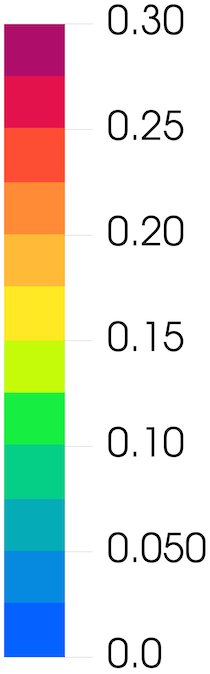}
    }
    \caption{Visualization of the norm of the displacement for a compressed square, with $\nu = 0.9, 0.99, 0.999, 0.9999$, left to right, for different elements and discretizations.}
    \label{fig:incompressible}
\end{figure}
For the last linear benchmark, we compared the performance of the four discretizations with the material approaching incompressibility. We apply a boundary displacement $[0.2, 0]$ on the left and $[-0.2, 0]$ on the right of a unit square. We perform a series of experiments in which  we keep the Young's modulus fixed at $0.1$ while changing the Poisson's ratio from $0.9$ to $0.9999$ ($1$ being the limit of incompressibility in 2D, i.e., area preservation). We compare the standard formulation~\eqref{eq:lin-elast} with a mixed formulation~\eqref{eq:incompressible-elast} that does not become unstable as $\nu\to1$. Note that since mixed formulations require  different  basis degrees for the displacement and the pressure. We performed our experiments using linear pressure bases and quadratic bases for the displacements. We mesh the square with a quad mesh with 4\,225 vertices and a tri mesh with 4\,229 vertices.

Figure~\ref{fig:incompressible} shows the norm of the displacement for this series of experiments. For the nearly incompressible regime (i.e., $\nu=0.9999$) it is remarkable that the quadrilateral element discretization leads to a symmetric and smooth (but incorrect) result for the linear case, while the triangular elements  producing an unstable output. The two quadratic discretizations produce visually similar results, close to those obtained with the stable mixed method. The only quantitative difference is that the residual error for the direct solver drops from 1e$-15$ (numerical zero) to 1e$-12$, indicating that the system is close to singular.

\clearpage
\subsection{Beam with Torsional Loads}\label{sec:torsion}
\begin{figure}\centering\footnotesize
    \parbox{.89\linewidth}{\centering
        \parbox{.45\linewidth}{\centering\footnotesize
            \parbox{0.12\linewidth}{\includegraphics[width=\linewidth]{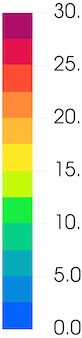}}\hfill
            \parbox{.86\linewidth}{\centering\footnotesize
                \parbox{.22\linewidth}{\centering $P_1$}\hfill
                \parbox{.22\linewidth}{\centering $Q_1$}\hfill
                \parbox{.22\linewidth}{\centering $P_2$}\hfill
                \parbox{.22\linewidth}{\centering $Q_2$}\par
                \includegraphics[width=0.22\linewidth]{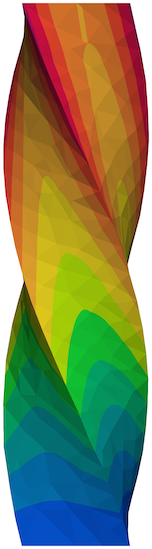}\hfill
                \includegraphics[width=0.22\linewidth]{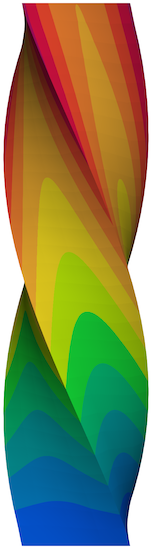}\hfill
                \includegraphics[width=0.22\linewidth]{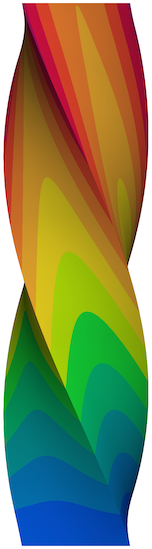}\hfill
                \includegraphics[width=0.22\linewidth]{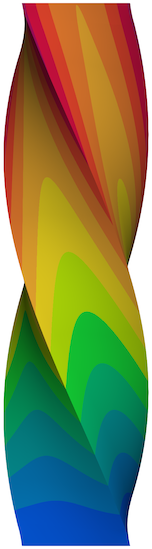}\\[-0.5ex]
                \parbox{.22\linewidth}{\centering\scriptsize 0:00:10}\hfill
                \parbox{.22\linewidth}{\centering\scriptsize 0:00:32}\hfill\textbf{}
                \parbox{.22\linewidth}{\centering\scriptsize 0:01:03}\hfill
                \parbox{.22\linewidth}{\centering\scriptsize 0:07:46}\\[2ex]
                \includegraphics[width=0.22\linewidth]{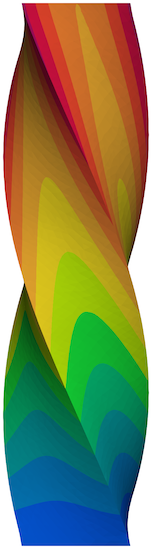}\hfill
                \includegraphics[width=0.22\linewidth]{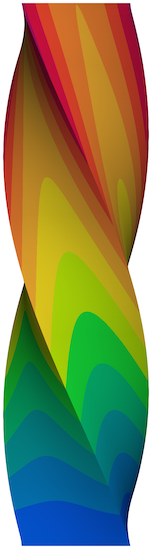}\hfill
                \includegraphics[width=0.22\linewidth]{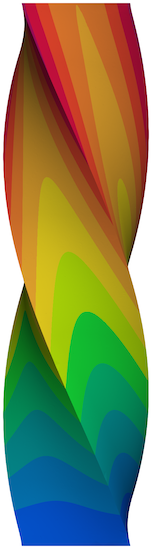}\hfill
                \includegraphics[width=0.22\linewidth]{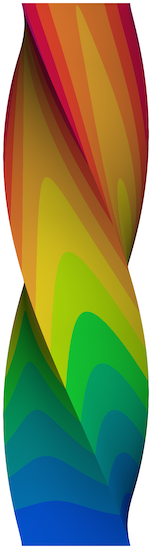}\\[-0.5ex]
                \parbox{.22\linewidth}{\centering\scriptsize 0:14:45}\hfill
                \parbox{.22\linewidth}{\centering\scriptsize 0:46:53}\hfill
                \parbox{.22\linewidth}{\centering\scriptsize 5:19:43}\hfill
                \parbox{.22\linewidth}{\centering\scriptsize 31:31:49}\par
            }
        }\hfill\hfill
        \parbox{.52\linewidth}{\centering\footnotesize
            \includegraphics[width=\linewidth]{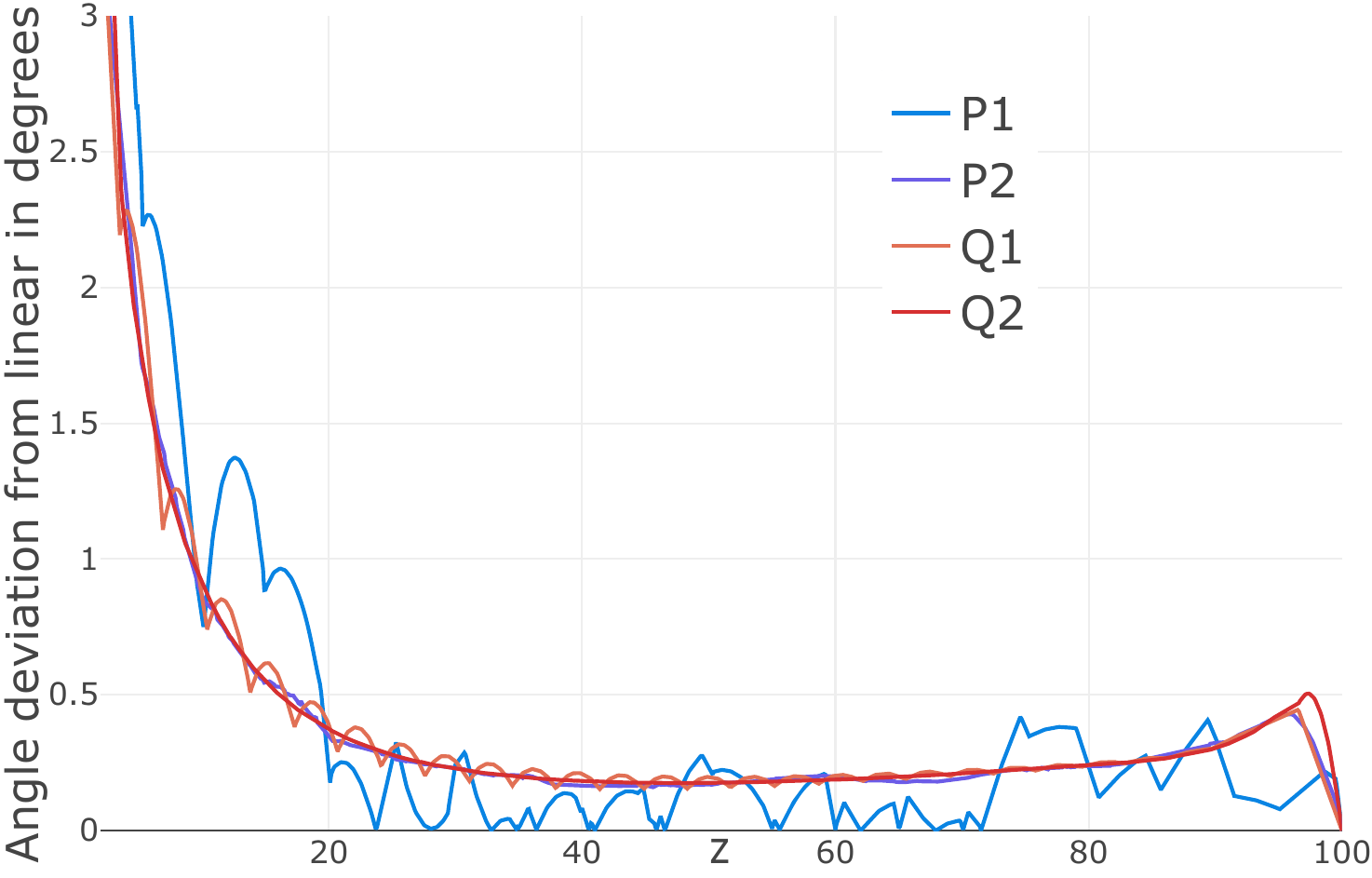}\par
            \includegraphics[width=\linewidth]{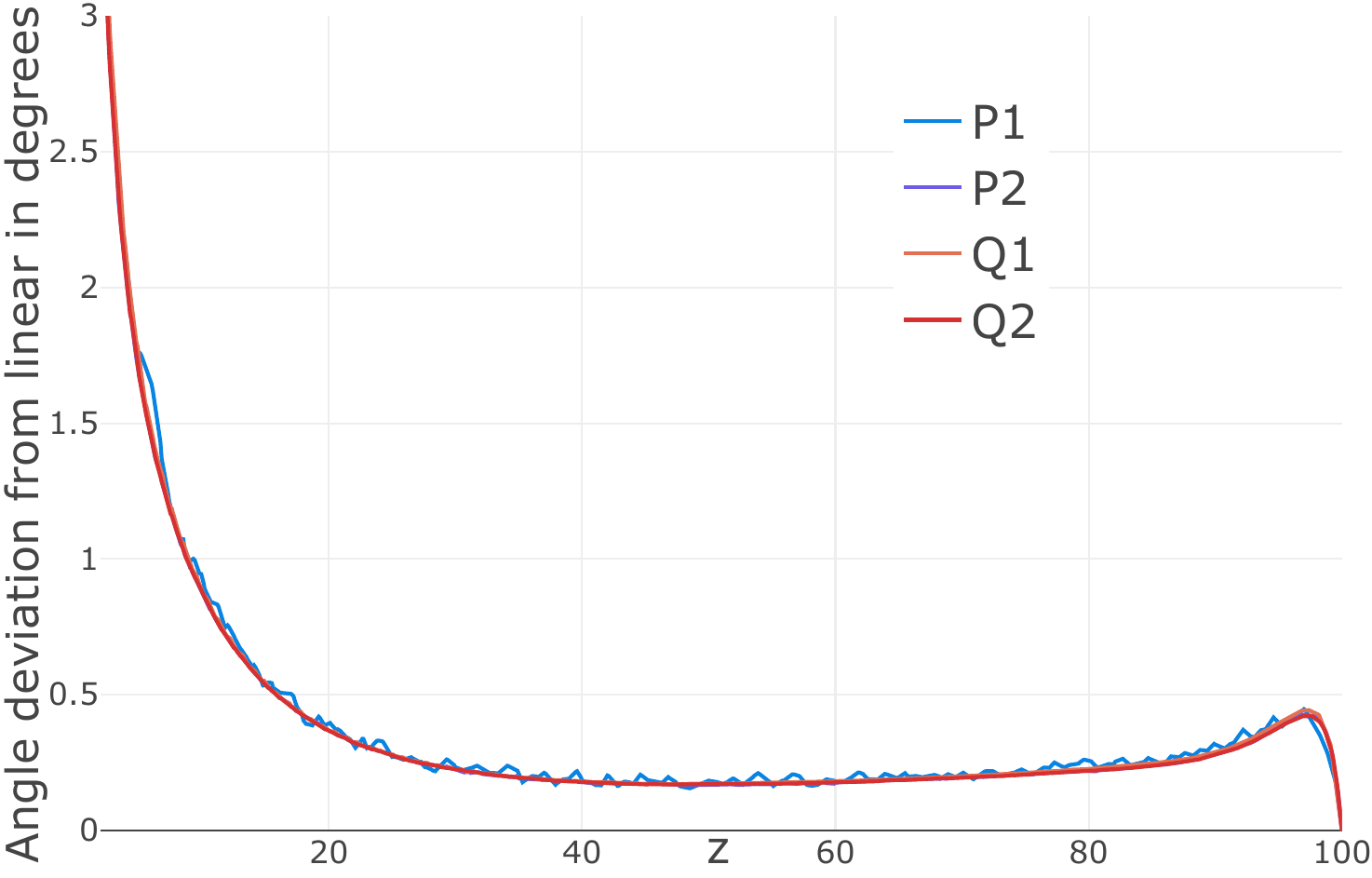}\par
        }}
    \caption{Nonlinear elastic deformation. Left: the deformed mesh color coded by the norm of the displacement and the running time. Right: the angle of rotation of the cross-section  deviation from linearly interpolated along the depth of the bar. We show the deviation from linear interpolation to make the differences between different elements more visible. Note that the linear interpolation is not the exact solution so we do not expect the line to go to zero.}
    \label{fig:torsion}
\end{figure}

\begin{wrapfigure}{r}{0pt}\centering
    \raisebox{0pt}[\dimexpr\height-0.7\baselineskip\relax]{
        \includegraphics[width=0.07514\linewidth]{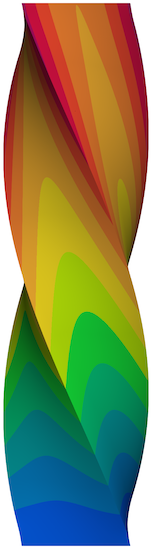}
    }
\end{wrapfigure}
We now compare the solutions for the Neo-Hookean~\eqref{eq:nl-elast} material model for our discretizations. We take a beam with a cross-section $[-10, 10]^2$ and length 100, $E=200$ $nu=0.35$, fix the bottom part and apply a rotation of $90$ degrees to the top. The rest of the surface is left free. \revision{To avoid ambiguities in the rotation we use five steps of incremental loading in the Newton solver.}
We run the experiment on two sets of tetrahedral and hexahedral meshes. Coarse meshes  have $739$ and $750$ while the dense have $50\,000$ and $58\,719$ vertices respectively.
The first three images of Figure~\ref{fig:torsion} show that the results are indistinguishable, except for small numerical fluctuations. Similar results are for the plot (Figure~\ref{fig:torsion} last plot) of the rotation angle along a line starting at the point $[9.5, 9.5,0]$ parallel to the beam axis. Note that the dense $Q_2$ solution required more that 44GB of RAM for the solver, while $P_2$ required around 21GB.

%We omit the results for the $Q_2$ discretization because it exceeded our maximum budged of 30 hours running time.

The reason for the high memory consumption is the size of the element matrix, which has $(27\times 3)^ = 6\,561$ entries (compared to 144 for $P_1$, 576 for $Q_1$, and 900 for $P_2$). \revision{Note that the difference in running time does not come from the number of iterations of the Newton solver: for $P_1$, $Q_1$, $P_2$, $Q_2$ we obtained 16, 17, 20, 17 iterations respectively for the coarse mesh,  and 18, 16, 17, 18 for the dense mesh.}

We have repeated this experiment using quadratic B-spline bases on the coarse mesh. The result is similar to $P_2$ and $Q_2$, see the inset figure. For this particular example, we measure the \emph{solve} time of the three discretizations: the spline solve is 3 times faster than $P_2$ (0.51s versus 1.50s) and 9 times faster than $Q_2$ (0.51s vs 4.62s) while having roughly the same number of iterations: 16. Note that the assembly time (using full integration which  could be improved using~\cite{Schillinger:2014:RBE}) for spline is similar to $Q_2$ and is 12 times slower than $P_2$ (20.54s versus 1.70s). \revision{While using splines on regular grids is natural, the extension to irregular meshes requires the use of T- or U- splines~\cite{beer:2015:isogeometric,Wei:2018:BBS}, increasing the implementation complexity.}

\subsection{High Stress}\label{sec:high-stresses}

\begin{figure}\centering\footnotesize
    \parbox{\linewidth}{\centering
        \includegraphics[width=.22\linewidth]{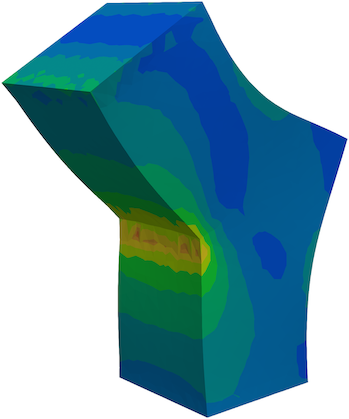}\hfill
        \includegraphics[width=.22\linewidth]{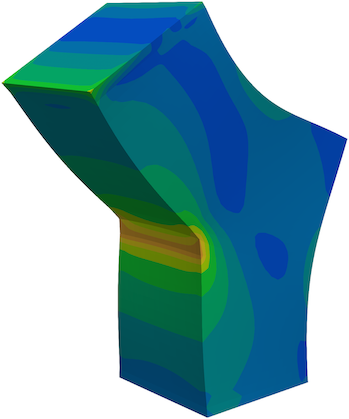}\hfill
        \includegraphics[width=.22\linewidth]{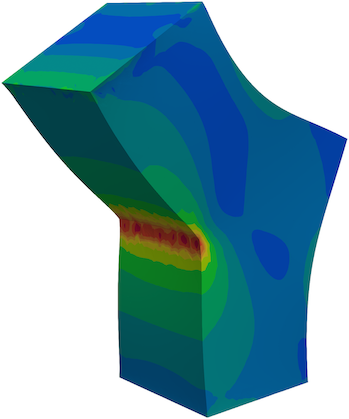}\hfill
        \includegraphics[width=.22\linewidth]{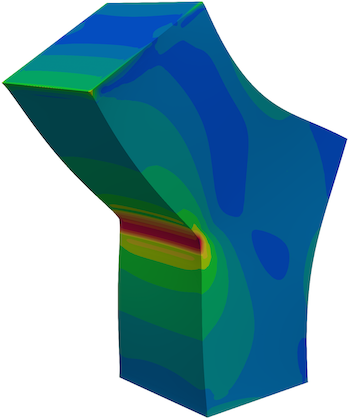}\hfill\hfill
        \includegraphics[width=.08\linewidth]{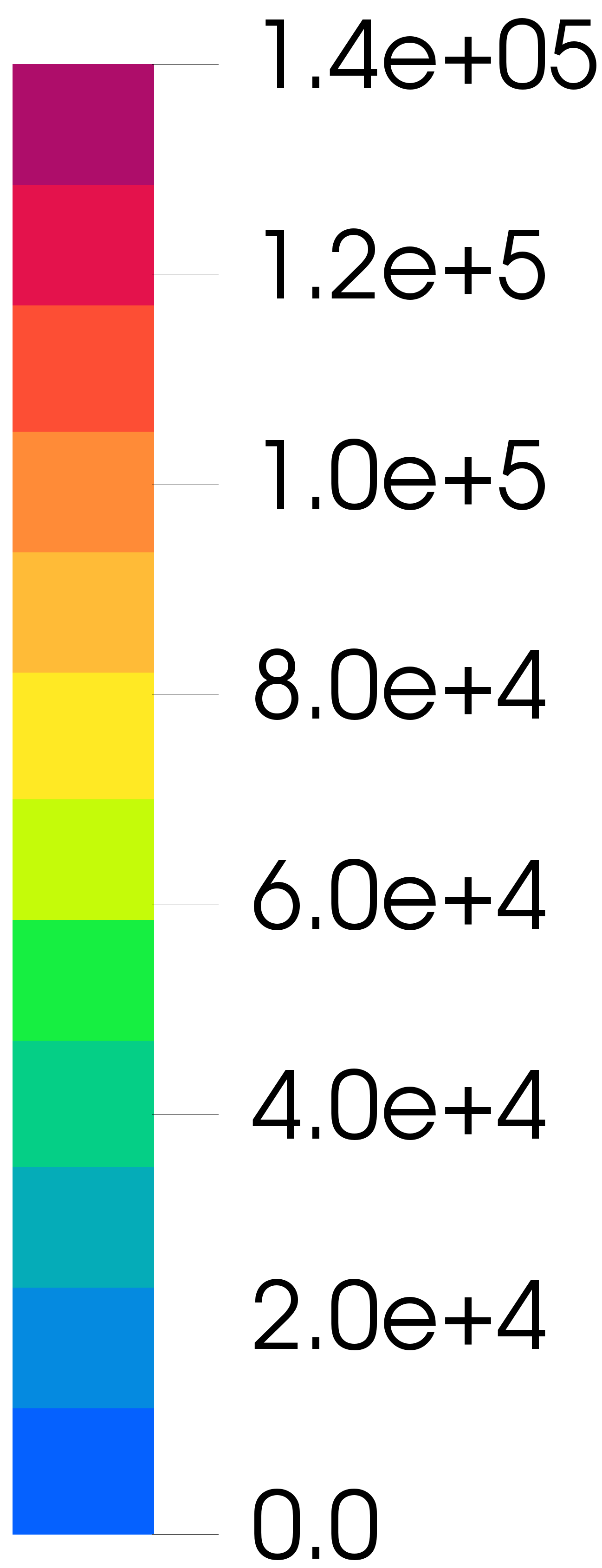}\par
        \parbox{.22\linewidth}{\centering $P_1$ 0:01:07}\hfill
        \parbox{.22\linewidth}{\centering $Q_1$ 0:03:00}\hfill
        \parbox{.22\linewidth}{\centering $P_2$ 0:17:39}\hfill
        \parbox{.22\linewidth}{\centering $Q_2$ 1:34:02}\hfill\hfill
        \parbox{.08\linewidth}{~}
    }
    \caption{Von Mises stress and singular solution timings for the four discretizations for the Neo-Hookean material model. The stresses are averaged around vertices (including the vertices where the solution is singular).}
    \label{fig:L}
\end{figure}

As a final experiment, we run a simulation for an L-shaped domain with the Neo-Hookean material and
$E=210000$ $\nu=0.3$. Our goal is to study the differences in the stresses for singular solutions: the concave corner of L will have a stress singularity. We mesh our domains with $14\,155$ vertices for the tetrahedral mesh and $14\,161$ vertices for the hexahedral mesh. We fixed the bottom part of the domain (zero displacement) and rotate the top part by 120 degrees (Dirichlet constraint on the displacement), the rest of the boundary is let free (zero Neumann condition). Figure~\ref{fig:L} shows that linear tetrahedral elements underestimate the stress while linear hexahedral elements are somewhat better.
The quadratic discretizations are qualitatively similar: the hexahedral mesh does not have the spurious small stress oscillations of $P_2$
because the elements are aligned with the mesh, however the price to pay is significant, 17 minutes for $P_2$ compared to more than 1.5 hours for $Q_2$.

\section{Large Dataset}\label{sec:large-dataset}

Next, to evaluate the performance of different types of elements on a large diverse set of realistic domains, we compute  solutions for the Poisson~\eqref{eq:poisson} and linear elasticity~\eqref{eq:lin-elast}.  We use the \emph{method of manufactured solutions}~\cite{SALARI:2000:CVB}, that is, for an analytically defined  solution $u$ we compute the corresponding  right-hand side $b$ by plugging it into the PDE. The boundary condition $d$ is obtained by sampling $u$ on the boundary.
For the Poisson equation, we use the Franke function~\cite{Franke:1979:ACC,Cavoretto:2017:OPUM}
\begin{align*}
    u & (x_1,x_2,x_3) =
    \\&3/4\,{ e}^{-((9x_1-2)^2+(9x_2-2)^2+(9x_3-2)^2)/4}
    \\&+3/4\,{ e}^{-(9x_1+1)^2/49-(9x_2+1)/10-(9x_3+1)/10}
    \\&+1/2\,{ e}^{-((9x_1-7)^2+(9x_2-3)^2+(9x_3-5)^2)/4}
    \\&-1/5\,{ e}^{-(9x_1-4)^2-(9x_2-7)^2-(9x_3-5)^2},
\end{align*}
while for elasticity
\[
    u(x_1,x_2,x_3) =
    \frac{1}{80}\begin{pmatrix}
        x_1x_2+x_1^2+x_2^3+6x_3      \\
        x_1x_3-x_3^3+x_1x_2^2+3x_1^4 \\
        x_1x_2x_3+x_2^2x_3^2-2x_1,
    \end{pmatrix}
\]
with Lam\`e parameters $E=200$ and $\nu = 0.35$. In addition to standard tensor product bases for hexahedra, we compare to the popular \emph{serendipity} bases~\cite{Zienkiewicz:2005:TFE}[Chapter 6], which have only 20 nodes per element instead of 27.

We use two sources for our data: (1) the Hexalab dataset containing results of 16 state-of-the-art hexahedral meshing techniques~\cite{Bracci:2019:HNA}, (2) the Thingi10k dataset~\cite{Zhou:2016:TKA} consisting of triangulated surfaces. For each dataset, we produce a tetrahedral mesh dataset from the surfaces of the hexahedral meshes \revisionn{(generated with MeshGems~\cite{code:meshgems})} using TetWild~\cite{Hu:2018:TMI} with a matching number of vertices. Note that since matching the number of vertices is a heuristic process, we discard all meshes where the difference in the number of vertices is larger than 5\% of the total number of vertices.
To ensure that we are solving a similar problem on the two tessellations we remove meshes whose Hausdorff distance between the surfaces of corresponding hexahedral and tetrahedral meshes differs more than $10^{-3}$ of the diagonal of the bounding box of the hexahedral mesh surfaces. Finally, we discard all meshes whose ratio between boundary and total vertices is more than 30\%. \revision{Since the Hexalab dataset is small, we opted for doing one step of uniform refinement to increase the number of interior vertices instead of discarding them}.
In summary, the two datasets are:
\begin{enumerate}
    \item 237 Hexalab hexahedral meshes and 237 tetrahedral meshes generated with Tetwild.
    \item 3\,200 hexahedral meshes generated with MeshGems~\cite{code:meshgems} and 3\,200 tetrahedral meshes generated with Tetwild both obtained from the surfaces in the Thingi10k dataset.
\end{enumerate}

\begin{figure}\centering\footnotesize
    \parbox{\linewidth}{\centering
        \rotatebox{90}{\centering Percentage}
        \parbox{0.95\linewidth}{
            \parbox{.45\linewidth}{\centering
                Max
                \includegraphics[width=\linewidth]{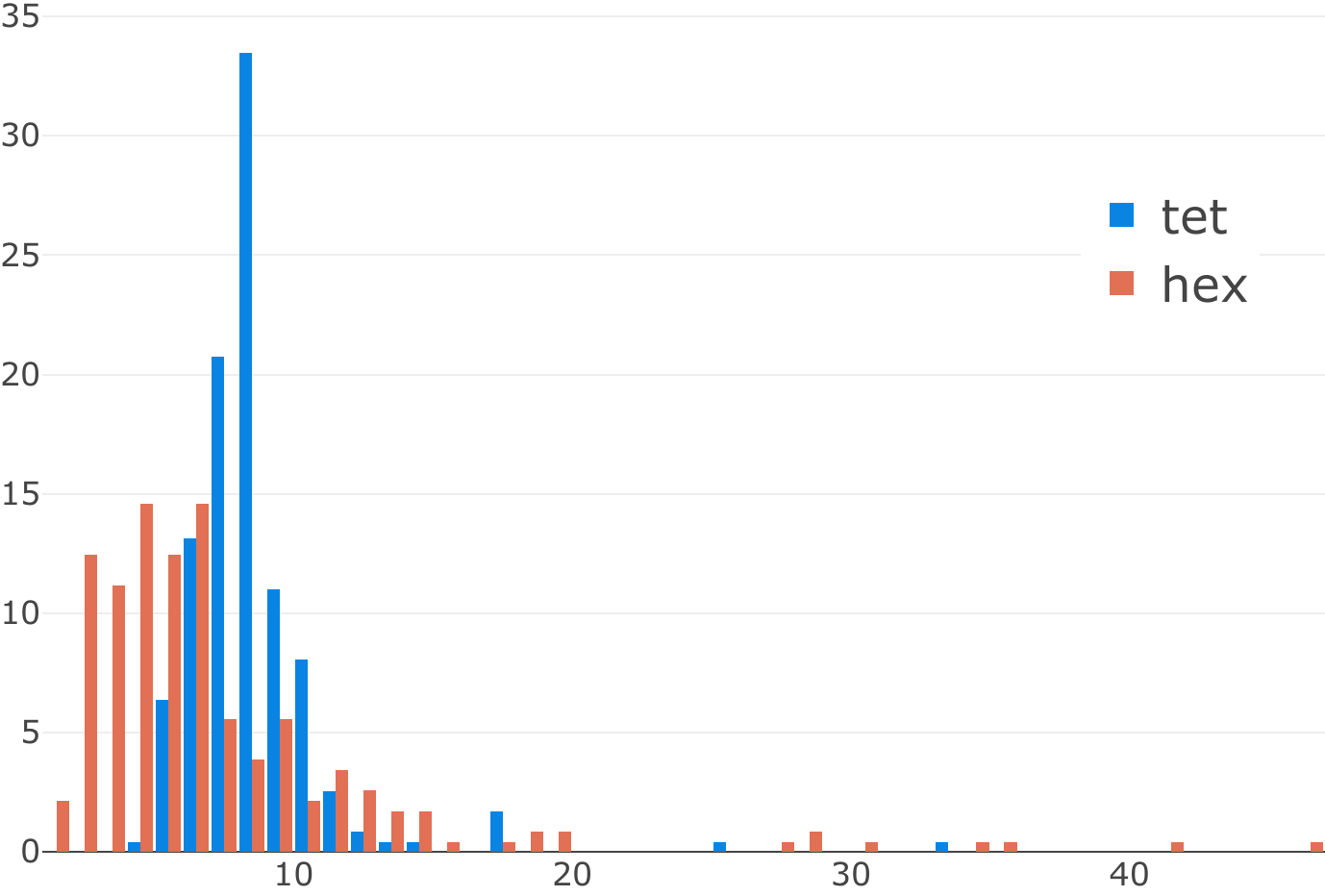}\par
                Avg
                \includegraphics[width=\linewidth]{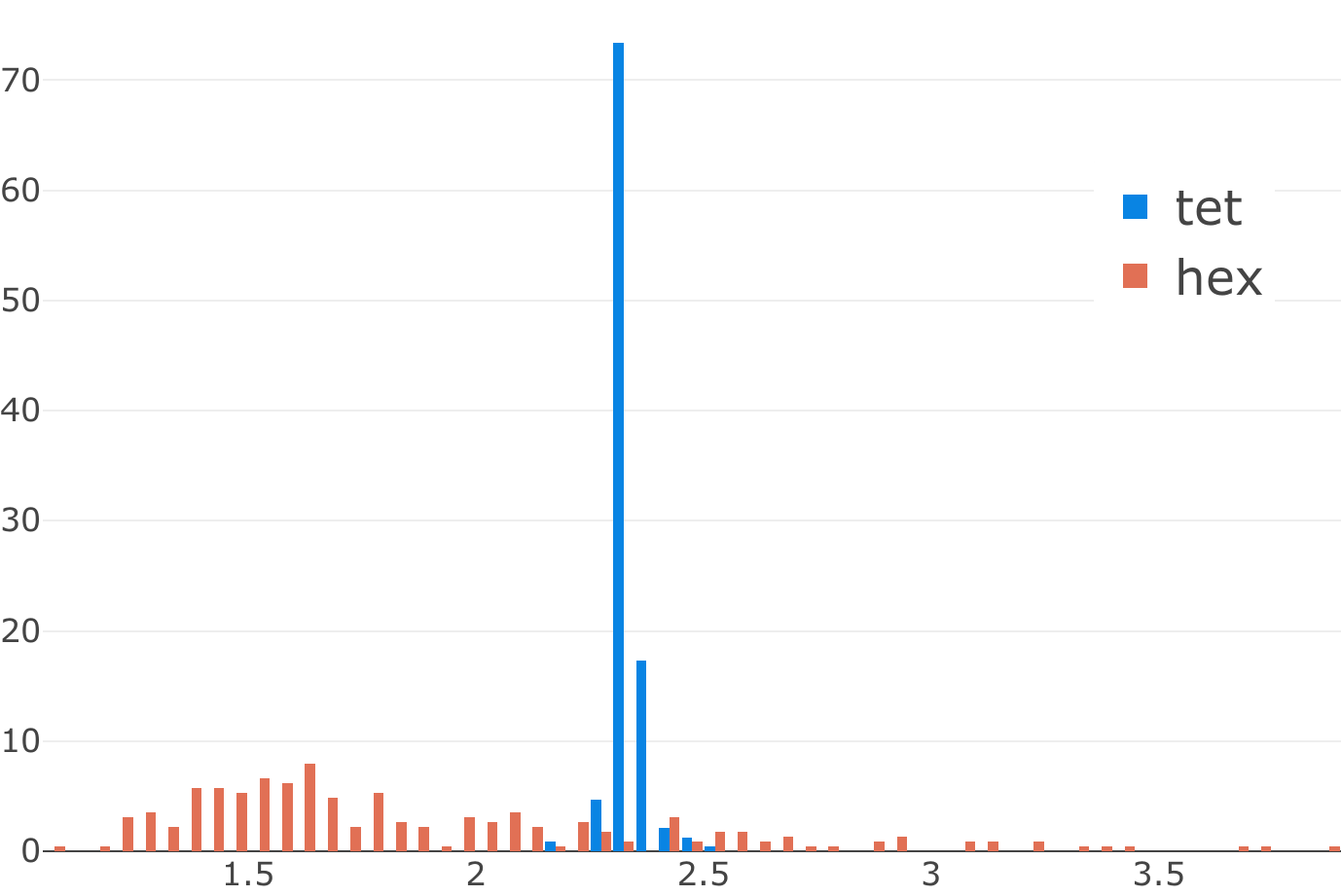}\par}\hfill
            \parbox{.45\linewidth}{\centering
                Max
                \includegraphics[width=\linewidth]{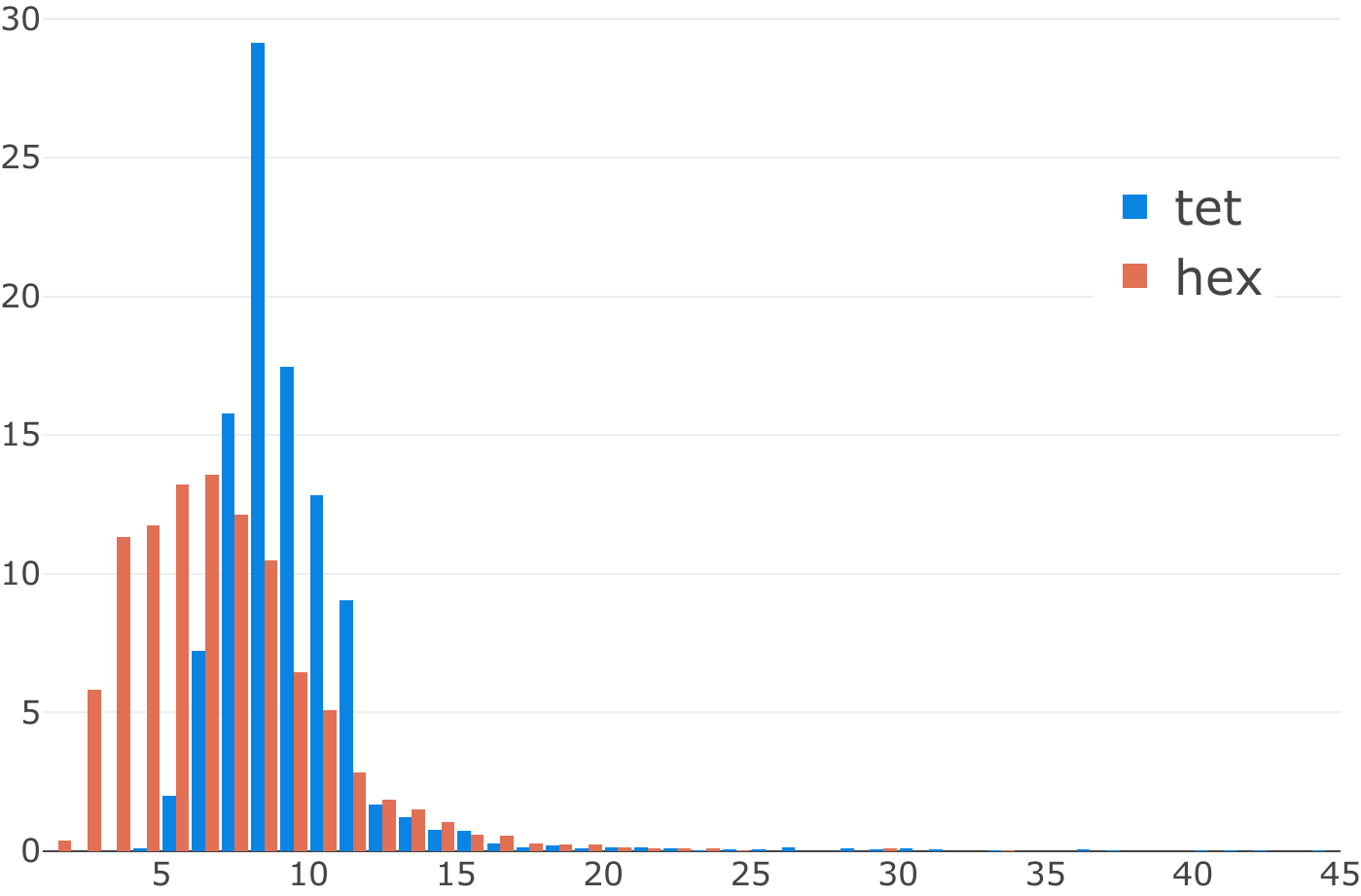}\par
                Avg\\
                \includegraphics[width=\linewidth]{figs/aspect_hist/10k_max}\par}\par%\vspace{-0.7em}
            \parbox{.45\linewidth}{\centering Hexalab}\hfill
            \parbox{.45\linewidth}{\centering Thingi10k}
        }\\
        Aspect ratio
    }
    \caption{\revision{Histogram of the maximal and mean aspect ratios for the Hexalab (left) and Thingi10k (right) datasets.}}
    \label{fig:dataset-aspect}
\end{figure}

\begin{figure}\centering\footnotesize
    \parbox{\linewidth}{\centering
        \includegraphics[width=.49\linewidth]{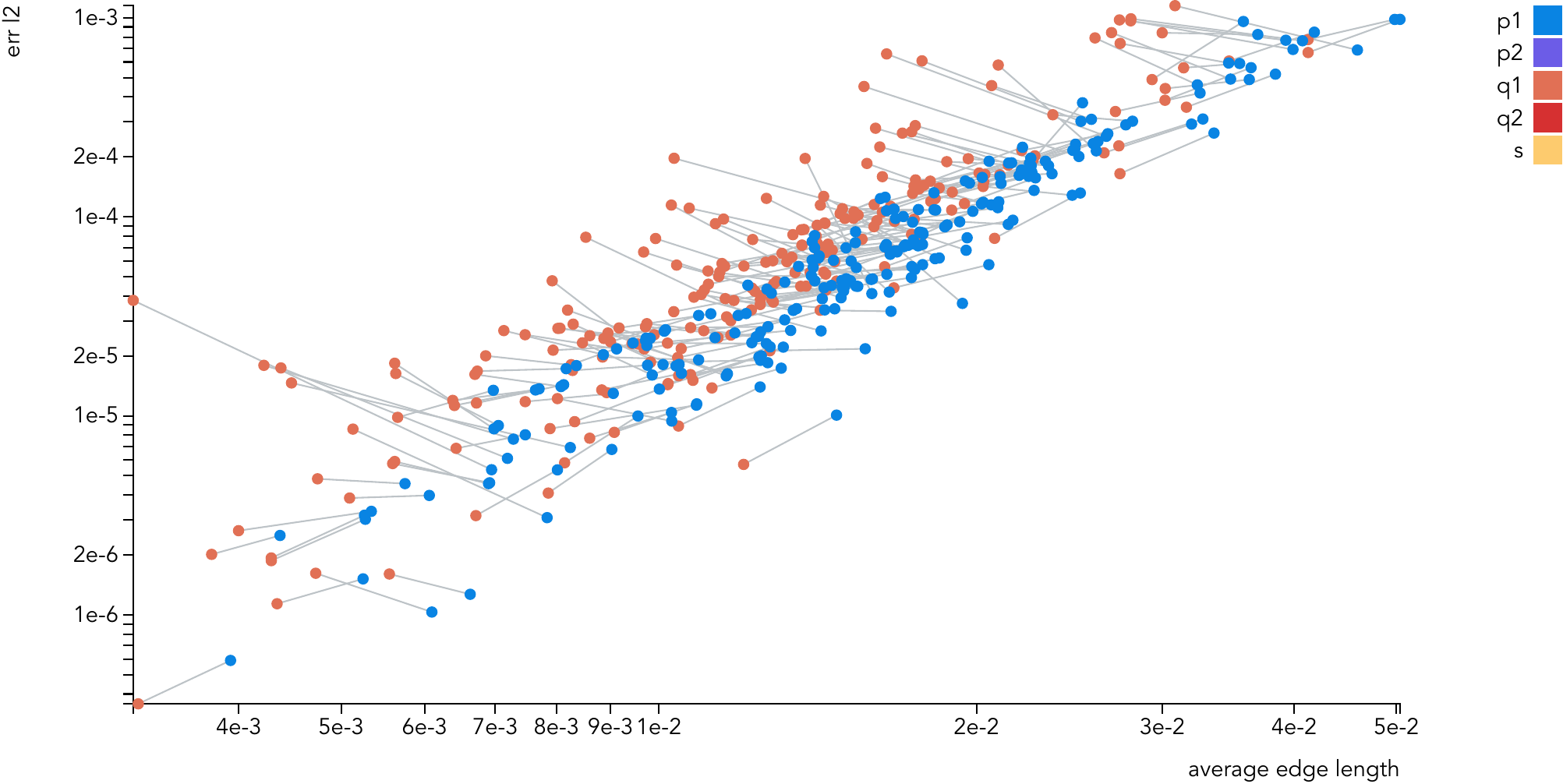}\hfill
        \includegraphics[width=.49\linewidth]{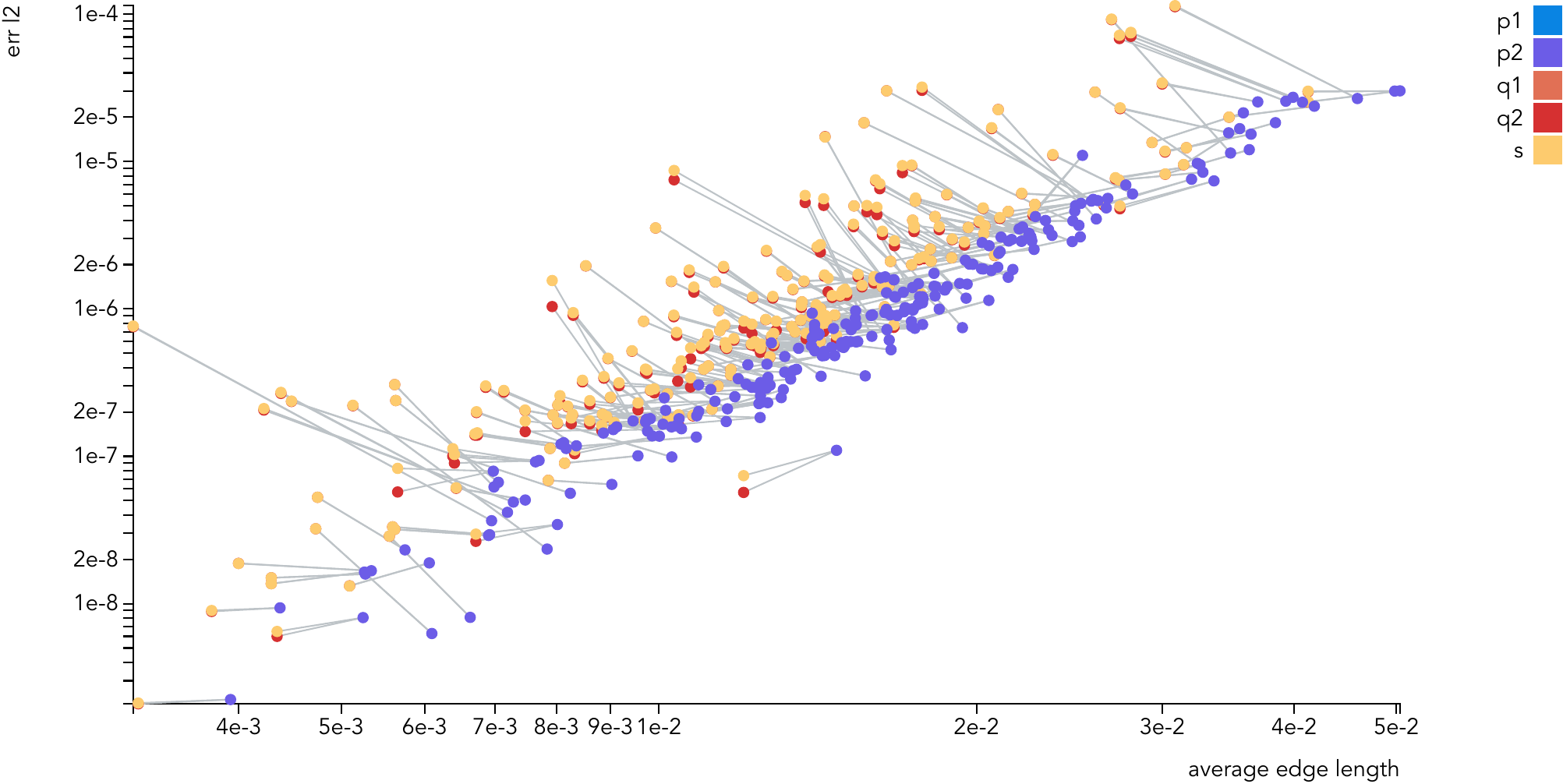}\par
        Poisson equation~\eqref{eq:poisson}\\[4ex]
        %%%%%%%%%%%%%%%%%%%%%%%%%%%%%%
        \includegraphics[width=.49\linewidth]{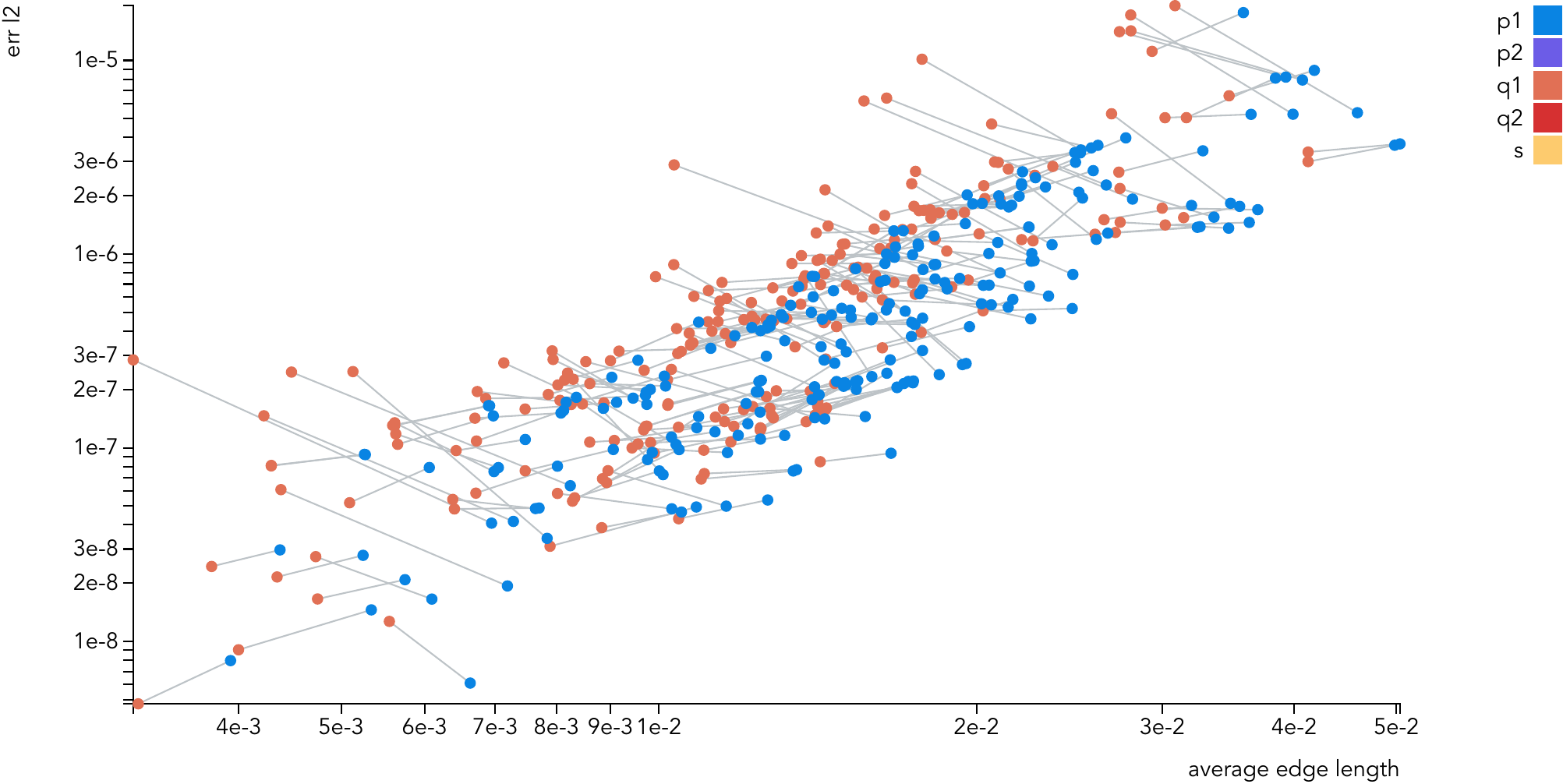}\hfill
        \includegraphics[width=.49\linewidth]{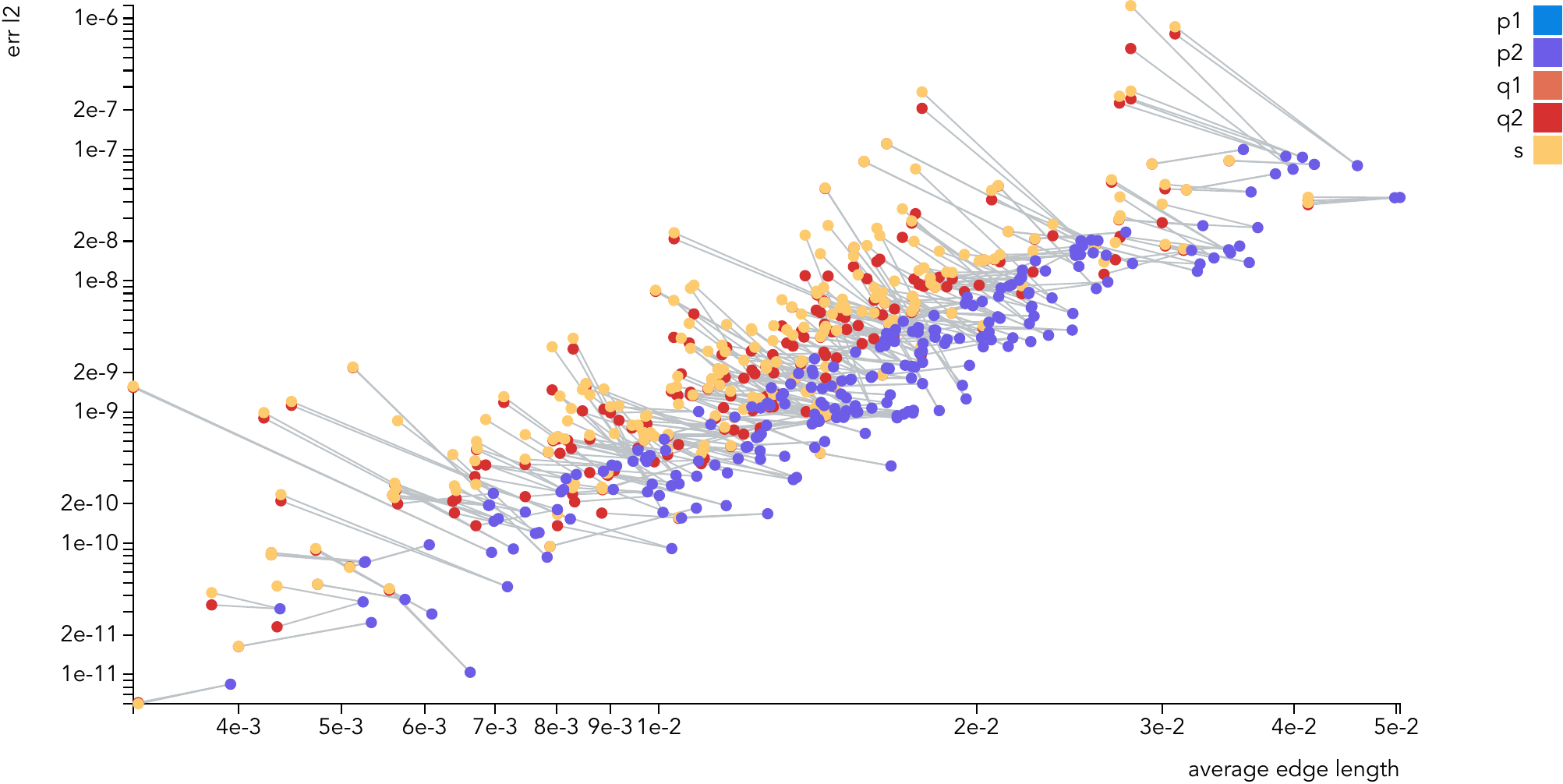}\par
        linear elasticity~\eqref{eq:lin-elast}
        %%%%%%%%%%%%%%%%%%%%%%%%%%%%%%
    }
    \caption{$L_2$ error vs. average mesh size \revisionn{for the Hexalab dataset} for linear (left) and quadratic (right) elements. The lines connect two points belonging to the same model. }
    \label{fig:hexalb_err_vs_h}
\end{figure}

\begin{figure}\centering\footnotesize
    \parbox{\linewidth}{\centering
        \includegraphics[width=.49\linewidth]{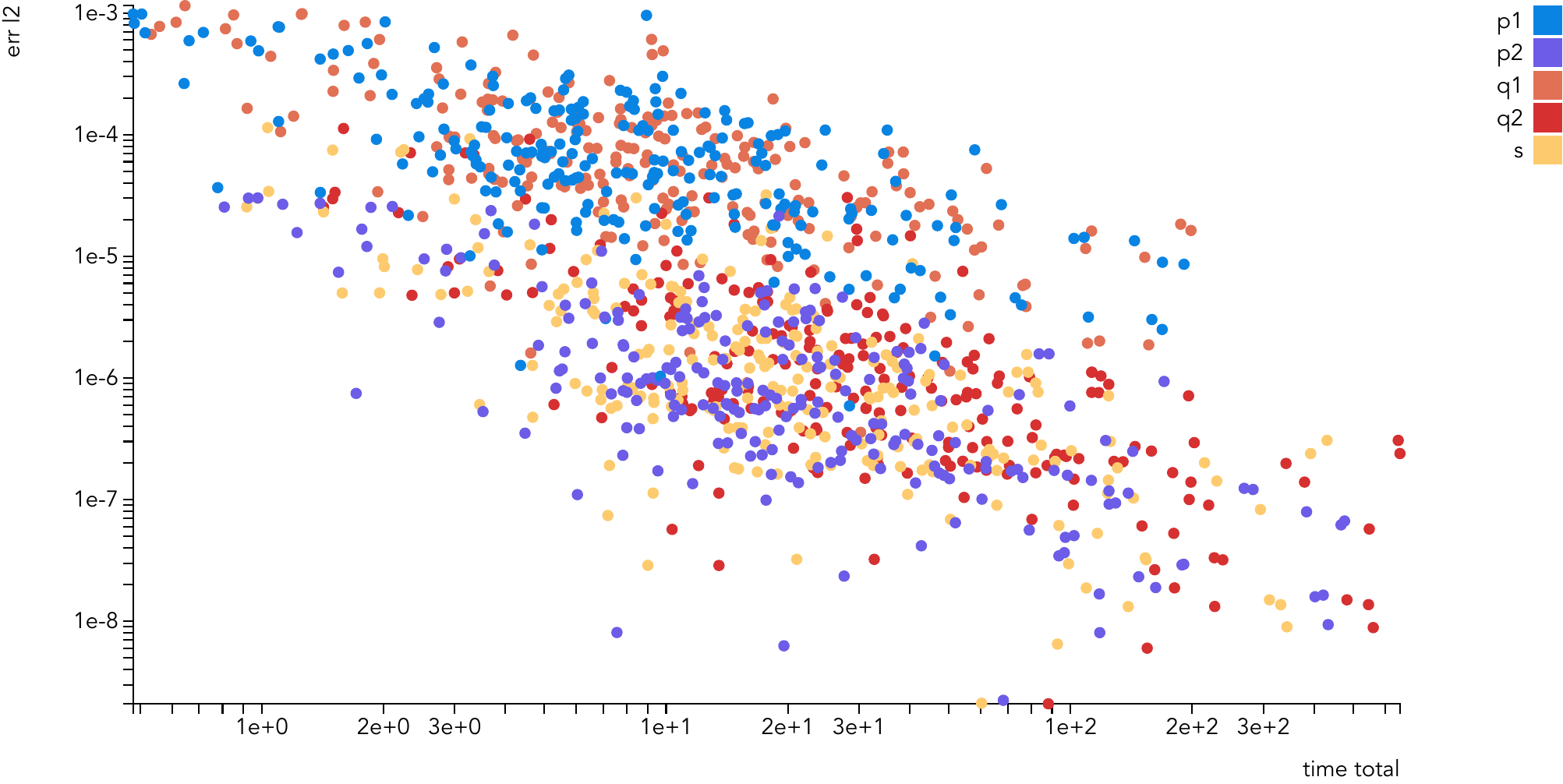}\hfill
        \includegraphics[width=.49\linewidth]{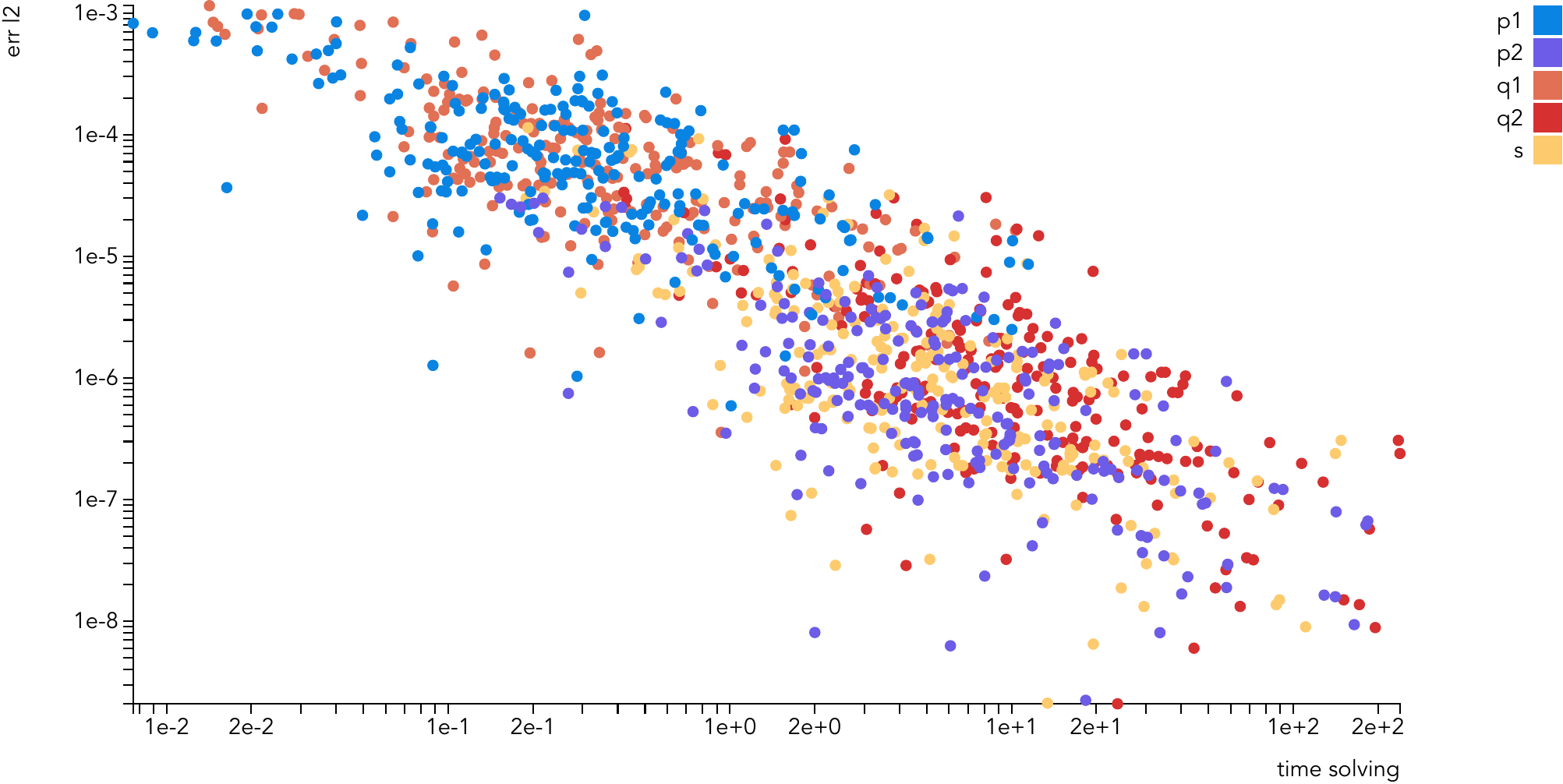}\par
        Poisson equation~\eqref{eq:poisson}\\[4ex]
        %%%%%%%%%%%%%%%%%%%%%%%%%%%%%%
        \includegraphics[width=.49\linewidth]{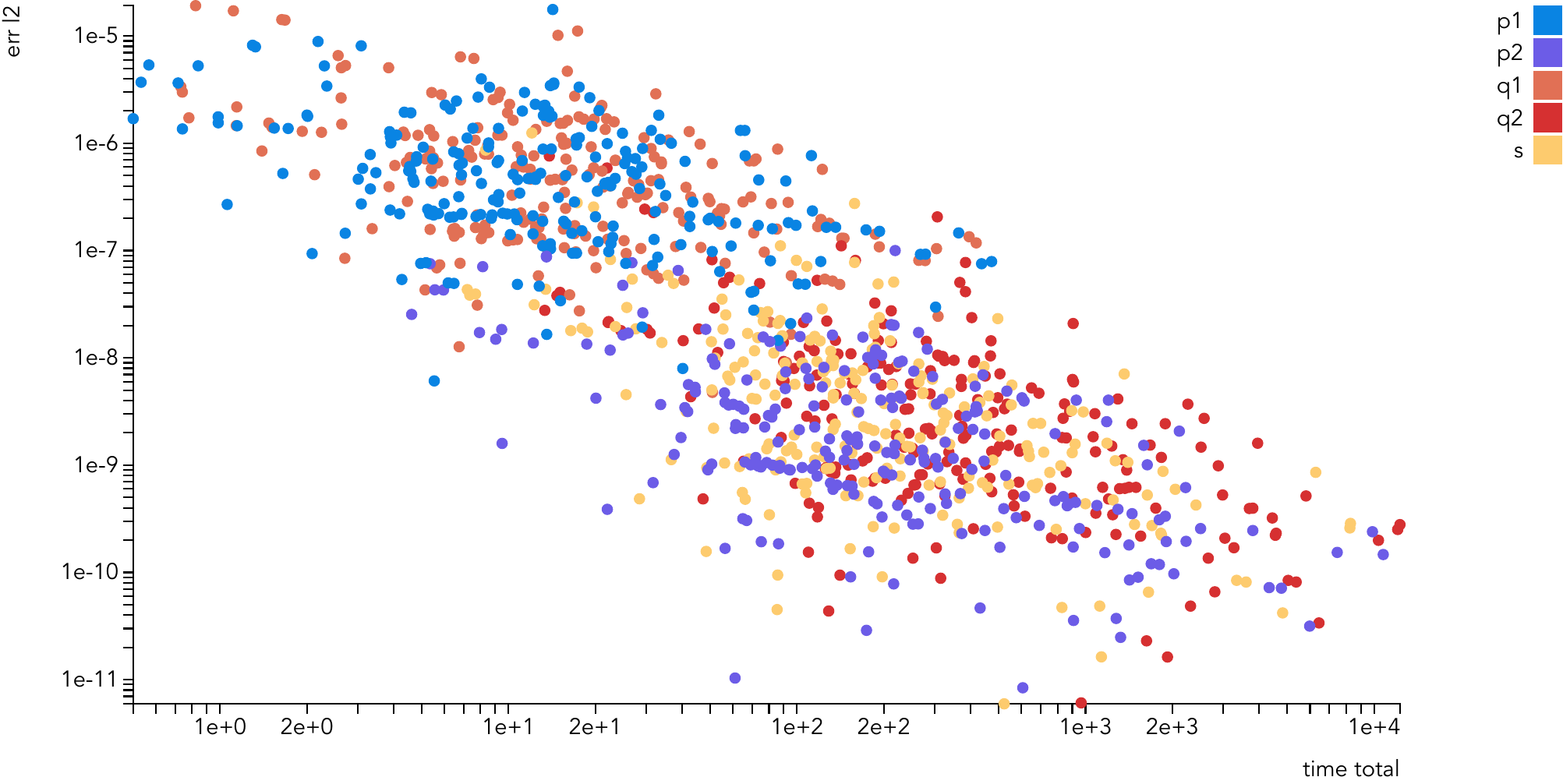}\hfill
        \includegraphics[width=.49\linewidth]{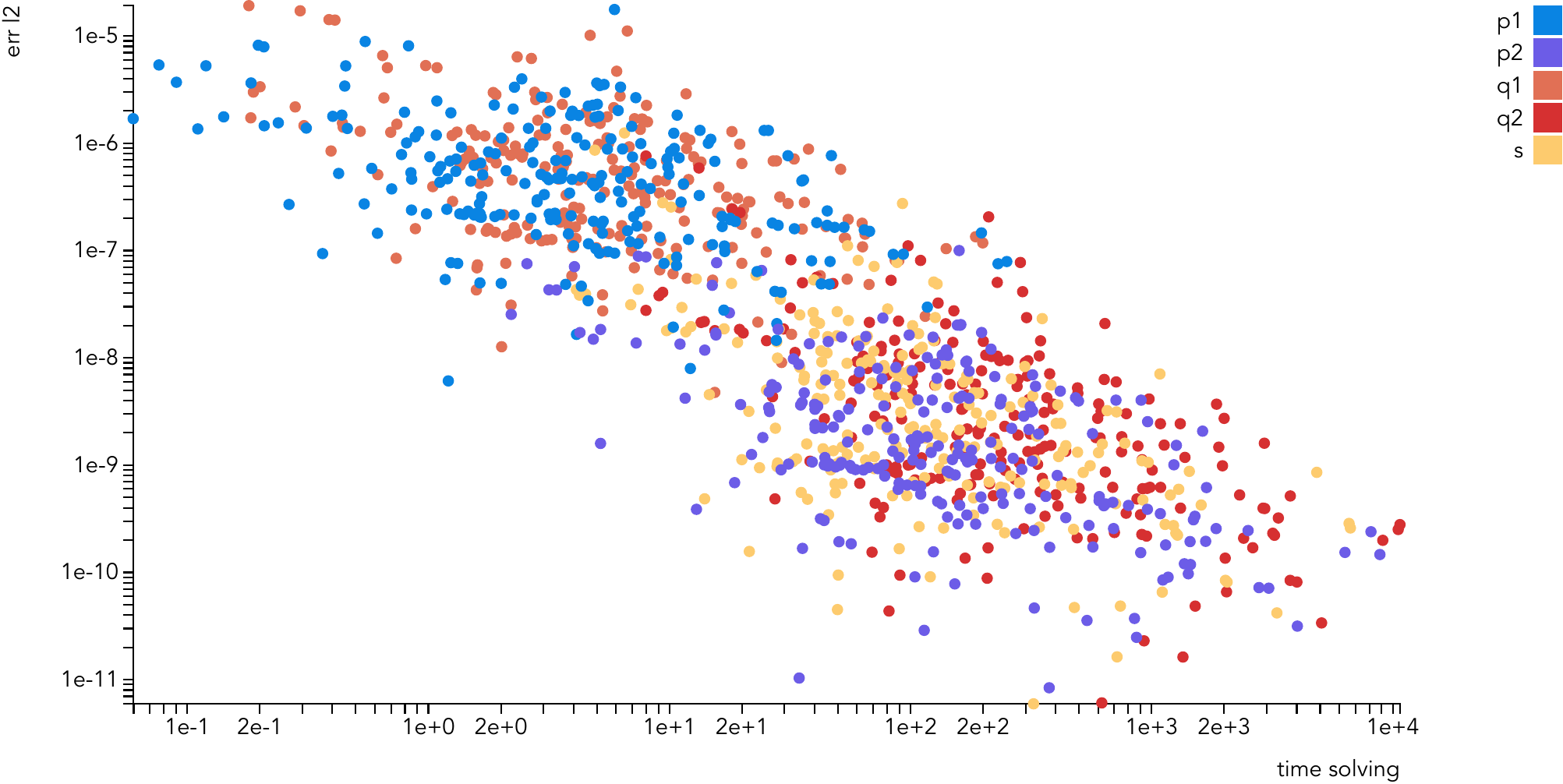}\par
        Linear elasticity~\eqref{eq:lin-elast}
        %%%%%%%%%%%%%%%%%%%%%%%%%%%%%%
    }
    \caption{$L_2$ error vs. total time (left) and solver time (right) \revisionn{for the Hexalab dataset}.}
    \label{fig:hexalb_time}
\end{figure}

For conciseness, we report only the most significant results. Many other metrics (e.g., $H^1$-error, the time required to assemble bases, nonzero entries of the matrix, etc.) can be found in the \interactiveplot{}.

We remark that, while Tetwild \emph{guarantees} to produce valid tetrahedral meshes, Meshgems and the methods used in the Hexalab dataset do not provide any guarantee. We observe that out of the 237 Hexalab meshes, 8 (3.4\%)  contain at least one inverted element (2 from \cite{Livesu:2016:SDA} and 6 from \cite{Gao:2016:SVD}). For the Thingi10k dataset, Meshgems produces 577 (18.0\%) meshes with at least one invalid element. To check if a hexahedron has a negative volume we sample it with $10^3$ uniformly spaced samples, evaluate the Jacobian at each point, and mark it as flipped if at least one evaluation is negative. \revision{Another important quality measure is the aspect ratio of the elements (Section~\ref{sec:high-aspect}). Figure~\ref{fig:dataset-aspect} shows that both our datasets contain reasonably well-shaped elements.}

All experiments are run on a cluster node with 2 Xeon E5-2690v4 2.6GHz CPUs and 250GB memory, each with max 128GB of reserved memory and 8 threads. For all experiments we use the Hypre~\cite{code:hypre} algebraic multigrid iterative solver and the PolyFEM library for the finite element assembly.

\begin{figure}\centering\footnotesize
    \parbox{\linewidth}{\centering
        \includegraphics[width=.49\linewidth]{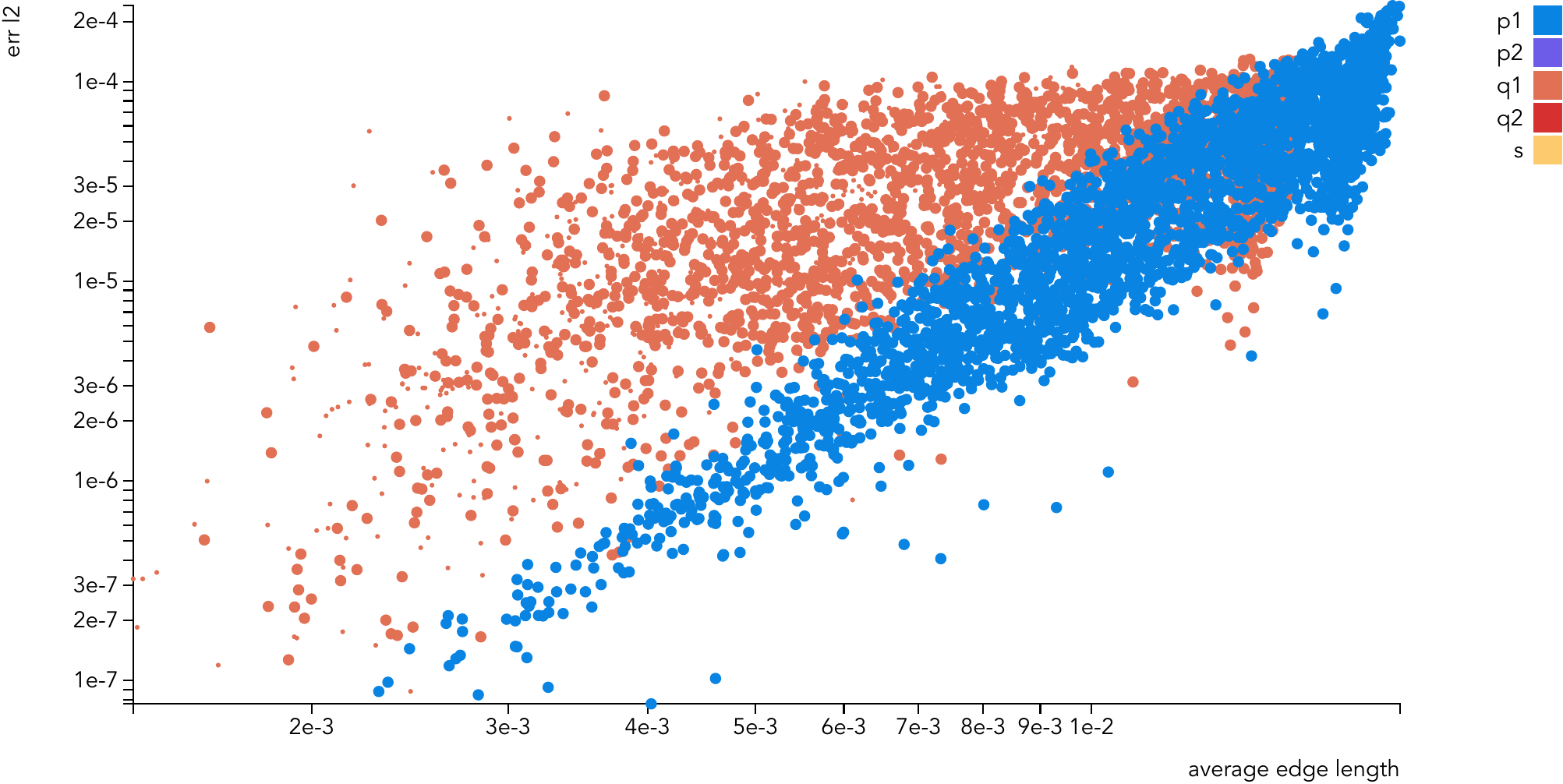}\hfill
        \includegraphics[width=.49\linewidth]{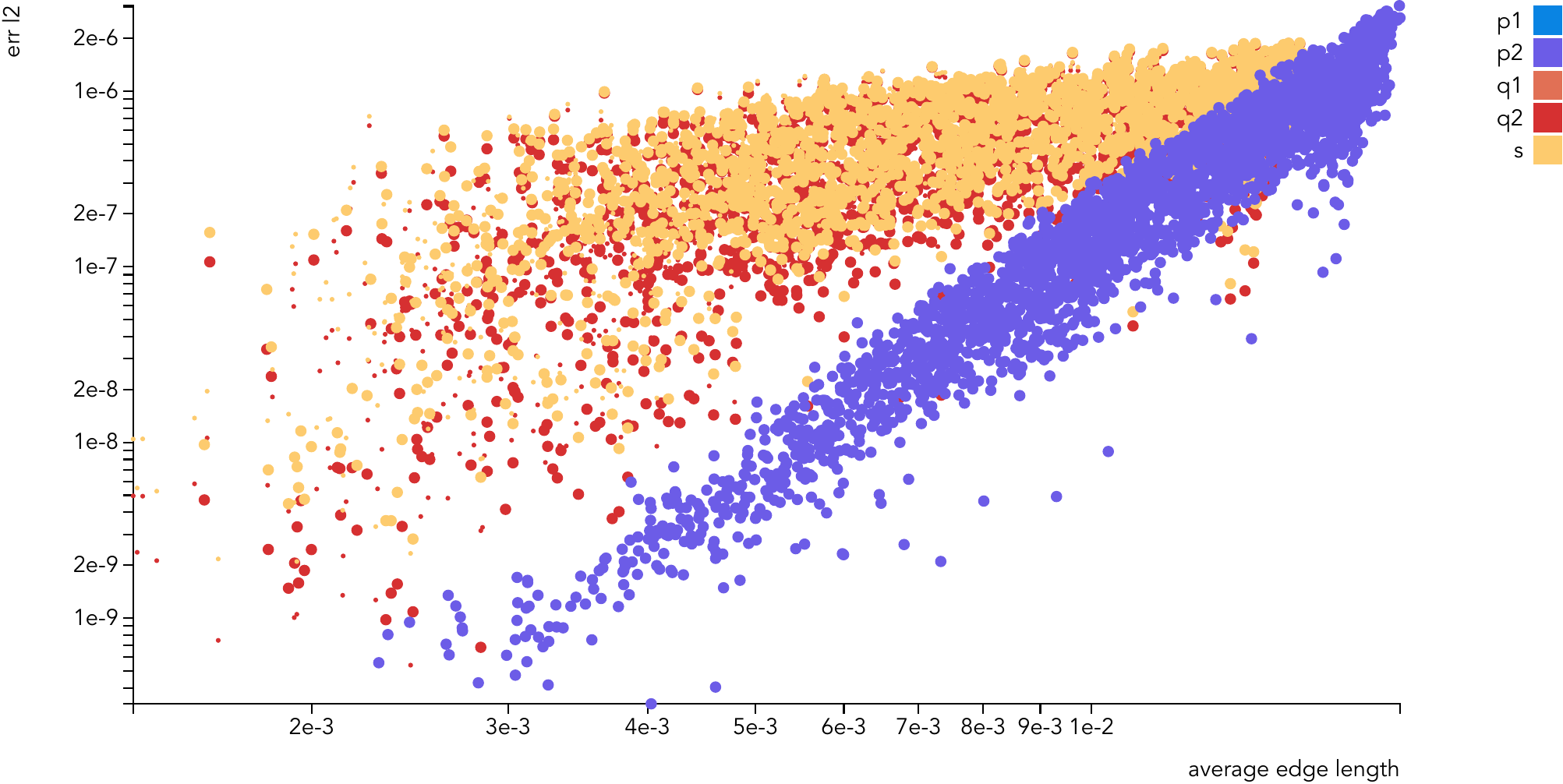}\par
        Poisson equation~\eqref{eq:poisson}\\[4ex]
        %%%%%%%%%%%%%%%%%%%%%%%%%%%%%%
        \includegraphics[width=.49\linewidth]{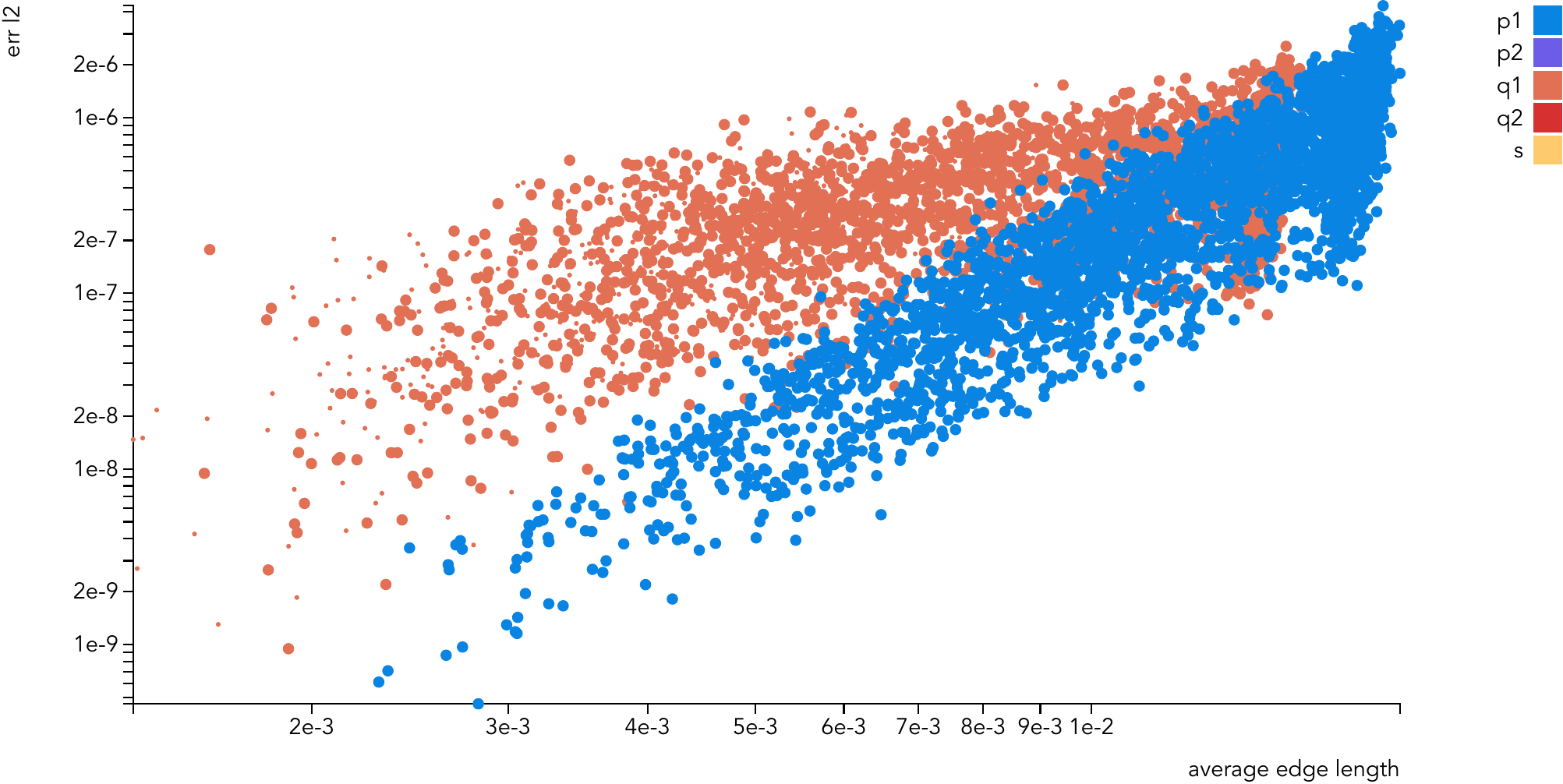}\hfill
        \includegraphics[width=.49\linewidth]{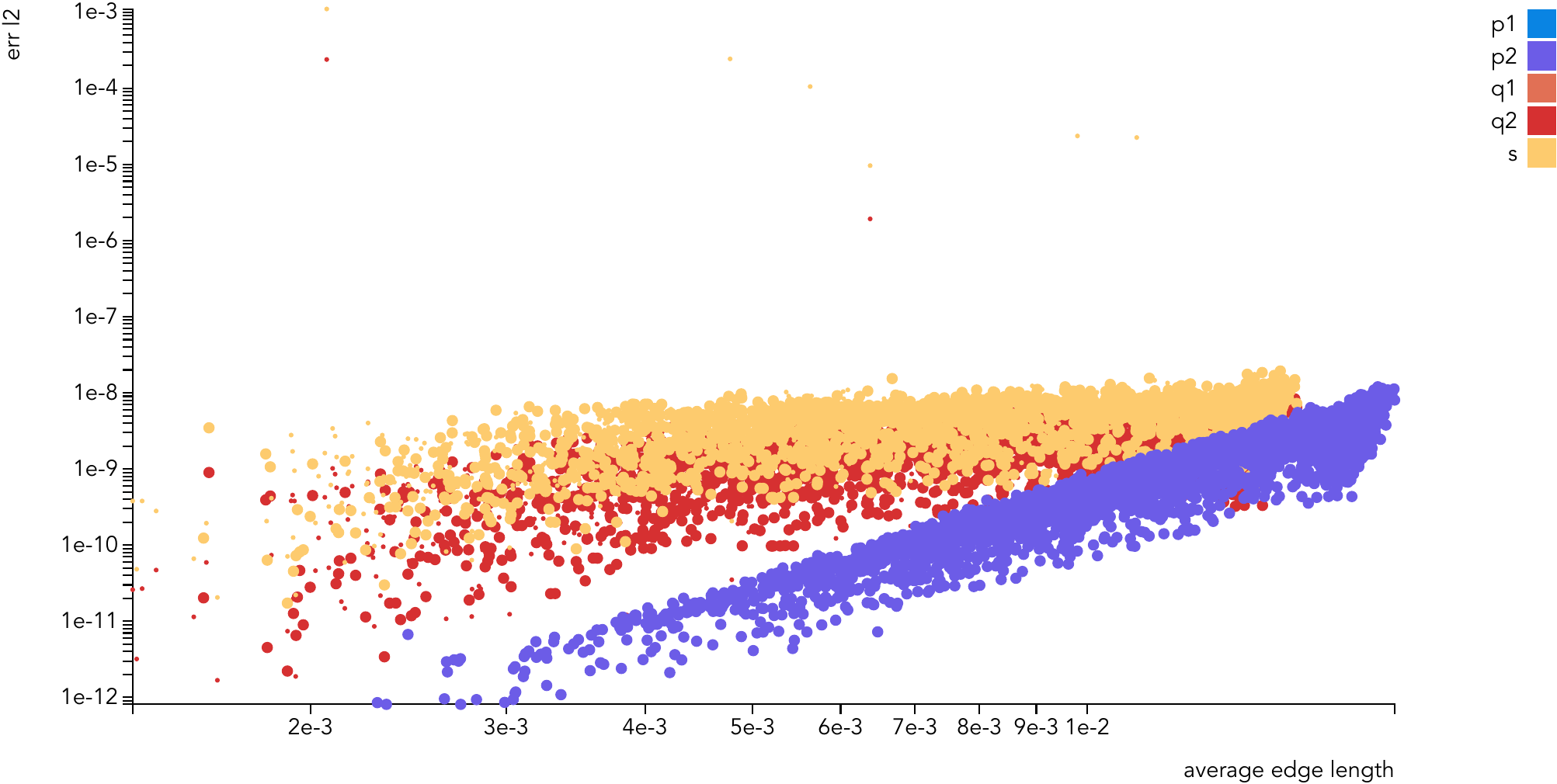}\par
        linear elasticity~\eqref{eq:lin-elast}
        %%%%%%%%%%%%%%%%%%%%%%%%%%%%%%
    }
    \caption{$L_2$ error vs. average mesh size for linear (left) and quadratic (right) elements. The smaller dot sizes indicate models with inverted elements. }
    \label{fig:thingi10k_err_vs_h}
\end{figure}

\begin{figure}\centering\footnotesize
    \parbox{\linewidth}{\centering
        \includegraphics[width=.49\linewidth]{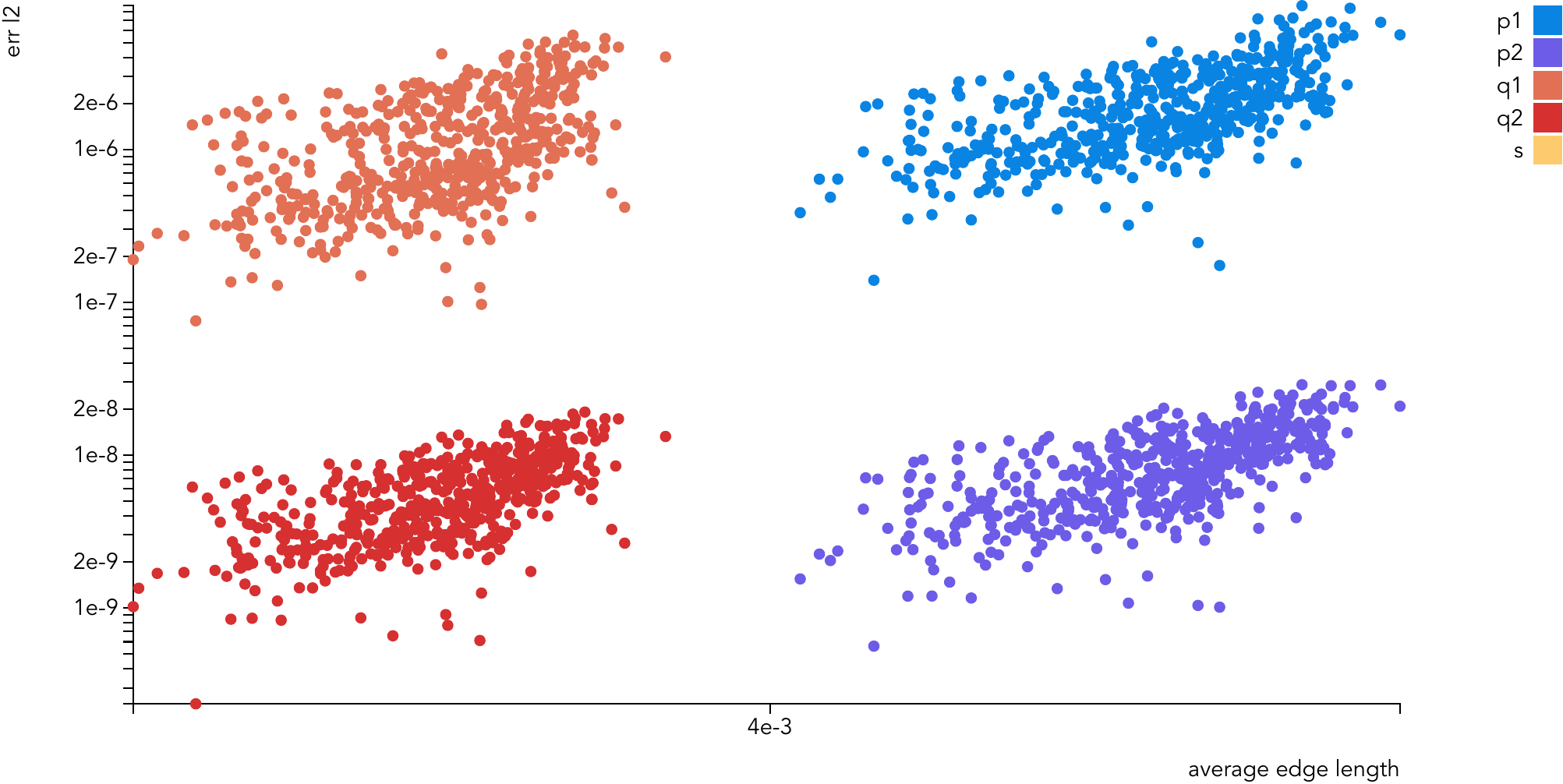}\hfill
        \includegraphics[width=.49\linewidth]{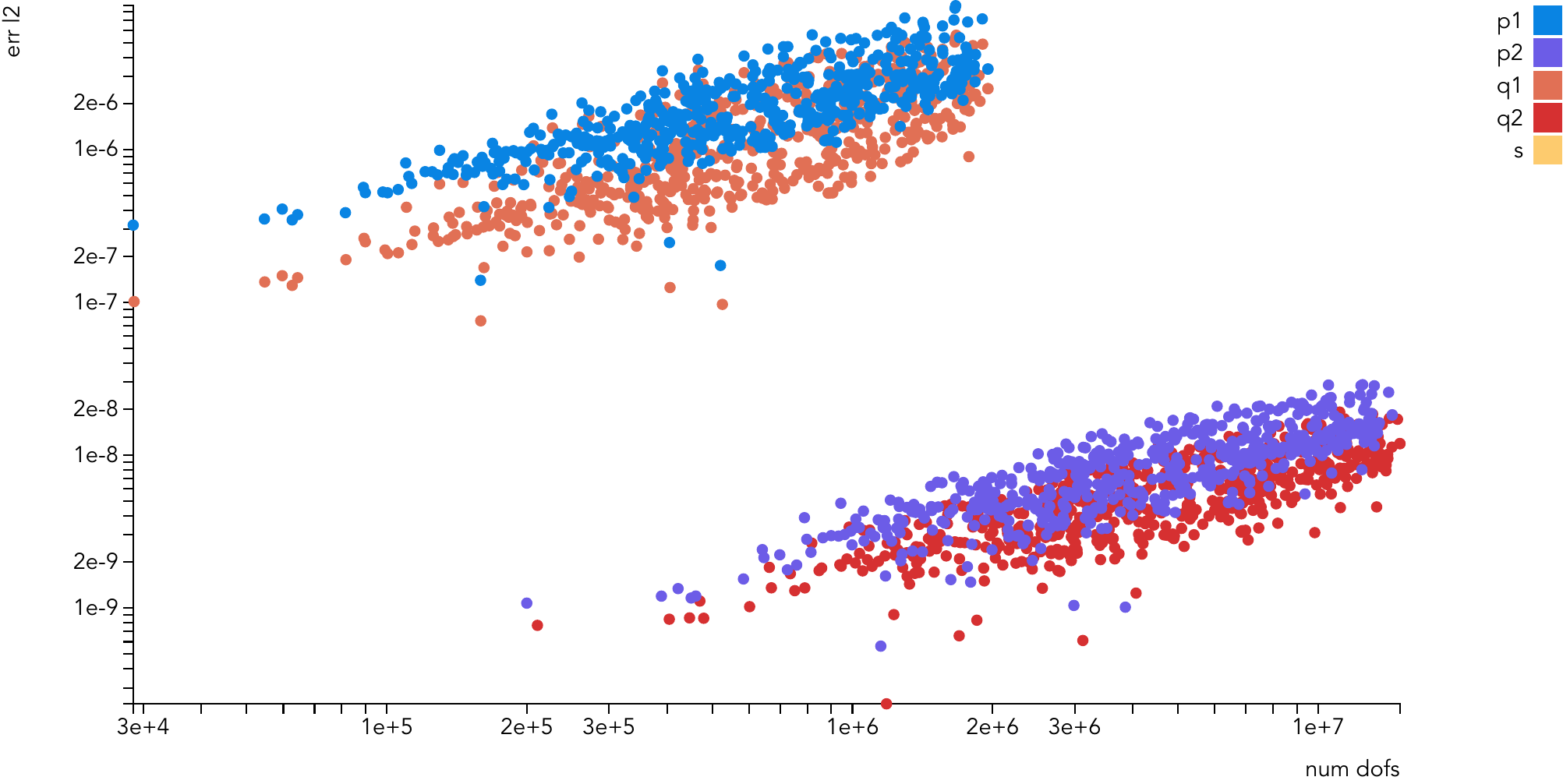}\par
    }
    \caption{The $L_2$ error  vs. average edge length (left) and number of degrees of freedom (right)  for the Poisson equation~\eqref{eq:poisson} on 580 ``uniform'' hexahedral meshes.}
    \label{fig:uniform}
\end{figure}

\begin{figure}\centering\footnotesize
    \parbox{\linewidth}{\centering
        \includegraphics[width=.49\linewidth]{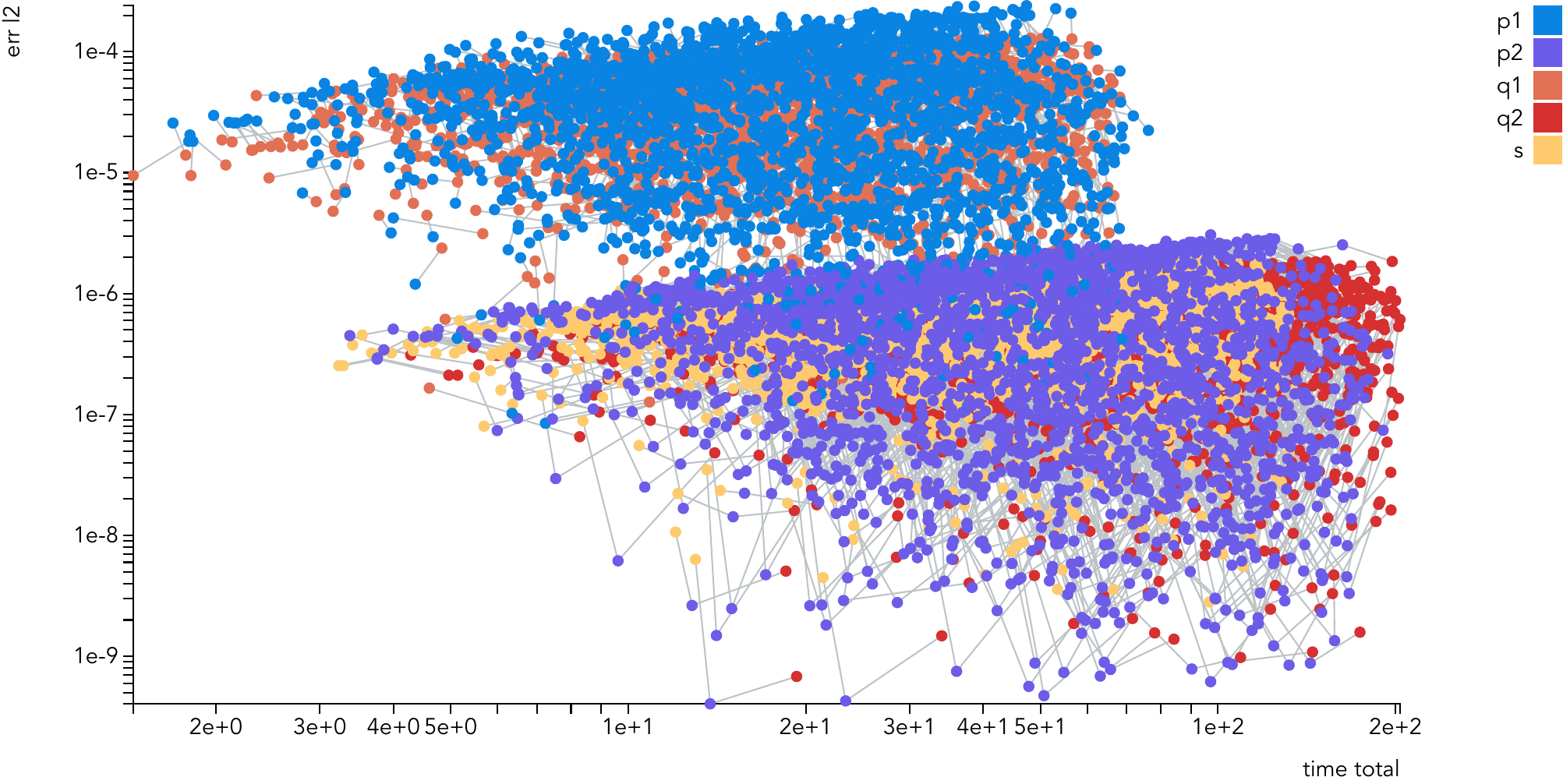}\hfill
        \includegraphics[width=.49\linewidth]{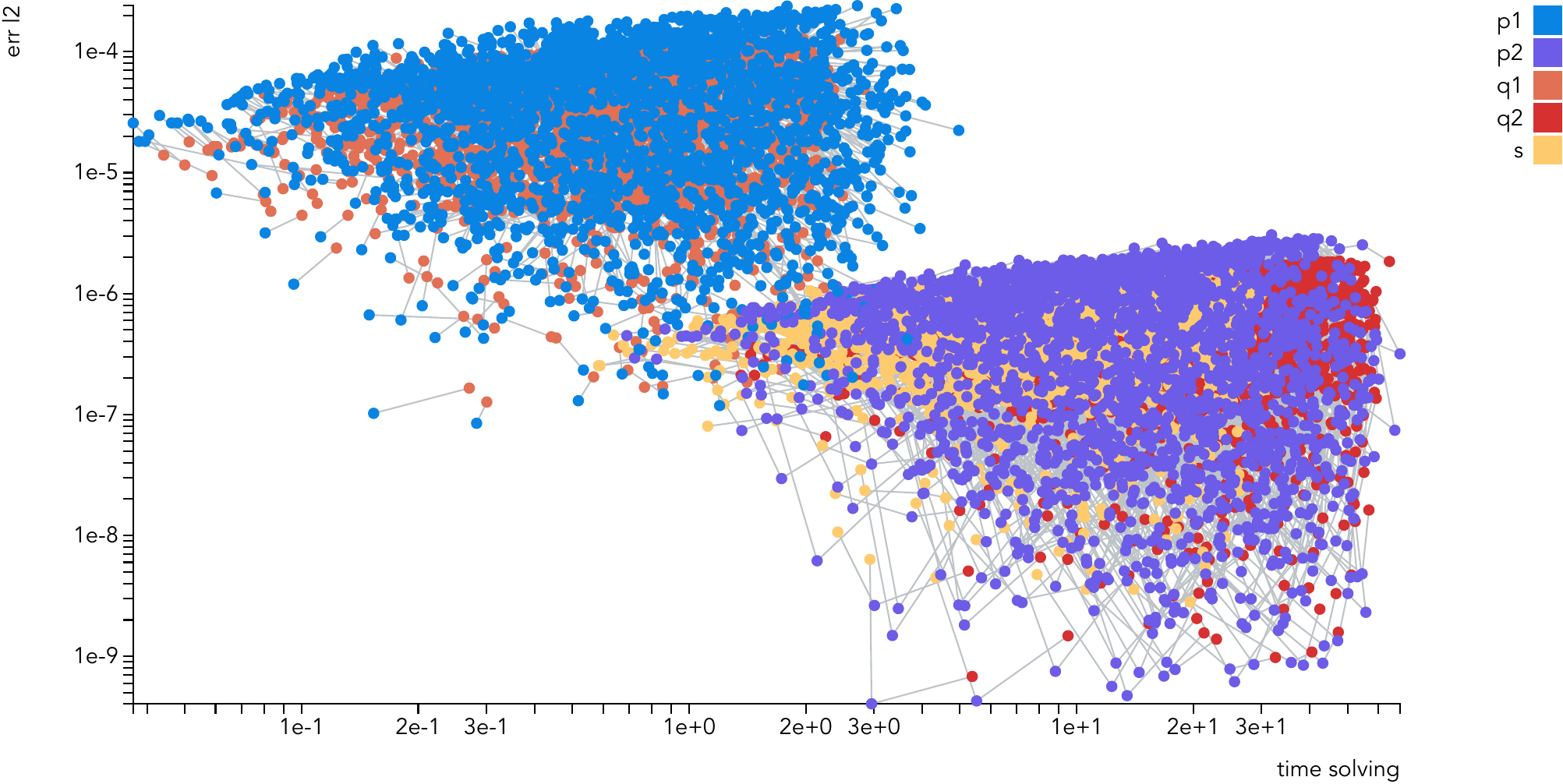}\par
        Poisson equation~\eqref{eq:poisson}\\[4ex]
        %%%%%%%%%%%%%%%%%%%%%%%%%%%%%%
        \includegraphics[width=.49\linewidth]{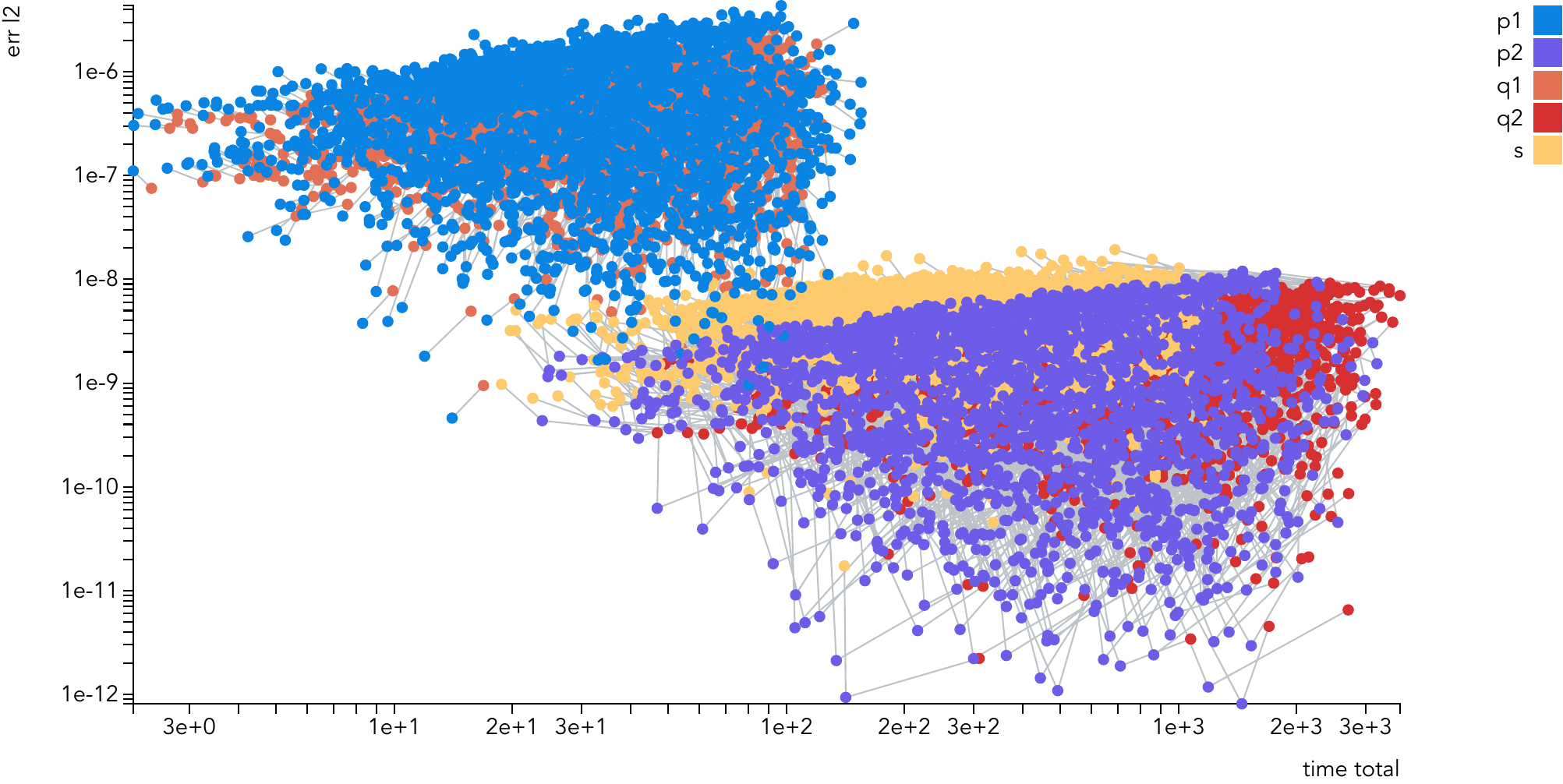}\hfill
        \includegraphics[width=.49\linewidth]{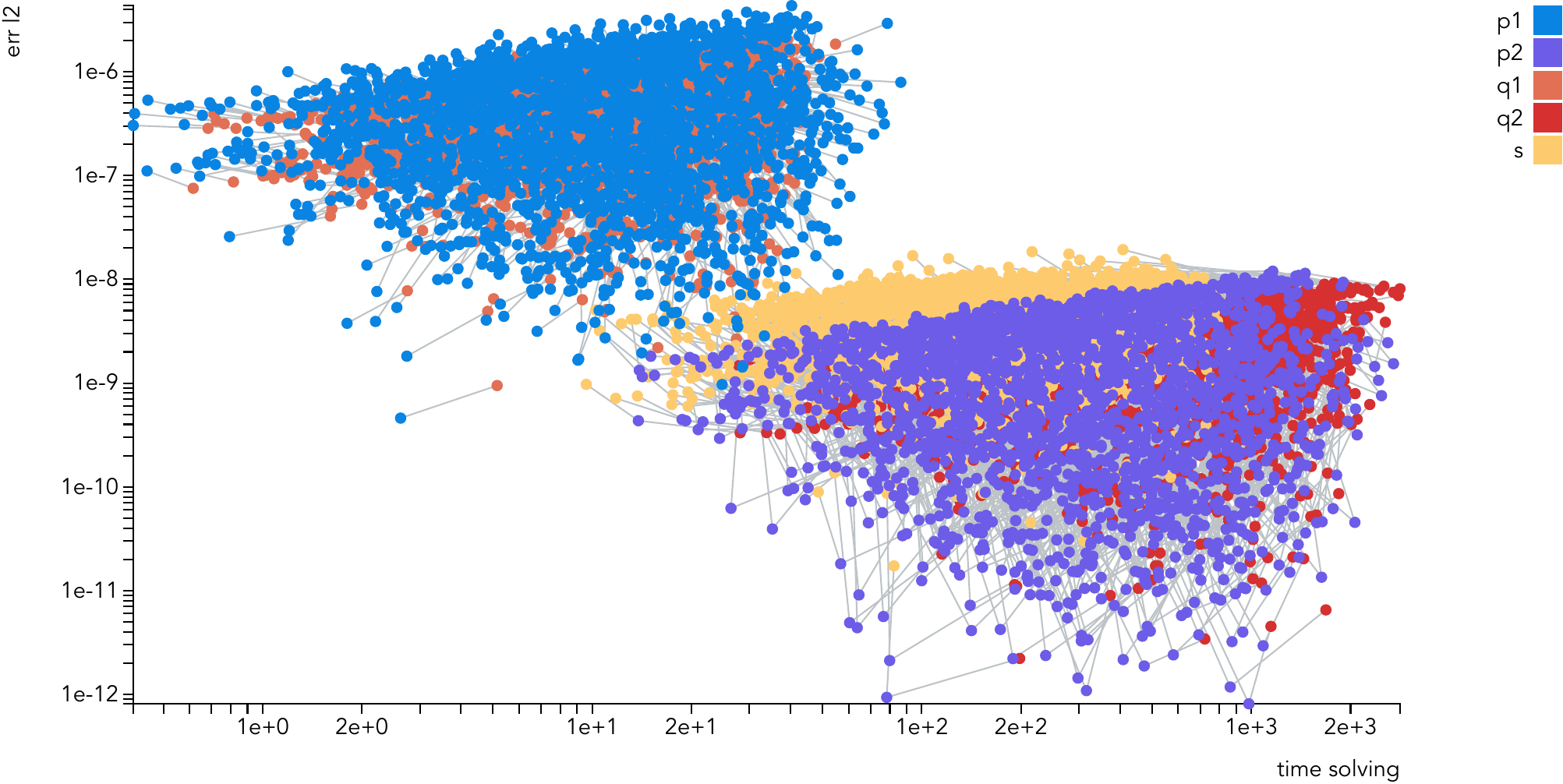}\par
        linear elasticity~\eqref{eq:lin-elast}
        %%%%%%%%%%%%%%%%%%%%%%%%%%%%%%
    }
    \caption{Total time (left) and solve time (right) vs. the $L_2$ error.}
    \label{fig:thingi10k_time}
\end{figure}

\paragraph{Hexalab}
To avoid clutter in the plots we omit the results obtained from meshes with inverted elements leading to plots with 229 points. For the complete statistics see the \interactiveplot{}.
We first compare the error of the method with respect to the average edge length, Figure~\ref{fig:hexalb_err_vs_h}. We confirm the results of Section~\ref{sec:classical-problems} for the state-of-the-art hexahedral meshing methods; the accuracy of the solution on a hexahedral and tetrahedral mesh is comparable, in our experiments, for both Poisson and linear elasticity.
Figure~\ref{fig:hexalb_time} shows the total and solve time required to reach a certain error, where we draw the same conclusions: the results of the four discretizations are similar. The plots show, as expected, that for a given mesh serendipity elements are faster but less accurate than $Q_2$ elements. However this advantage is not consistent enough to change the conclusion related to quadratic tetrahedral elements. Statistics for the individual hexahedral meshing method are available in the \interactiveplot{}.

\paragraph{Thingi10k}
We repeated the same experiment on $3\,200$ hexahedral meshes generated with MeshGems. For this large dataset it is interesting to note that qualitative behavior of  the edge length vs. error curve (Figure~\ref{fig:thingi10k_err_vs_h}) is different between hexahedra and tetrahedra: the curve for tetrahedral elements exhibit the expected convergence, while the curve for hexahedra is more flat. This effect comes from the fact that MeshGems is an octree-based method with a tendency to create highly anisotropic meshes. This effect can be mitigated by limiting the difference between the minimal and maximal refinement levels in the octree used to construct the mesh. This leads to more uniform element size, and the results become similar to the results for tetrahedral meshes and the Hexalab dataset, Figure~\ref{fig:uniform}.% Another interesting curiosity is that the Poisson equation is almost not affected by the presence of inverted elements (top) while elasticity is (bottom).\todo{because?}\DP{vote for dropping}

We also compared running and solve times (Figure~\ref{fig:thingi10k_time}) and, as expected, serendipity elements are faster than $Q_2$ elements but have a larger error. Tetrahedral elements are between the two hexahedral elements: their accuracy is similar to $Q_2$ and a running/solving time similar to serendipity.

\section{Discussion and conclusions}
\label{sec:conclusions}

\begin{figure}\centering\footnotesize
    \includegraphics[width=0.7\linewidth]{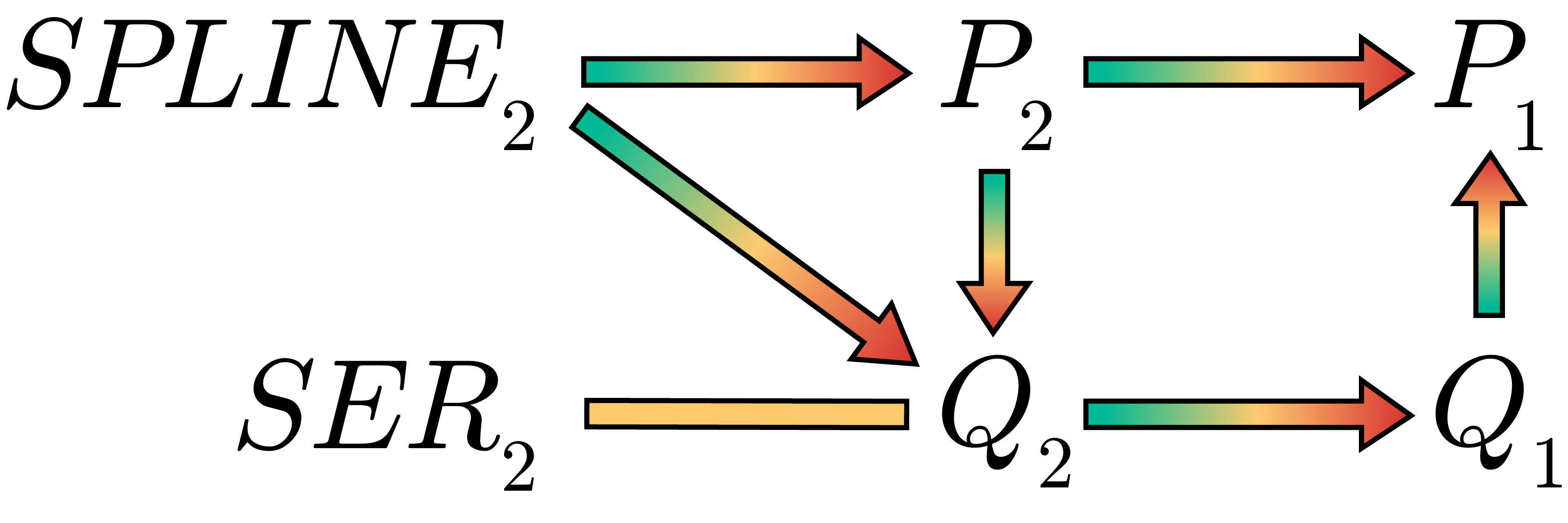}
    \caption{The arrows indicate which method is inferior (red side); the yellow box indicates that the methods are comparable.}
    \label{fig:diagram_elements}
\end{figure}

We presented a large-scale, quantitative study of several common types of finite elements applied to five elliptic PDEs. Our results are consistent on all elliptic PDEs we tried.

We summarize our findings in Figure~\ref{fig:diagram_elements}, which allows us to draw the following conclusions for the five elliptic PDEs we considered in our study:
\begin{enumerate}
    \item Consistently with well-known observations, $P_1$ elements are less efficient (more time spent to obtain a solution with given accuracy) than all other options in all our experiments (Sections~\ref{sec:classical-problems} and~\ref{sec:large-dataset}).
    \item $Q_2$ elements are slightly more accurate than quadratic serendipity elements $\mathrm{SER}_2$ but are slightly more expensive for a fixed mesh (Section~\ref{sec:large-dataset}).
    \item $P_2$ elements are generally more efficient than  $P_1$, $Q_1$, $Q_2$, $\mathrm{SER}_2$, that is, we can obtain a given target error in less time, if we can chose the mesh resolution optimal for the desired error level. We were not able to identify any disadvantages for these elements for the range of problems and geometries  we have considered (Sections~\ref{sec:classical-problems} and~\ref{sec:large-dataset}).
    \item Quadratic spline elements $\mathrm{SPLINE}_2$ (on a regular lattice) are more efficient than  $Q_2$ elements (Section \ref{sec:torsion}). $\mathrm{SPLINE}_2$ are also more efficient compared to $P_2$ (3x faster solving time for the same accuracy) but with a much longer assembly time (12 times slower, which could be reduced with more advanced integration techniques \cite{Schillinger:2014:RBE}). Their use, however, is restricted by the current meshing technology, as they require meshes with regular grid structure almost everywhere for optimal performance.  When these elements are mixed with standard $Q_2$ elements to handle general hexahedral meshes with singular vertices and edges ~\cite{Schneider:2018:PF,Wei:2018:BBS}, their performance advantage is considerably reduced.
\end{enumerate}

For the five elliptic PDEs we considered, unstructured tetrahedral meshes with quadratic Lagrangian basis are a good choice for a ``black-box'' analysis pipeline: robust tetrahedral meshing algorithms that can process thousands of real-world models exist \cite{Hu:2018:TMI}, and $p$-refinement can be used to compensate for the rare badly shaped triangles introduced by the meshing algorithms \cite{Schneider:2018:DSA}.

We leave the extension of this study to non-elliptic PDEs, multiphysics, and collision response as future work. Another important potential extension is the study of bases with orders higher than 2, as is typically the case in IGA setting, or an extension to spectral elements.
Another venue for future work is to analyze the impact of the existing different per-element optimizations (e.g., reduced quadrature, hourglass control, special quadrature rules that exploit the tensor-product structure of $Q$ elements, etc.). However, we note that these different optimizations will mostly impact the performance of the assembly and will have little influence on the solve time, which dominates the runtime for sufficiently large problems.

\revisionn{Finally, we acknowledge that the choice of elements is only one of the sources of error in numerical simulations; in realistic scenarios, material models, boundary conditions, or domain shape also play role in the accuracy of a simulation. Extending our study to account for these sources of error is an interesting venue for future work.}

% \section*{Acknowledgments}
\begin{acks}
% missing IIS-1320635, DMS-1436591
We thank the NYU IT High Performance Computing for resources, services, and staff expertise.
This work was partially supported by the NSF CAREER award under Grant No. 1652515, the NSF grants OAC-1835712, OIA-1937043, CHS-1908767, CHS-1901091, DMS-1436591, NSERC DGECR-2021-00461 and RGPIN-2021-03707, the SNSF grant P2TIP2\_175859, a Sloan Fellowship, a gift from Adobe Research, a gift from nTopology, and a gift from Advanced Micro Devices, Inc.
\end{acks}

%%
%% The next two lines define the bibliography style to be used, and
%% the bibliography file.
% \bibliographystyle{ACM-Reference-Format}
% \bibliography{biblio,software}
%%% -*-BibTeX-*-
%%% Do NOT edit. File created by BibTeX with style
%%% ACM-Reference-Format-Journals [18-Jan-2012].

%%
%% If your work has an appendix, this is the place to put it.
% \appendix

\end{document}